\documentclass[12pt]{article}

\usepackage{amsfonts,amsmath,amsthm,amssymb,amscd}
\usepackage[dvips]{graphicx}
\usepackage{mathrsfs}

\theoremstyle{definition}
\newtheorem{theorem}{Theorem}
\newtheorem{lemma}[theorem]{Lemma}
\newtheorem{proposition}[theorem]{Proposition}
\newtheorem{corollary}[theorem]{Corollary}

\numberwithin{equation}{section}
\numberwithin{theorem}{section}

\begin{document}

\begin{center}
{\bf{\Large Cubic theta functions and modular forms of level six }}
\end{center}

\begin{center}
By Kazuhide Matsuda
\end{center}

\begin{center}
Faculty of Fundamental Science, National Institute of Technology, Niihama College,\\
7-1 Yagumo-chou, Niihama, Ehime, Japan, 792-8580. \\
E-mail: matsuda@sci.niihama-nct.ac.jp  \\
Fax: 81-0897-37-7809 
\end{center}

\noindent
{\bf Abstract}
The aim of the research presented in this paper is to derive the systems of ordinary differential equations (ODEs) satisfied by modular forms of level six and to construct extensions of the differential field of the cubic theta functions, generalizing the classical Ramanujan and Halphen fields. We treat both modular forms that appear in the literature and others that do not. We find Riccati equations satisfied by level six modular forms, explore applications to number theory, and find new relations among the Eisenstein series.

\section{Introduction}
\label{intro}
\par
Much research has been devoted to ordinary differential equations (ODEs) satisfied by modular forms
(see, for example, Pol \cite{Pol}). 
This paper is concerned with generalizing two systems of ODEs: Ramanujan's and Halphen's (Berndt \cite{Berndt}; Halphen \cite{Halphen}). 
Ramanujan's system treats modular forms of level one and Halphen's system deals with modular forms of level two. 
\par
Throughout this paper, 
the upper half plane $\mathbb{H}^2$ is defined by 
$
\mathbb{H}^2=
\{
\tau\in\mathbb{C} \,\, | \,\, \Im \tau>0
\}. 
$
Let $q=\exp(2\pi i  \tau)$. We define the {\it Dedekind eta function} by 
\begin{equation*}
\displaystyle
\eta(\tau)=q^{\frac{1}{24}} \prod_{n=1}^{\infty} (1-q^n)=q^{\frac{1}{24}} (q; q)_{\infty}. 
\end{equation*}
\par 
Ramanujan's system is given by 
\begin{equation}
\label{eqn:Ramanujan-ODE}
q\frac{d E_2}{dq}=\frac{(E_2)^2-E_4  }{12  }, \,\,
q\frac{d E_4}{dq}=\frac{E_2 E_4-E_6  }{3  }, \,\,
q\frac{d E_6}{dq}=\frac{E_2 E_6-(E_4)^2  }{2  }, 
\end{equation}
where the Eisenstein series $E_2, E_4$, and $E_6$ are defined by 
\begin{align*}
E_2(q)&=E_2(\tau)=1-24\sum_{n=1}^{\infty} \sigma_1(n) q^n, \,\,
E_4(q)=E_4(\tau)=1+240\sum_{n=1}^{\infty} \sigma_3(n) q^n, \\
E_6(q)&=E_6(\tau)=1-504\sum_{n=1}^{\infty} \sigma_5(n) q^n, \,\, \sigma_k(n)=\sum_{d|n} d^k \,\,(k=1,2,3,4,5,\ldots).
\end{align*}
We note that Ramanujan's system is equivalent to Chazy's third-order nonlinear ODE,
\begin{equation*}
y^{\prime \prime \prime}=
2y y^{\prime \prime}-3(y^{\prime})^2. 
\end{equation*}
\par
Halphen's system is given by 
\begin{equation*}
X^{\prime}+Y^{\prime}=2XY, \,
Y^{\prime}+Z^{\prime}=2YZ, \,
Z^{\prime}+X^{\prime}=2ZX, \,
\end{equation*} 
the solution of which is expressed by the logarithmic derivatives of theta constants. 
Ohyama \cite{Ohyama0} showed that 
Halphen's differential field is an extension of Ramanujan's differential field, 
For which the Galois group is the symmetric group $S_3.$
\par
Huber \cite{Huber} and the present writer \cite{Matsuda-1} have derived the systems of ODEs satisfied by the cubic theta functions,
\begin{align*}
a(q)=&\sum_{m,n\in\mathbb{Z}} q^{m^2+mn+n^2}, \,\,
b(q)=\sum_{m,n\in\mathbb{Z}} \omega^{n-m} q^{m^2+mn+n^2}, \\
c(q)=&\sum_{m,n\in\mathbb{Z}} q^{(n+\frac13)^2+(n+\frac13)(m+\frac13)+(m+\frac13)^2}, \,\,\omega=e^{\frac{2\pi i}{3}}, \,\,|q|<1.
\end{align*}
In particular, our systems of ODEs are given by 
\footnotesize{
\begin{align}
q\frac{d}{dq}a(q)=&\frac{3a^3(q)+a(q)E_2(q)-4b^3(q)}{12},   \,\,\,
q\frac{d}{dq}E_2(q)=\frac{-9a^4(q)+8a(q)b^3(q)+(E_2(q))^2}{12}, \label{eqn:Matsda-ODE-1}\\
q\frac{d}{dq}b^3(q)=&\frac{ -a^2(q)b^3(q)+E_2(q) b^3(q) }{4   },  \notag
\end{align}
}
\normalsize{
and
}
\footnotesize{
\begin{align}
q\frac{d}{dq}a(q)=&\frac{-3a^3(q)+a(q)E_2(q^3)+4c^3(q)}{12}, \,\,\,
q\frac{d}{dq}E_2(q^3)=\frac{-9a^4(q)+8a(q)c^3(q)+9(E_2(q^3))^2}{36}, \label{eqn:Matsda-ODE-2}\\
q\frac{d}{dq}c^3(q)=&\frac{ a^2(q)c^3(q)+3E_2(q^3) c^3(q) }{4   }. \notag
\end{align}
}
\par
\normalsize
{
The aim of the research presented in this paper is to derive the systems of ODEs satisfied by modular forms of level six 
and 
construct extensions of the differential field of the cubic theta functions, 
$
\mathbb{C}(a(q), E_2(q), b^3(q) ). 
$
}
\par
This paper is organized as follows: 
Section \ref{sec:notations} fixes our notation. 
Section \ref{sec:properties} reviews Farkas and Kra's theory of theta functions with rational characteristics. 
Section \ref{sec:Weierstrass} treats Weierstrass' elliptic function theory. 
Section \ref{sec:derivative-level-six} describes derivative formulas of level six. 
\par
Sections \ref{sec:proof:level-six-E_2(q)}, \ref{sec:proof:level-six-E_2(q^6)}, and \ref{sec:proof:level-six-E_2(q)-split} derive the systems of ODEs satisfied by modular forms of level six. 
Our main theorems are Theorems \ref{thm-level6-E_2(q)-(1,2/3)}, \ref{thm-level6-E_2(q^6)-(2/3,1)}, and 
\ref{thm-level6-E_2(q)-(0,1/3)-(0,2/3)}. 
Sections \ref{sec:proof:level-six-E_2(q)} and \ref{sec:proof:level-six-E_2(q)-split} show applications to number theory 
by evaluating $a^2(q)a(q), \, a(q)a(q^2)a(q^4)$, and similar expressions. 
\par
Sections \ref{sec:proof:level-six-E_2(q)} and \ref{sec:proof:level-six-E_2(q^6)} treat modular forms of level six that also appear in Cooper \cite[pp. 375]{Cooper}: 
 
\begin{equation*}
\left(
\frac{\eta(2\tau)  \eta^6(3\tau)}
{
\eta^2(\tau) \eta^3(6\tau)
}, \,
\frac{\eta^6(\tau)  \eta(6\tau)}
{
\eta^3(2\tau) \eta^2(3\tau)
}
\right), \,\,
\mathrm{and} \,\,
\left(
\frac{\eta(3\tau)  \eta^6(2\tau)}
{
\eta^3(\tau) \eta^2(6\tau)
}, \,
\frac{\eta(\tau)  \eta^6(6\tau)}
{
\eta^2(2\tau) \eta^3(3\tau)
}
\right).
\end{equation*}
These differential fields are the same and can be expressed by 
$
\mathbb{C} \left( a(q), a(q^2), E_2(q) \right); 
$ 
see Theorem \ref{thm:cubic-ODE-(l,m,n)}. 
\par
On the other hand, 
Section \ref{sec:proof:level-six-E_2(q)-split} deals with the modular forms
\begin{equation*}
\left(
\frac{\eta(\tau/2)  \eta^6(3\tau)}
{
\eta^2(\tau) \eta^3(3\tau/2)
}, \,
\frac{\eta(\tau)  \eta^3(3\tau/2) \eta^3(6\tau)  }
{
\eta(\tau/2) \eta(2\tau) \eta^3(3\tau) 
}
\right),
\end{equation*}
which are not treated by Cooper \cite{Cooper}. 
The differential field of Theorem \ref{thm-level6-E_2(q)-(0,1/3)-(0,2/3)} can be also expressed by 
$
\mathbb{C} \left( a(q^{\frac12}), a(q), a(q^2), E_2(q) \right),
$ 
which is a two-dimensional extension of 
$
\mathbb{C} \left( a(q), a(q^2), E_2(q) \right); 
$ 
see Theorem \ref{thm:cubic-ODE-(k, l,m,n)}. 
\par
The differential field $\mathbb{C} \left( a(q^{\frac12}), a(q), a(q^2), E_2(q) \right)$ contains Halphen's differential field.
(For the details of Halphen's differential field, see Matsuda \cite{Matsuda-3}.) 
This paper derives additional relations among the Eisenstein series, 
\begin{equation*}
A=E_2(q^{\frac12}), B=E_2(q), C=E_2(q^2), D=E_2(q^{\frac32}), E=E_2(q^3), F=E_2(q^6). 
\end{equation*}
\par
Section \ref{sec:Riccati-level-6} derives Riccati equations satisfied by modular functions of level six: 
\begin{equation*}
f(\tau)=
\frac
{
\eta^4(2\tau)
\eta^8(3\tau)
}
{
\eta^8(\tau)
\eta^4(6\tau)
},    \,\, \,\,
g(\tau)
=
\frac
{
\eta^4(\tau) \eta^8(6\tau)
}
{
\eta^8(2\tau) \eta^4(3\tau)
}
\end{equation*}
and 
\begin{equation*}
t(\tau)=g(\tau/2)=
\frac
{
\eta^4(\tau/2) \eta^8(3\tau)
}
{
\eta^8(\tau)  \eta^4(3\tau/2)
}, \,\,\,\,
u(\tau)
=
\frac
{
\eta^4(\tau) \eta^4(3\tau/2) \eta^4(6\tau)  
}
{
\eta^4(\tau/2) \eta^4(2\tau)  \eta^4(3\tau)
}. 
\end{equation*}

\subsection*{Acknowledgments}
This work was supported by JSPS KAKENHI Grant Number JP17K14213.

\section{Notations}
\label{sec:notations}
For the positive integers $j,k,$ and $n\in\mathbb{N},$ let
$d_{j,k}(n)$ denote the number of positive divisors $d$ of $n$ such that $d\equiv j \,\,\mathrm{mod} \,k,$. 
Let 
$d_{j,k}^{*}(n)$ denote the number of positive divisors $d$ of $n$ such that $d\equiv j \,\,\mathrm{mod} \,k$ and $n/d$ is odd.  
Then 
$
d_{j,k}^{*}(n)=d_{j,k}(n)-d_{j,k}(n/2). 
$ 
Moreover, 
for $k$ and $n\in\mathbb{N},$ if 
$\sigma_k(n)$ is the sum of the $k$-th powers of the positive divisors of $n,$ 
then
$d_{j,k}(n)=\sigma_k(n)=0$ for $n\in\mathbb{Q}\setminus\mathbb{N}_0,$ 
where $ \mathbb{N}_0=\{0,1,2,3,\ldots,\}.$ 
\par
For each positive integer, $n\in\mathbb{N},$ define the Legendre symbol by 
\begin{equation*}
\left(\frac{n}{3}\right)
=
\begin{cases}
+1, &\text{if $ n \equiv 1 \,\,\mathrm{mod} \,3,$}  \\
-1, &\text{if $n\equiv -1 \,\,\mathrm{mod} \,3,$}  \\
0, &\text{if $n\equiv 0 \,\,\mathrm{mod} \,3.$}
\end{cases}
\end{equation*}
\par
The cubic theta functions can be also expressed by 
\begin{align*}
a(q)=&1+6\sum_{n=1}^{\infty} (d_{1,3}(n)-d_{2,3}(n) ) q^n, \\
b^3(q)=&
\frac
{
(q;q)_{\infty}^9
}
{
(q^3;q^3)_{\infty}^3
}
=
1-9
\sum_{n=1}^{\infty}q^n \left( \sum_{d|n} d^2 \left( \frac{d}{3} \right)  \right),  \,\,
c^3(q)=
27
q\frac{(q^3;q^3)_{\infty}^9}{(q;q)_{\infty}^3}
=27\sum_{n=1}^{\infty} q^n \left( \sum_{d|n} d^2 \left( \frac{n/d}{3} \right)  \right). 
\end{align*}
\par
Moreover, we introduce Ramanujan's theta functions:
\begin{equation*}
\varphi(q)=\sum_{n\in\mathbb{Z}} q^{n^2} \,\,\mathrm{and} \,\,  \psi(q)=\sum_{n=0}^{\infty} q^{\frac{n(n+1)}{2}} \,\,
\mathrm{with} \,\,|q|<1. 
\end{equation*}

\section{Properties of the theta functions}
\label{sec:properties}

\subsection{Definitions}
Following the work of Farkas and Kra \cite{Farkas-Kra}, 
we introduce the {\it theta function with characteristics,} 
which is defined by 
\begin{align*}
\theta 
\left[
\begin{array}{c}
\epsilon \\
\epsilon^{\prime}
\end{array}
\right] (\zeta, \tau) 
=
\theta 
\left[
\begin{array}{c}
\epsilon \\
\epsilon^{\prime}
\end{array}
\right] (\zeta) 
:=&\sum_{n\in\mathbb{Z}} \exp
\left(2\pi i\left[ \frac12\left(n+\frac{\epsilon}{2}\right)^2 \tau+\left(n+\frac{\epsilon}{2}\right)\left(\zeta+\frac{\epsilon^{\prime}}{2}\right) \right] \right), 
\end{align*}
where $\epsilon, \epsilon^{\prime}\in\mathbb{R}, \, \zeta\in\mathbb{C},$ and $\tau\in\mathbb{H}^{2}.$ 
The {\it theta constants} are given by 
\begin{equation*}
\theta 
\left[
\begin{array}{c}
\epsilon \\
\epsilon^{\prime}
\end{array}
\right]
:=
\theta 
\left[
\begin{array}{c}
\epsilon \\
\epsilon^{\prime}
\end{array}
\right] (0, \tau).
\end{equation*}
In terms of Jacobi’s original $\vartheta$ notation for theta functions,
\begin{equation}
\label{eqn:theta-null}
\vartheta_2=
\theta 
\left[
\begin{array}{c}
1 \\
0
\end{array}
\right], \,\,
\vartheta_3=
\theta 
\left[
\begin{array}{c}
0 \\
0
\end{array}
\right], \,\,
\vartheta_4=
\theta 
\left[
\begin{array}{c}
0 \\
1
\end{array}
\right]. 
\end{equation}
Furthermore, 
we denote the derivative coefficients of the theta function by 
\begin{equation*}
\theta^{\prime} 
\left[
\begin{array}{c}
\epsilon \\
\epsilon^{\prime}
\end{array}
\right]
:=\left.
\frac{\partial}{\partial \zeta} 
\theta 
\left[
\begin{array}{c}
\epsilon \\
\epsilon^{\prime}
\end{array}
\right] (\zeta, \tau)
\right|_{\zeta=0}, 
\,
\theta^{\prime \prime} 
\left[
\begin{array}{c}
\epsilon \\
\epsilon^{\prime}
\end{array}
\right]
:=\left.
\frac{\partial^2 }{\partial \zeta^2} 
\theta 
\left[
\begin{array}{c}
\epsilon \\
\epsilon^{\prime}
\end{array}
\right] (\zeta, \tau)
\right|_{\zeta=0},
\end{equation*}
and
\begin{equation*}
\theta^{\prime \prime \prime} 
\left[
\begin{array}{c}
\epsilon \\
\epsilon^{\prime}
\end{array}
\right]
:=\left.
\frac{\partial^3 }{\partial \zeta^3} 
\theta 
\left[
\begin{array}{c}
\epsilon \\
\epsilon^{\prime}
\end{array}
\right] (\zeta, \tau)
\right|_{\zeta=0},\,\,
\theta^{(n)} 
\left[
\begin{array}{c}
\epsilon \\
\epsilon^{\prime}
\end{array}
\right]
:=\left.
\frac{\partial^n}{\partial \zeta^n} 
\theta 
\left[
\begin{array}{c}
\epsilon \\
\epsilon^{\prime}
\end{array}
\right] (\zeta, \tau)
\right|_{\zeta=0}, \,\,\,(n=1,2,3,4,\ldots). 
\end{equation*}
In particular, Jacobi's derivative formula is given by 
\begin{equation}
\label{eqn:Jacobi-derivative}
\theta^{\prime} 
\left[
\begin{array}{c}
1 \\
1
\end{array}
\right] 
=
-\pi 
\theta
\left[
\begin{array}{c}
0 \\
0
\end{array}
\right] 
\theta
\left[
\begin{array}{c}
1 \\
0
\end{array}
\right] 
\theta
\left[
\begin{array}{c}
0 \\
1
\end{array}
\right].  
\end{equation}

\subsection{Basic properties}
We first note that 
for $m,n\in\mathbb{Z},$ 
\begin{equation}
\label{eqn:integer-char}
\theta 
\left[
\begin{array}{c}
\epsilon \\
\epsilon^{\prime}
\end{array}
\right] (\zeta+n+m\tau, \tau) =
\exp(2\pi i)\left[\frac{n\epsilon-m\epsilon^{\prime}}{2}-m\zeta-\frac{m^2\tau}{2}\right]
\theta 
\left[
\begin{array}{c}
\epsilon \\
\epsilon^{\prime}
\end{array}
\right] (\zeta,\tau),
\end{equation}
and 
\begin{equation}
\theta 
\left[
\begin{array}{c}
\epsilon +2m\\
\epsilon^{\prime}+2n
\end{array}
\right] 
(\zeta,\tau)
=\exp(\pi i \epsilon n)
\theta 
\left[
\begin{array}{c}
\epsilon \\
\epsilon^{\prime}
\end{array}
\right] 
(\zeta,\tau).
\end{equation}
Furthermore, 
it is easy to see that 
\begin{equation*}
\theta 
\left[
\begin{array}{c}
-\epsilon \\
-\epsilon^{\prime}
\end{array}
\right] (\zeta,\tau)
=
\theta 
\left[
\begin{array}{c}
\epsilon \\
\epsilon^{\prime}
\end{array}
\right] (-\zeta,\tau)
\,\,
\mathrm{and}
\,\,
\theta^{\prime} 
\left[
\begin{array}{c}
-\epsilon \\
-\epsilon^{\prime}
\end{array}
\right] (\zeta,\tau)
=
-
\theta^{\prime} 
\left[
\begin{array}{c}
\epsilon \\
\epsilon^{\prime}
\end{array}
\right] (-\zeta,\tau).
\end{equation*}
\par
For $m,n\in\mathbb{R},$ 
we see that 
\begin{align}
\label{eqn:real-char}
&\theta 
\left[
\begin{array}{c}
\epsilon \\
\epsilon^{\prime}
\end{array}
\right] \left(\zeta+\frac{n+m\tau}{2}, \tau\right)   \notag\\
&=
\exp(2\pi i)\left[
-\frac{m\zeta}{2}-\frac{m^2\tau}{8}-\frac{m(\epsilon^{\prime}+n)}{4}
\right]
\theta 
\left[
\begin{array}{c}
\epsilon+m \\
\epsilon^{\prime}+n
\end{array}
\right] 
(\zeta,\tau). 
\end{align}
We note that 
$\theta 
\left[
\begin{array}{c}
\epsilon \\
\epsilon^{\prime}
\end{array}
\right] \left(\zeta, \tau\right)$ has only one zero in the fundamental parallelogram, 
which is given by 
$$
\zeta=\frac{1-\epsilon}{2}\tau+\frac{1-\epsilon^{\prime}}{2}. 
$$

\subsection{Jacobi's triple product identity}
All the theta functions have infinite product expansions, which are given by 
\begin{align}
\theta 
\left[
\begin{array}{c}
\epsilon \\
\epsilon^{\prime}
\end{array}
\right] (\zeta, \tau) &=\exp\left(\frac{\pi i \epsilon \epsilon^{\prime}}{2}\right) x^{\frac{\epsilon^2}{4}} z^{\frac{\epsilon}{2}}    \notag  \\
                           &\quad 
                           \displaystyle \times\prod_{n=1}^{\infty}(1-x^{2n})(1+e^{\pi i \epsilon^{\prime}} x^{2n-1+\epsilon} z)(1+e^{-\pi i \epsilon^{\prime}} x^{2n-1-\epsilon}/z),  \label{eqn:Jacobi-triple}
\end{align}
where $x=\exp(\pi i \tau)$ and $z=\exp(2\pi i \zeta).$ 
Therefore, it follows from Jacobi's derivative formula (\ref{eqn:Jacobi-derivative}) that 
\begin{equation*}
\label{eqn:Jacobi}
\theta^{\prime} 
\left[
\begin{array}{c}
1 \\
1
\end{array}
\right](0,\tau) 
=
-2\pi 
q^{\frac18}
\prod_{n=1}^{\infty}(1-q^n)^3, \,\,q=\exp(2\pi i \tau). 
\end{equation*}

\subsection{Spaces of $N$-th order $\theta$-functions}

Following Farkas and Kra \cite{Farkas-Kra}, 
we define 
$\mathcal{F}_{N}\left[
\begin{array}{c}
\epsilon \\
\epsilon^{\prime}
\end{array}
\right] $ to be the set of entire functions $f$ that satisfy the two functional equations, 
$$
f(\zeta+1)=\exp(\pi i \epsilon) \,\,f(\zeta),
$$
and 
$$
f(\zeta+\tau)=\exp(-\pi i)[\epsilon^{\prime}+2N\zeta+N\tau] \,\,f(\zeta), \quad \zeta\in\mathbb{C},  \,\,\tau \in\mathbb{H}^2,
$$ 
where 
$N$ is a positive integer and 
$\left[
\begin{array}{c}
\epsilon \\
\epsilon^{\prime}
\end{array}
\right] \in\mathbb{R}^2.$ 
This set of functions is called the space of {\it $N$-th order $\theta$-functions with characteristic }
$\left[
\begin{array}{c}
\epsilon \\
\epsilon^{\prime}
\end{array}
\right]. $ 
Note that 
$$
\dim \mathcal{F}_{N}\left[
\begin{array}{c}
\epsilon \\
\epsilon^{\prime}
\end{array}
\right] =N;
$$
for the proof, see Farkas and Kra \cite[pp.133]{Farkas-Kra}.

\subsection{The heat equation}
The theta function satisfies the following heat equation: 
\begin{equation}
\label{eqn:heat}
\frac{\partial^2}{\partial \zeta^2}
\theta
\left[
\begin{array}{c}
\epsilon \\
\epsilon^{\prime}
\end{array}
\right](\zeta,\tau)
=
4\pi i
\frac{\partial}{\partial \tau}
\theta
\left[
\begin{array}{c}
\epsilon \\
\epsilon^{\prime}
\end{array}
\right](\zeta,\tau). 
\end{equation}

\begin{lemma}
\label{lem:first-order-deri-square}
{\it
For every $(z,\tau)\in\mathbb{C}\times\mathbb{H}^2,$ 
we have 
\begin{equation*}
\left\{
\frac{
\theta^{\prime}
\left[
\begin{array}{c}
\epsilon \\
\epsilon^{\prime}
\end{array}
\right](z)
}
{
\theta
\left[
\begin{array}{c}
\epsilon \\
\epsilon^{\prime}
\end{array}
\right](z)
}
\right\}^2
=
4 \pi i 
\frac{d}{d \tau}
\log
\theta
\left[
\begin{array}{c}
\epsilon \\
\epsilon^{\prime}
\end{array}
\right](z)
-
\frac
{
d^2
}
{
dz^2
}
\log
\theta
\left[
\begin{array}{c}
\epsilon \\
\epsilon^{\prime}
\end{array}
\right](z). 
\end{equation*}
}
\end{lemma}

\begin{proof}
The lemma can be proved by direct calculation. 
\end{proof}

\section{Weierstrass elliptic function theory  }
\label{sec:Weierstrass}

We introduce Weierstrass $\wp$-function and $\sigma$-function:
\begin{align*}
\wp(z;\omega_1, \omega_2)=\wp(z)
=&
\frac{1}{z^2}+\sum_
{
\tiny{
\begin{matrix}
(m,n)\in\mathbb{Z}^2 \\ 
(m,n)\neq(0,0)
\end{matrix}
}
}
\left(
\frac{1}{(z-m\omega_1-n\omega_2)^2}-\frac{1}{(m\omega_1+n\omega_2)^2}
\right),   \\
\sigma(z;\omega_1, \omega_2)=\sigma(z)
=&
z
\prod_
{
\tiny{
\begin{matrix}
(m,n)\in\mathbb{Z}^2 \\ 
(m,n)\neq(0,0)
\end{matrix}
}
}\left( 1-\frac{z}{m\omega_1+n\omega_2}  \right)
\exp
\left(
\frac{z}{m\omega_1+n\omega_2}
+
\frac{z^2}{2(m\omega_1+n\omega_2)^2}
\right),
\end{align*}
where $z,\omega_1,\omega_2$ are complex numbers with $\omega_2/ \omega_1 \not\in\mathbb{R}.$
\par
From Whittaker and Watson \cite[pp. 437, 459]{WW}, we recall the following formulas:
\begin{align*}
\wp^{\prime}(z)^2=&4\wp(z)^3-g_2\wp(z)-g_3,  \\
g_2(\omega_1, \omega_2)=&
60G_4(\omega_1, \omega_2)
=
60
\sum_
{
\tiny{
\begin{matrix}
(m,n)\in\mathbb{Z}^2 \\ 
(m,n)\neq(0,0)
\end{matrix}
}
}
\frac{1}{(m\omega_1+n\omega_2)^4},  \\
g_3(\omega_1, \omega_2)=&
140G_6(\omega_1, \omega_2)
=
140
\sum_
{
\tiny{
\begin{matrix}
(m,n)\in\mathbb{Z}^2 \\ 
(m,n)\neq(0,0)
\end{matrix}
}
}
\frac{1}{(m\omega_1+n\omega_2)^6}. 
\end{align*}
These imply that
\begin{equation}
\wp^{\prime \prime}(z)=6\wp^2(z)-\frac12g_2. 
\end{equation}
In particular, we note that 
\begin{equation}
\label{eqn:g_2-g_3-E_4-E_6}
g_2(1,\tau)=\frac{4\pi^4}{3} E_4(\tau)=\frac{4\pi^4}{3} E_4(q), \,\,
g_3(1,\tau)=\frac{8\pi^6}{27}E_6(\tau)=\frac{8\pi^6}{27}E_6(q),
\end{equation}
and
\begin{equation}
\label{eqn:Delta-J}
q(q;q)_{\infty}^{24}=\frac{(E_4(q))^3-(E_6(q))^2}{1728}, \,\,
j(q)=\frac{1728  (E_4(q))^3  }{(E_4(q))^3-(E_6(q))^2   }.
\end{equation}
Moreover, from Farkas and Kra \cite[pp. 124]{Farkas-Kra}, we recall 
\begin{align*}
\wp(z;1, \tau)=&
\frac{1}{3}
\frac
{
\theta^{\prime  \prime \prime} 
\left[
\begin{array}{c}
1 \\
1
\end{array}
\right]
}
{
\theta^{\prime } 
\left[
\begin{array}{c}
1 \\
1
\end{array}
\right]
}-
\frac{d^2}{dz^2} \log 
\theta
\left[
\begin{array}{c}
1 \\
1
\end{array}
\right](z,\tau),  \\
\sigma(z;\omega_1, \omega_2)=&
\exp
\left(
\frac{\eta_1 z^2}{2\omega_1}
\right)
\frac{\omega_1}
{ 
\theta^{\prime}
\left[
\begin{array}{c}
1 \\
1
\end{array}
\right] }
\theta
\left[
\begin{array}{c}
1 \\
1
\end{array}
\right]
\left(
\frac{z}{\omega_1},
\tau
\right),
\end{align*}
where $\omega_2/\omega_1=\tau\in\mathbb{H}^2.$

\section{Derivative formulas of level six   }
\label{sec:derivative-level-six}

\begin{theorem}(The Residue Theorem for Elliptic Functions)
\label{thm-fundamental-elliptic-function}
{\it
The sum of all the residues of an elliptic function in the fundamental parallelogram is zero. 
}
\end{theorem}

\begin{proposition}
\label{prop:1st-derivative-(1,1/3)-(1/3,1)}
{\it
For every $\tau\in\mathbb{H}^2,$ we have 
\begin{equation}
\label{eqn:1st-deri-(1,1/3)-(1,2/3)}
2
\frac
{
\theta^{\prime }
\left[
\begin{array}{c}
1 \\
\frac13
\end{array}
\right]
}
{
\theta
\left[
\begin{array}{c}
1 \\
\frac13
\end{array}
\right]
}
-
\frac
{
\theta^{\prime}
\left[
\begin{array}{c}
1 \\
\frac23
\end{array}
\right]
}
{
\theta
\left[
\begin{array}{c}
1 \\
\frac23
\end{array}
\right]
}
+
\frac
{
\theta^{\prime}
\left[
\begin{array}{c}
1 \\
1
\end{array}
\right]
\theta
\left[
\begin{array}{c}
1 \\
\frac23
\end{array}
\right]
}
{
\theta
\left[
\begin{array}{c}
1 \\
0
\end{array}
\right]
\theta
\left[
\begin{array}{c}
1 \\
\frac13
\end{array}
\right]
}
=0,
\end{equation}
and
\begin{equation}
\label{eqn:1st-deri-(1/3,1)-(2/3,1)}
2
\frac
{
\theta^{\prime }
\left[
\begin{array}{c}
\frac13 \\
1
\end{array}
\right]
}
{
\theta
\left[
\begin{array}{c}
\frac13  \\
1
\end{array}
\right]
}
-
\frac
{
\theta^{\prime}
\left[
\begin{array}{c}
\frac23\\
1
\end{array}
\right]
}
{
\theta
\left[
\begin{array}{c}
\frac23 \\
1
\end{array}
\right]
}
+
\zeta_3^2
\frac
{
\theta^{\prime}
\left[
\begin{array}{c}
1 \\
1
\end{array}
\right]
\theta
\left[
\begin{array}{c}
\frac23 \\
1
\end{array}
\right]
}
{
\theta
\left[
\begin{array}{c}
0 \\
1
\end{array}
\right]
\theta
\left[
\begin{array}{c}
\frac13\\
1
\end{array}
\right]
}
=0, \,\,\zeta_3=\exp(2 \pi i/3).
\end{equation}
}
\end{proposition}

\begin{proof}
Consider the following elliptic functions:
\begin{equation*}
\varphi(z)
=
\frac
{
\theta^{2}
\left[
\begin{array}{c}
1 \\
\frac13
\end{array}
\right](z)
\theta
\left[
\begin{array}{c}
1 \\
-\frac23
\end{array}
\right](z)
}
{
\theta^{2}
\left[
\begin{array}{c}
1 \\
1
\end{array}
\right](z)
\theta
\left[
\begin{array}{c}
1 \\
0
\end{array}
\right](z)
}, \,\,
\text{and}
\,\,
\psi(z)
=
\frac
{
\theta^{2}
\left[
\begin{array}{c}
\frac13 \\
1
\end{array}
\right](z)
\theta
\left[
\begin{array}{c}
-\frac23 \\
1
\end{array}
\right](z)
}
{
\theta^{2}
\left[
\begin{array}{c}
1 \\
1
\end{array}
\right](z)
\theta
\left[
\begin{array}{c}
0 \\
1
\end{array}
\right](z)
}. 
\end{equation*}
\ par
In the fundamental parallelogram, 
the poles of $\varphi(z)$ are $z=0$ and $z=1/2.$ 
Direct calculation yields 
\begin{equation*}
\mathrm{Res}(\varphi(z),0)
=
\frac
{
\theta^2
\left[
\begin{array}{c}
1 \\
\frac13
\end{array}
\right]
\theta
\left[
\begin{array}{c}
1 \\
\frac23
\end{array}
\right]
}
{
\theta^{\prime}
\left[
\begin{array}{c}
1 \\
1
\end{array}
\right]^2
\theta
\left[
\begin{array}{c}
1 \\
0
\end{array}
\right]
}
\left\{
2
\frac
{
\theta^{\prime }
\left[
\begin{array}{c}
1 \\
\frac13
\end{array}
\right]
}
{
\theta
\left[
\begin{array}{c}
1 \\
\frac13
\end{array}
\right]
}
-
\frac
{
\theta^{\prime}
\left[
\begin{array}{c}
1 \\
\frac23
\end{array}
\right]
}
{
\theta
\left[
\begin{array}{c}
1 \\
\frac23
\end{array}
\right]
}
\right\},
\end{equation*}
and
\begin{equation*}
\mathrm{Res}\left( \varphi(z) ,\frac12 \right)
=
\frac
{
\theta^2
\left[
\begin{array}{c}
1 \\
\frac23
\end{array}
\right]
\theta
\left[
\begin{array}{c}
1 \\
\frac13
\end{array}
\right]
}
{
\theta^{\prime}
\left[
\begin{array}{c}
1 \\
1
\end{array}
\right]
\theta^2
\left[
\begin{array}{c}
1 \\
0
\end{array}
\right]
}.
\end{equation*}
Equation (\ref{eqn:1st-deri-(1,1/3)-(1,2/3)}) follows from the residue theorem. Equation (\ref{eqn:1st-deri-(1/3,1)-(2/3,1)}) can be proved by using $\psi(z)$ in the same way. 
 
\end{proof}

\begin{corollary}
{\it
For $q\in\mathbb{C}$ with $|q|<1,$ we have 
\begin{align*}
\frac{(q ;q)_{\infty}^6 (q^6;q^6)_{\infty}  }{(q^2;q^2)_{\infty}^3 (q^3;q^3)_{\infty}^2  } =&1-6\sum_{n=1}^{\infty} (d_{1,6}(n)-d_{5,6}(n)) q^n+18 \sum_{n=1}^{\infty} ( d_{2,6}(n)-d_{4,6}(n) )q^n,  \\
q \frac{(q ; q )_{\infty}  (q^6; q^6)_{\infty}^6   }{  (q^2; q^2)_{\infty}^2 (q^3; q^3)_{\infty}^3   }=&   \sum_{n=1}^{\infty}(d_{1,6}(n)-d_{5,6}(n))q^n-2  \sum_{n=1}^{\infty}  (d_{1,3}(n)-d_{2,3}(n)) q^{2n}.
\end{align*}
}
\end{corollary}

\begin{proposition}
\label{prop:1st-derivative-(1,1/3)-(0,1/3)-(0,2/3)}
{\it
For every $\tau\in\mathbb{H}^2,$ we have 
\begin{equation}
\label{eqn:1st-deri-(1,1/3)-(0,1/3)}
2
\frac
{
\theta^{\prime }
\left[
\begin{array}{c}
1 \\
\frac13
\end{array}
\right]
}
{
\theta
\left[
\begin{array}{c}
1 \\
\frac13
\end{array}
\right]
}
+
\frac
{
\theta^{\prime}
\left[
\begin{array}{c}
0 \\
\frac13
\end{array}
\right]
}
{
\theta
\left[
\begin{array}{c}
0 \\
\frac13
\end{array}
\right]
}
-
\frac
{
\theta^{\prime}
\left[
\begin{array}{c}
1 \\
1
\end{array}
\right]
\theta
\left[
\begin{array}{c}
0 \\
\frac13
\end{array}
\right]
}
{
\theta
\left[
\begin{array}{c}
0 \\
1
\end{array}
\right]
\theta
\left[
\begin{array}{c}
1 \\
\frac13
\end{array}
\right]
}
=0,
\end{equation}
and
\begin{equation}
\label{eqn:1st-deri-(1,1/3)-(0,2/3)}
2
\frac
{
\theta^{\prime }
\left[
\begin{array}{c}
1 \\
\frac13
\end{array}
\right]
}
{
\theta
\left[
\begin{array}{c}
1  \\
\frac13
\end{array}
\right]
}
-
\frac
{
\theta^{\prime}
\left[
\begin{array}{c}
0\\
\frac23
\end{array}
\right]
}
{
\theta
\left[
\begin{array}{c}
0 \\
\frac23
\end{array}
\right]
}
-
\frac
{
\theta^{\prime}
\left[
\begin{array}{c}
1 \\
1
\end{array}
\right]
\theta
\left[
\begin{array}{c}
0 \\
\frac23
\end{array}
\right]
}
{
\theta
\left[
\begin{array}{c}
0 \\
0
\end{array}
\right]
\theta
\left[
\begin{array}{c}
1\\
\frac13
\end{array}
\right]
}
=0. 
\end{equation}
}
\end{proposition}

\begin{proof}
Consider the following elliptic functions:
\begin{equation*}
\varphi(z)
=
\frac
{
\theta^{2}
\left[
\begin{array}{c}
1 \\
\frac13
\end{array}
\right](z)
\theta
\left[
\begin{array}{c}
0 \\
\frac13
\end{array}
\right](z)
}
{
\theta^{2}
\left[
\begin{array}{c}
1 \\
1
\end{array}
\right](z)
\theta
\left[
\begin{array}{c}
0 \\
1
\end{array}
\right](z)
}, \,\,
\text{and}
\,\,
\psi(z)
=
\frac
{
\theta^{2}
\left[
\begin{array}{c}
1 \\
\frac13
\end{array}
\right](z)
\theta
\left[
\begin{array}{c}
0 \\
-\frac23
\end{array}
\right](z)
}
{
\theta^{2}
\left[
\begin{array}{c}
1 \\
1
\end{array}
\right](z)
\theta
\left[
\begin{array}{c}
0 \\
0
\end{array}
\right](z)
}. 
\end{equation*}
The proposition can be proved in the same way as Proposition \ref{prop:1st-derivative-(1,1/3)-(1/3,1)}. 
\end{proof}

\begin{corollary}
{\it
For $q\in\mathbb{C}$ with $|q|<1,$ we have 
\begin{align*}
\frac{(q ;q)_{\infty}^6 (q^{\frac32};q^{\frac32})_{\infty}  }{(q^{\frac12};q^{\frac12})_{\infty}^3 (q^3;q^3)_{\infty}^2  } =&
1+3\sum_{n=1}^{\infty} (d_{1,3}^{*}(n)-d_{2,3}^{*}(n)) q^{\frac{n}{2}}+6 \sum_{n=1}^{\infty} ( d_{1,3}(n)-d_{2,3}(n) )q^n,  \\
\frac{(q^{\frac12} ; q^{\frac12} )_{\infty}^3  (q^2; q^2)_{\infty}^3   (q^3; q^3)_{\infty}   }{  (q; q)_{\infty}^3 (q^{\frac32}; q^{\frac32})_{\infty}  (q^6; q^6)_{\infty}   }=& 
1-3  \sum_{n=1}^{\infty}(d_{1,6}^{*}(n)+d_{2,6}^{*}(n)  -d_{4,6}^{*}(n)-d_{5,6}^{*}(n))q^{\frac{n}{2}}    \\
&\hspace{50mm}+6  \sum_{n=1}^{\infty}  (d_{1,3}(n)-d_{2,3}(n)) q^{n}.
\end{align*}
}
\end{corollary}

\begin{theorem}
\label{thm:analogue-Jacobi-(1,2/3)-(2/3,1)}
{\it
For every $\tau\in\mathbb{H}^2,$ we have
\begin{equation*}
\theta^{\prime}
\left[
\begin{array}{c}
1 \\
\frac23
\end{array}
\right]
=
\frac13
\frac
{
\theta^{\prime}
\left[
\begin{array}{c}
1 \\
1
\end{array}
\right]
\theta^3
\left[
\begin{array}{c}
1 \\
\frac13
\end{array}
\right]
}
{
\theta^2
\left[
\begin{array}{c}
1 \\
\frac23
\end{array}
\right]
\theta
\left[
\begin{array}{c}
1 \\
0
\end{array}
\right]
}, \,\,
\theta^{\prime}
\left[
\begin{array}{c}
\frac23 \\
1
\end{array}
\right]
=
\frac13
\frac
{
\theta^{\prime}
\left[
\begin{array}{c}
1 \\
1
\end{array}
\right]
\theta^3
\left[
\begin{array}{c}
\frac13 \\
1
\end{array}
\right]
}
{
\theta^2
\left[
\begin{array}{c}
\frac23 \\
1
\end{array}
\right]
\theta
\left[
\begin{array}{c}
0 \\
1
\end{array}
\right]
}, 
\end{equation*}
and
\begin{equation*}
\theta^{\prime}
\left[
\begin{array}{c}
0 \\
\frac13
\end{array}
\right]
=
\frac13
\frac
{
\theta^{\prime}
\left[
\begin{array}{c}
1 \\
1
\end{array}
\right]
\theta^3
\left[
\begin{array}{c}
1 \\
\frac13
\end{array}
\right]
}
{
\theta^2
\left[
\begin{array}{c}
0 \\
\frac13
\end{array}
\right]
\theta
\left[
\begin{array}{c}
0 \\
1
\end{array}
\right]
}, 
\,\,
\theta^{\prime}
\left[
\begin{array}{c}
0 \\
\frac23
\end{array}
\right]
=
\frac13
\frac
{
\theta^{\prime}
\left[
\begin{array}{c}
1 \\
1
\end{array}
\right]
\theta^3
\left[
\begin{array}{c}
1 \\
\frac13
\end{array}
\right]
}
{
\theta^2
\left[
\begin{array}{c}
0 \\
\frac23
\end{array}
\right]
\theta
\left[
\begin{array}{c}
0 \\
0
\end{array}
\right]
}.
\end{equation*}
}
\end{theorem}

\begin{proof}
Consider the following elliptic function:
\begin{equation*}
\frac
{
\theta^3
\left[
\begin{array}{c}
1 \\
1
\end{array}
\right](z)
}
{
\theta^2
\left[
\begin{array}{c}
1 \\
\frac23
\end{array}
\right](z)
\theta
\left[
\begin{array}{c}
1 \\
-\frac13
\end{array}
\right](z)
}, \,\,
\frac
{
\theta^3
\left[
\begin{array}{c}
1 \\
1
\end{array}
\right](z)
}
{
\theta^2
\left[
\begin{array}{c}
\frac23 \\
1
\end{array}
\right](z)
\theta
\left[
\begin{array}{c}
-\frac13 \\
1
\end{array}
\right](z)
}, 
\end{equation*}
and
\begin{equation*}
\frac
{
\theta^3
\left[
\begin{array}{c}
1 \\
1
\end{array}
\right](z)
}
{
\theta^2
\left[
\begin{array}{c}
0 \\
\frac13
\end{array}
\right](z)
\theta
\left[
\begin{array}{c}
1 \\
\frac13
\end{array}
\right](z)
}, 
\,\,
\frac
{
\theta^3
\left[
\begin{array}{c}
1 \\
1
\end{array}
\right](z)
}
{
\theta^2
\left[
\begin{array}{c}
0 \\
\frac23
\end{array}
\right](z)
\theta
\left[
\begin{array}{c}
1 \\
-\frac13
\end{array}
\right](z)
}.
\end{equation*}
The theorem can be proved in the same way as Proposition \ref{prop:1st-derivative-(1,1/3)-(1/3,1)}. 
\end{proof}

\begin{corollary}
{\it
For $q\in\mathbb{C}$ with $|q|<1,$ we have 
\begin{align*}
\frac{   (q^2; q^2)_{\infty}   (q^3; q^3)_{\infty}^6  }{   (q; q)_{\infty}^2 (q^6; q^6)_{\infty}^3  }  
=&1+2\sum_{n=1}^{\infty} (d_{1,6}(n)+d_{2,6}(n)-d_{4,6}(n)-d_{5,6}(n)) q^n,  \\
\frac{ (q^3; q^3)_{\infty} (q^2; q^2)_{\infty}^6  }{ (q; q)_{\infty}^3 (q^6; q^6)_{\infty}^2 }
=&1+ 3 \sum_{n=1}^{\infty} (d_{1,6}(n)-d_{5,6}(n)) q^n,   \\
q^{\frac12} \frac{ (q^{\frac12}; q^{\frac12})_{\infty}  (q^3;q^3)_{\infty}^6 }{  (q;q)_{\infty}^2  (q^{\frac32}; q^{\frac32})_{\infty}^3 }=&\sum_{n=1}^{\infty}(d_{1,3}^{*}(n)-d_{2,3}^{*}(n)) q^{\frac{n}{2}},  \\
q^{\frac12} \frac{   (q;q)_{\infty}   (q^{\frac32}; q^{\frac32})_{\infty}^3    (q^6;q^6)_{\infty}^3    }{   (q^{\frac12}; q^{\frac12})_{\infty}   (q^2;q^2)_{\infty}    (q^3;q^3)_{\infty}^3 }
=&\sum_{n=1}^{\infty}(d_{1,6}^{*}(n)+d_{2,6}^{*}(n)    - d_{4,6}^{*}(n)-d_{5,6}^{*}(n)  ) q^{\frac{n}{2}}.
\end{align*}
}
\end{corollary}

\begin{proof}
The product-series identities follow from the formulas 
\begin{equation*}
\theta^{\prime}
\left[
\begin{array}{c}
1 \\
\frac23
\end{array}
\right]
/
\theta
\left[
\begin{array}{c}
1 \\
\frac23
\end{array}
\right], \,\,
\theta^{\prime}
\left[
\begin{array}{c}
\frac23 \\
1
\end{array}
\right]
/
\theta
\left[
\begin{array}{c}
\frac23 \\
1
\end{array}
\right], \,\,
\theta^{\prime}
\left[
\begin{array}{c}
0 \\
\frac13
\end{array}
\right]
/
\theta
\left[
\begin{array}{c}
0 \\
\frac13
\end{array}
\right], \,\,\text{and} \,\,
\theta^{\prime}
\left[
\begin{array}{c}
0 \\
\frac23
\end{array}
\right]
/
\theta
\left[
\begin{array}{c}
0 \\
\frac23
\end{array}
\right]. 
\end{equation*}
\end{proof}

\begin{corollary}
\label{coro:analogue-Jacobi-(1,1/3)-(1/3,1)}
{\it
For every $\tau\in\mathbb{H}^2,$ we have 
\begin{equation*}
\frac{
\theta^{\prime}
\left[
\begin{array}{c}
1 \\
\frac13
\end{array}
\right]
}
{
\theta
\left[
\begin{array}{c}
1 \\
\frac13
\end{array}
\right]
}
=
\frac16
\frac{
\theta^{\prime}
\left[
\begin{array}{c}
1 \\
1
\end{array}
\right]
\theta^3
\left[
\begin{array}{c}
1 \\
\frac13
\end{array}
\right]
}
{
\theta
\left[
\begin{array}{c}
1 \\
0
\end{array}
\right]
\theta^3
\left[
\begin{array}{c}
1 \\
\frac23
\end{array}
\right]
}
-
\frac12
\frac{
\theta^{\prime}
\left[
\begin{array}{c}
1 \\
1
\end{array}
\right]
\theta
\left[
\begin{array}{c}
1 \\
\frac23
\end{array}
\right]
}
{
\theta
\left[
\begin{array}{c}
1 \\
0
\end{array}
\right]
\theta
\left[
\begin{array}{c}
1 \\
\frac13
\end{array}
\right]
},
\end{equation*}
and
\begin{equation*}
\frac{
\theta^{\prime}
\left[
\begin{array}{c}
\frac13 \\
1
\end{array}
\right]
}
{
\theta
\left[
\begin{array}{c}
\frac13 \\
1
\end{array}
\right]
}
=
\frac16
\frac{
\theta^{\prime}
\left[
\begin{array}{c}
1 \\
1
\end{array}
\right]
\theta^3
\left[
\begin{array}{c}
\frac13 \\
1
\end{array}
\right]
}
{
\theta
\left[
\begin{array}{c}
0 \\
1
\end{array}
\right]
\theta^3
\left[
\begin{array}{c}
\frac23 \\
1
\end{array}
\right]
}
-
\frac{ \zeta_3^2 }{2}
\frac{
\theta^{\prime}
\left[
\begin{array}{c}
1 \\
1
\end{array}
\right]
\theta
\left[
\begin{array}{c}
\frac23 \\
1
\end{array}
\right]
}
{
\theta
\left[
\begin{array}{c}
0 \\
1
\end{array}
\right]
\theta
\left[
\begin{array}{c}
\frac13 \\
1
\end{array}
\right]
}.
\end{equation*}
}
\end{corollary}

\begin{corollary}
\label{relation-f-g}
{\it
For $q\in\mathbb{C}$ with $|q|<1,$ 
set
\begin{equation*}
f(q)=\frac
{
(q^2; q^2)_{\infty}^4 (q^3; q^3)_{\infty}^8
}
{
(q;q)_{\infty}^8
(q^6; q^6)_{\infty}^4
}, \,\,
{\it and} \,\,
g(q)
=
q\frac{
(q;q)_{\infty}^4 (q^6; q^6)_{\infty}^8
}
{
(q^2; q^2)_{\infty}^8 (q^3; q^3)_{\infty}^4
}.
\end{equation*}
Then, we have 
\begin{equation*}
f=\frac{g-1}{9g-1}, \,\,{\it and} \,\,
g=\frac{f-1}{9f-1}.
\end{equation*}
}
\end{corollary}

\begin{proof}
We first set 
\begin{align*}
P_1(q)=&1+2\sum_{n=1}^{\infty} (d_{1,6}(n)+d_{2,6}(n)-d_{4,6}(n)-d_{5,6}(n)) q^n, \\
Q_1(q)=&1-6\sum_{n=1}^{\infty} (d_{1,6}(n)-d_{5,6}(n)) q^n+18 \sum_{n=1}^{\infty} ( d_{2,6}(n)-d_{4,6}(n) )q^n, \\
\end{align*}
and
\begin{align*}
P_2(q)=&1+ 3 \sum_{n=1}^{\infty} (d_{1,6}(n)-d_{5,6}(n)) q^n,  \\
Q_2(q)=& \sum_{n=1}^{\infty}(d_{1,6}(n)-d_{5,6}(n))q^n-2  \sum_{n=1}^{\infty}  (d_{1,3}(n)-d_{2,3}(n)) q^{2n},  \\
\end{align*}
which imply that 
\begin{equation}
\label{relation-(P_1,Q_1)-(P_2,Q_2)}
P_1=P_2-Q_2, \,\,  Q_1=P_2-9 Q_2, \,\, \mathrm{and} \,\, P_2=\frac{9 P_1-Q_1}{8}, \,\,Q_2=\frac{P_1-Q_1}{8}.
\end{equation}
Then, it follows that 
\begin{equation*}
f(q)=\frac{P_1}{Q_1}=\frac{P_2-Q_2}{P_2-9 Q_2}=\frac{g(q)-1  }{ 9  g(q)-1 }, 
\,\,\mathrm{and}\,\,
g(q)=\frac{Q_2}{P_2}=\frac{P_1-Q_1}{9P_1-Q_1}=\frac{f(q)-1}{9f(q)-1}.
\end{equation*}
\end{proof}

\begin{corollary}
\label{coro:analogue-Jacobi-(1,1/3)-(0,1/3)-(0,2/3)}
{\it
For every $\tau\in\mathbb{H}^2,$ we have 
\begin{equation*}
\frac{
\theta^{\prime}
\left[
\begin{array}{c}
1 \\
\frac13
\end{array}
\right]
}
{
\theta
\left[
\begin{array}{c}
1 \\
\frac13
\end{array}
\right]
}
=
-
\frac16
\frac{
\theta^{\prime}
\left[
\begin{array}{c}
1 \\
1
\end{array}
\right]
\theta^3
\left[
\begin{array}{c}
1 \\
\frac13
\end{array}
\right]
}
{
\theta
\left[
\begin{array}{c}
0 \\
1
\end{array}
\right]
\theta^3
\left[
\begin{array}{c}
0 \\
\frac13
\end{array}
\right]
}
+
\frac12
\frac{
\theta^{\prime}
\left[
\begin{array}{c}
1 \\
1
\end{array}
\right]
\theta
\left[
\begin{array}{c}
0 \\
\frac13
\end{array}
\right]
}
{
\theta
\left[
\begin{array}{c}
0 \\
1
\end{array}
\right]
\theta
\left[
\begin{array}{c}
1 \\
\frac13
\end{array}
\right]
},
\end{equation*}
and
\begin{equation*}
\frac{
\theta^{\prime}
\left[
\begin{array}{c}
1 \\
\frac13
\end{array}
\right]
}
{
\theta
\left[
\begin{array}{c}
1 \\
\frac13
\end{array}
\right]
}
=
\frac16
\frac{
\theta^{\prime}
\left[
\begin{array}{c}
1 \\
1
\end{array}
\right]
\theta^3
\left[
\begin{array}{c}
1 \\
\frac13
\end{array}
\right]
}
{
\theta
\left[
\begin{array}{c}
0 \\
0
\end{array}
\right]
\theta^3
\left[
\begin{array}{c}
0 \\
\frac23
\end{array}
\right]
}
+
\frac{1 }{2}
\frac{
\theta^{\prime}
\left[
\begin{array}{c}
1 \\
1
\end{array}
\right]
\theta
\left[
\begin{array}{c}
0 \\
\frac23
\end{array}
\right]
}
{
\theta
\left[
\begin{array}{c}
0 \\
0
\end{array}
\right]
\theta
\left[
\begin{array}{c}
1 \\
\frac13
\end{array}
\right]
}.
\end{equation*}
}
\end{corollary}

\section{ODEs satisfied by modular forms of level six (1)}
\label{sec:proof:level-six-E_2(q)}
The main result of this section is as follows:

\begin{theorem}
\label{thm-level6-E_2(q)-(1,2/3)}
{\it
For $q\in\mathbb{C}$ with $|q|<1,$ 
set 
\begin{align*}
P(q)=&\frac{   (q^2; q^2)_{\infty}   (q^3; q^3)_{\infty}^6  }{   (q; q)_{\infty}^2 (q^6; q^6)_{\infty}^3  }  
=1+2\sum_{n=1}^{\infty} (d_{1,6}(n)+d_{2,6}(n)-d_{4,6}(n)-d_{5,6}(n)) q^n,   \\
Q(q)=&\frac{(q ;q)_{\infty}^6 (q^6;q^6)_{\infty}  }{(q^2;q^2)_{\infty}^3 (q^3;q^3)_{\infty}^2  } =1-6\sum_{n=1}^{\infty} (d_{1,6}(n)-d_{5,6}(n)) q^n+18 \sum_{n=1}^{\infty} ( d_{2,6}(n)-d_{4,6}(n) )q^n,  \\
R(q)=&E_2(q)=1-24\sum_{n=1}^{\infty} \sigma_1(n) q^n. 
\end{align*}
Then, we have 
\begin{align*}
\label{eqn-level6-E_2(q)}
q\frac{d}{dq}P=&\frac{27 P^3-36P^2 Q+5P Q^2+4 P R}{48},  \,\,
q\frac{d}{dq}Q=\frac{ -27 P^2 Q+24 P Q^2-Q^3+4 Q R }{48}, \\
q\frac{d}{dq}R=&\frac{ -729 P^4+972 P^3 Q-270 P^2 Q^2+12 P Q^3-Q^4+16 R^2}{192}. 
\end{align*}
}
\end{theorem}

\begin{proof}
The proof consists of Subsections 6.1-6.5.
\end{proof}

\subsubsection*{Remark}
Throughout this section, for $q\in\mathbb{C}$ with $|q|<1,$ 
we set 
\begin{align*}
P(q)=&\frac{   (q^2; q^2)_{\infty}   (q^3; q^3)_{\infty}^6  }{   (q; q)_{\infty}^2 (q^6; q^6)_{\infty}^3  }  
=1+2\sum_{n=1}^{\infty} (d_{1,6}(n)+d_{2,6}(n)-d_{4,6}(n)-d_{5,6}(n)) q^n,   \\
Q(q)=&\frac{(q ;q)_{\infty}^6 (q^6;q^6)_{\infty}  }{(q^2;q^2)_{\infty}^3 (q^3;q^3)_{\infty}^2  } =1-6\sum_{n=1}^{\infty} (d_{1,6}(n)-d_{5,6}(n)) q^n+18 \sum_{n=1}^{\infty} ( d_{2,6}(n)-d_{4,6}(n) )q^n,  \\
R(q)=&E_2(q)=1-24\sum_{n=1}^{\infty} \sigma_1(n) q^n. 
\end{align*}

\subsection{Notations}
For every $\tau\in\mathbb{H}^2,$ set $q=\exp(2\pi i \tau)$ and 
\begin{equation*}
X=
\frac
{
\theta^{\prime}
\left[
\begin{array}{c}
1 \\
\frac23
\end{array}
\right]
}
{
\theta
\left[
\begin{array}{c}
1 \\
\frac23
\end{array}
\right]
},  \,\,
Y=
\frac
{
\theta^{\prime}
\left[
\begin{array}{c}
1 \\
1
\end{array}
\right]
\theta
\left[
\begin{array}{c}
1 \\
\frac23
\end{array}
\right]
}
{
\theta
\left[
\begin{array}{c}
1 \\
0
\end{array}
\right]
\theta
\left[
\begin{array}{c}
1 \\
\frac13
\end{array}
\right]
},  \,\,
Z=
\frac
{
\theta^{\prime  \prime  \prime}
\left[
\begin{array}{c}
1 \\
1
\end{array}
\right]
}
{
\theta^{\prime  }
\left[
\begin{array}{c}
1 \\
1
\end{array}
\right]
}, \,\,
W=
\frac
{
\theta^{\prime}
\left[
\begin{array}{c}
1 \\
\frac13
\end{array}
\right]
}
{
\theta
\left[
\begin{array}{c}
1 \\
\frac13
\end{array}
\right]
}
=\frac{X-Y}{2}.
\end{equation*}
In particular, we note that 
\begin{equation*}
X=-\sqrt{3} \pi P(q), 
\,\,
Y=
-\frac{
\sqrt{3}
}
{
3
}
\pi Q(q), \,\,
Z=
-\pi^2 R(q).
\end{equation*}

\subsection{Relations between theta derivatives}

\begin{proposition}
\label{prop:2nd-derivative-(1,1/3)-(1,2/3)-(1,0)}
{\it
For every $\tau\in\mathbb{H}^2,$ we have 
\begin{equation}
\label{eqn:2nd-deri-(1,1/3)}
3
\frac
{
\theta^{\prime \prime}
\left[
\begin{array}{c}
1 \\
\frac13
\end{array}
\right]
}
{
\theta
\left[
\begin{array}{c}
1 \\
\frac13
\end{array}
\right]
}
-
\frac
{
\theta^{\prime \prime \prime}
\left[
\begin{array}{c}
1 \\
1
\end{array}
\right]
}
{
\theta^{\prime}
\left[
\begin{array}{c}
1 \\
1
\end{array}
\right]
}
+
6
\left\{
\frac
{
\theta^{\prime}
\left[
\begin{array}{c}
1 \\
\frac13
\end{array}
\right]
}
{
\theta
\left[
\begin{array}{c}
1 \\
\frac13
\end{array}
\right]
}
\right\}^2=0,  
\end{equation}
\begin{equation}
\label{eqn:2nd-deri-(1,2/3)}
2
\frac
{
\theta^{\prime \prime}
\left[
\begin{array}{c}
1 \\
\frac23
\end{array}
\right]
}
{
\theta
\left[
\begin{array}{c}
1 \\
\frac23
\end{array}
\right]
}
+
\frac
{
\theta^{\prime \prime}
\left[
\begin{array}{c}
1 \\
\frac13
\end{array}
\right]
}
{
\theta
\left[
\begin{array}{c}
1 \\
\frac13
\end{array}
\right]
}
-
\frac
{
\theta^{\prime  \prime \prime}
\left[
\begin{array}{c}
1 \\
1
\end{array}
\right]
}
{
\theta^{\prime}
\left[
\begin{array}{c}
1 \\
1
\end{array}
\right]
}
-4
\frac
{
\theta^{\prime}
\left[
\begin{array}{c}
1 \\
\frac13
\end{array}
\right]
}
{
\theta
\left[
\begin{array}{c}
1 \\
\frac13
\end{array}
\right]
}
\frac
{
\theta^{\prime}
\left[
\begin{array}{c}
1 \\
\frac23
\end{array}
\right]
}
{
\theta
\left[
\begin{array}{c}
1 \\
\frac23
\end{array}
\right]
}
+
2
\left\{
\frac
{
\theta^{\prime}
\left[
\begin{array}{c}
1 \\
\frac23
\end{array}
\right]
}
{
\theta
\left[
\begin{array}{c}
1 \\
\frac23
\end{array}
\right]
}
\right\}^2
=0,
\end{equation}
and
\begin{equation}
\label{eqn:2nd-deri-(1,0)}
\frac
{
\theta^{\prime \prime}
\left[
\begin{array}{c}
1 \\
\frac13
\end{array}
\right]
}
{
\theta
\left[
\begin{array}{c}
1 \\
\frac13
\end{array}
\right]
}
+
\frac
{
\theta^{\prime \prime}
\left[
\begin{array}{c}
1 \\
\frac23
\end{array}
\right]
}
{
\theta
\left[
\begin{array}{c}
1 \\
\frac23
\end{array}
\right]
}
+
\frac
{
\theta^{\prime \prime}
\left[
\begin{array}{c}
1 \\
0
\end{array}
\right]
}
{
\theta
\left[
\begin{array}{c}
1 \\
0
\end{array}
\right]
}
-
\frac
{
\theta^{\prime  \prime \prime}
\left[
\begin{array}{c}
1 \\
1
\end{array}
\right]
}
{
\theta^{\prime}
\left[
\begin{array}{c}
1 \\
1
\end{array}
\right]
}
+2
\frac
{
\theta^{\prime}
\left[
\begin{array}{c}
1 \\
\frac13
\end{array}
\right]
}
{
\theta
\left[
\begin{array}{c}
1 \\
\frac13
\end{array}
\right]
}
\frac
{
\theta^{\prime}
\left[
\begin{array}{c}
1 \\
\frac23
\end{array}
\right]
}
{
\theta
\left[
\begin{array}{c}
1 \\
\frac23
\end{array}
\right]
}
=0.
\end{equation}
}
\end{proposition}

\begin{proof}
Consider the following elliptic functions: 
\begin{equation*}
\frac
{
\theta^3
\left[
\begin{array}{c}
1 \\
\frac13
\end{array}
\right]
(z)
}
{
\theta^3
\left[
\begin{array}{c}
1 \\
1
\end{array}
\right]
(z)
}, \quad
\frac
{
\theta^2
\left[
\begin{array}{c}
1 \\
\frac23
\end{array}
\right]
(z)
\theta
\left[
\begin{array}{c}
1 \\
-\frac13
\end{array}
\right]
(z)
}
{
\theta^3
\left[
\begin{array}{c}
1 \\
1
\end{array}
\right]
(z)
}, \,\,
\frac
{
\theta
\left[
\begin{array}{c}
1 \\
\frac13
\end{array}
\right]
(z)
\theta
\left[
\begin{array}{c}
1 \\
\frac23
\end{array}
\right]
(z)
\theta
\left[
\begin{array}{c}
1 \\
0
\end{array}
\right]
(z)
}
{
\theta^3
\left[
\begin{array}{c}
1 \\
1
\end{array}
\right]
(z)
}. 
\end{equation*}
The proposition can be proved in the same way as Proposition \ref{prop:1st-derivative-(1,1/3)-(1/3,1)}. 
\end{proof}

\begin{proposition}
\label{prop:3rd-derivative-(1,1/3)-(1,2/3)}
{\it
For every $\tau\in\mathbb{H}^2,$ we have 
\begin{equation*}
\frac
{
\theta^{(3)}
\left[
\begin{array}{c}
1 \\
\frac23
\end{array}
\right]
}
{
\theta
\left[
\begin{array}{c}
1 \\
\frac23
\end{array}
\right]
}
+
6
\frac
{
\theta^{(2)}
\left[
\begin{array}{c}
1 \\
\frac23
\end{array}
\right]
}
{
\theta
\left[
\begin{array}{c}
1 \\
\frac23
\end{array}
\right]
}
\frac
{
\theta^{(1)}
\left[
\begin{array}{c}
1 \\
\frac23
\end{array}
\right]
}
{
\theta
\left[
\begin{array}{c}
1 \\
\frac23
\end{array}
\right]
}
+
3
\frac
{
\theta^{(1)}
\left[
\begin{array}{c}
1 \\
\frac23
\end{array}
\right]
}
{
\theta
\left[
\begin{array}{c}
1 \\
\frac23
\end{array}
\right]
}
\frac
{
\theta^{(2)}
\left[
\begin{array}{c}
1 \\
0
\end{array}
\right]
}
{
\theta
\left[
\begin{array}{c}
1 \\
0
\end{array}
\right]
}
-
4
\frac
{
\theta^{(1)}
\left[
\begin{array}{c}
1 \\
\frac23
\end{array}
\right]
}
{
\theta
\left[
\begin{array}{c}
1 \\
\frac23
\end{array}
\right]
}
\frac
{
\theta^{(3)}
\left[
\begin{array}{c}
1 \\
1
\end{array}
\right]
}
{
\theta^{(1)}
\left[
\begin{array}{c}
1 \\
1
\end{array}
\right]
}
+
2
\left(
\frac
{
\theta^{(1)}
\left[
\begin{array}{c}
1 \\
\frac23
\end{array}
\right]
}
{
\theta
\left[
\begin{array}{c}
1 \\
\frac23
\end{array}
\right]
}
\right)^3=0
\end{equation*}
and
\begin{align*}
&
\frac
{
\theta^{(3)}
\left[
\begin{array}{c}
1 \\
\frac23
\end{array}
\right]
}
{
\theta
\left[
\begin{array}{c}
1 \\
\frac23
\end{array}
\right]
}
+
\frac
{
\theta^{(3)}
\left[
\begin{array}{c}
1 \\
\frac13
\end{array}
\right]
}
{
\theta
\left[
\begin{array}{c}
1 \\
\frac13
\end{array}
\right]
}
+
3
\frac
{
\theta^{(2)}
\left[
\begin{array}{c}
1 \\
\frac23
\end{array}
\right]
}
{
\theta
\left[
\begin{array}{c}
1 \\
\frac23
\end{array}
\right]
}
\frac
{
\theta^{(1)}
\left[
\begin{array}{c}
1 \\
\frac23
\end{array}
\right]
}
{
\theta
\left[
\begin{array}{c}
1 \\
\frac23
\end{array}
\right]
}
+
6
\frac
{
\theta^{(1)}
\left[
\begin{array}{c}
1 \\
\frac13
\end{array}
\right]
}
{
\theta
\left[
\begin{array}{c}
1 \\
\frac13
\end{array}
\right]
}
\frac
{
\theta^{(2)}
\left[
\begin{array}{c}
1 \\
\frac23
\end{array}
\right]
}
{
\theta
\left[
\begin{array}{c}
1 \\
\frac23
\end{array}
\right]
}
+
6
\frac
{
\theta^{(1)}
\left[
\begin{array}{c}
1 \\
\frac13
\end{array}
\right]
}
{
\theta
\left[
\begin{array}{c}
1 \\
\frac13
\end{array}
\right]
}
\left\{
\frac
{
\theta^{(1)}
\left[
\begin{array}{c}
1 \\
\frac23
\end{array}
\right]
}
{
\theta
\left[
\begin{array}{c}
1 \\
\frac23
\end{array}
\right]
}
\right\}^2    \\
&
+
6
\frac
{
\theta^{(2)}
\left[
\begin{array}{c}
1 \\
\frac13
\end{array}
\right]
}
{
\theta
\left[
\begin{array}{c}
1 \\
\frac13
\end{array}
\right]
}
\frac
{
\theta^{(1)}
\left[
\begin{array}{c}
1 \\
\frac23
\end{array}
\right]
}
{
\theta
\left[
\begin{array}{c}
1 \\
\frac23
\end{array}
\right]
}
+
6
\left\{
\frac
{
\theta^{(1)}
\left[
\begin{array}{c}
1 \\
\frac13
\end{array}
\right]
}
{
\theta
\left[
\begin{array}{c}
1 \\
\frac13
\end{array}
\right]
}
\right\}^2
\frac
{
\theta^{(1)}
\left[
\begin{array}{c}
1 \\
\frac23
\end{array}
\right]
}
{
\theta
\left[
\begin{array}{c}
1 \\
\frac23
\end{array}
\right]
}
+
3
\frac
{
\theta^{(2)}
\left[
\begin{array}{c}
1 \\
\frac13
\end{array}
\right]
}
{
\theta
\left[
\begin{array}{c}
1 \\
\frac13
\end{array}
\right]
}
\frac
{
\theta^{(1)}
\left[
\begin{array}{c}
1 \\
\frac13
\end{array}
\right]
}
{
\theta
\left[
\begin{array}{c}
1 \\
\frac13
\end{array}
\right]
}   \\
&
-
4
\frac
{
\theta^{(3)}
\left[
\begin{array}{c}
1 \\
1
\end{array}
\right]
}
{
\theta^{(1)}
\left[
\begin{array}{c}
1 \\
1
\end{array}
\right]
}
\left\{
\frac
{
\theta^{(1)}
\left[
\begin{array}{c}
1 \\
\frac13
\end{array}
\right]
}
{
\theta
\left[
\begin{array}{c}
1 \\
\frac13
\end{array}
\right]
}
+
\frac
{
\theta^{(1)}
\left[
\begin{array}{c}
1 \\
\frac23
\end{array}
\right]
}
{
\theta
\left[
\begin{array}{c}
1 \\
\frac23
\end{array}
\right]
}
\right\}=0.
\end{align*}
}
\end{proposition}

\begin{proof}
Consider the following elliptic functions:
\begin{equation*}
\frac
{
\theta^3
\left[
\begin{array}{c}
1 \\
\frac23
\end{array}
\right](z)
\theta
\left[
\begin{array}{c}
1 \\
0
\end{array}
\right](z)
}
{
\theta^4
\left[
\begin{array}{c}
1 \\
1
\end{array}
\right](z)
}, \quad
\frac
{
\theta^2
\left[
\begin{array}{c}
1 \\
\frac13
\end{array}
\right](z)
\theta^2
\left[
\begin{array}{c}
1 \\
\frac23
\end{array}
\right](z)
}
{
\theta^4
\left[
\begin{array}{c}
1 \\
1
\end{array}
\right](z)
}.
\end{equation*}
\end{proof}

\begin{proposition}
\label{prop:4th-derivative-(1,1/3)-(1,2/3)}
{\it
For every $\tau\in\mathbb{H}^2,$ we have 
\begin{align*}
&
4
\frac
{
\theta^{(4)}
\left[
\begin{array}{c}
1 \\
\frac23
\end{array}
\right]
}
{
\theta
\left[
\begin{array}{c}
1 \\
\frac23
\end{array}
\right]
}
+
\frac
{
\theta^{(4)}
\left[
\begin{array}{c}
1 \\
\frac13
\end{array}
\right]
}
{
\theta
\left[
\begin{array}{c}
1 \\
\frac13
\end{array}
\right]
}
-
\frac
{
\theta^{(5)}
\left[
\begin{array}{c}
1 \\
1
\end{array}
\right]
}
{
\theta^{(1)}
\left[
\begin{array}{c}
1 \\
1
\end{array}
\right]
}   
+
48
\frac
{
\theta^{(1)}
\left[
\begin{array}{c}
1 \\
\frac23
\end{array}
\right]
}
{
\theta
\left[
\begin{array}{c}
1 \\
\frac23
\end{array}
\right]
}
\frac
{
\theta^{(3)}
\left[
\begin{array}{c}
1 \\
\frac23
\end{array}
\right]
}
{
\theta
\left[
\begin{array}{c}
1 \\
\frac23
\end{array}
\right]
}  
+
16
\frac
{
\theta^{(1)}
\left[
\begin{array}{c}
1 \\
\frac13
\end{array}
\right]
}
{
\theta
\left[
\begin{array}{c}
1 \\
\frac13
\end{array}
\right]
}
\frac
{
\theta^{(3)}
\left[
\begin{array}{c}
1 \\
\frac23
\end{array}
\right]
}
{
\theta
\left[
\begin{array}{c}
1 \\
\frac23
\end{array}
\right]
}   \\
&+
16
\frac
{
\theta^{(1)}
\left[
\begin{array}{c}
1 \\
\frac23
\end{array}
\right]
}
{
\theta
\left[
\begin{array}{c}
1 \\
\frac23
\end{array}
\right]
}
\frac
{
\theta^{(3)}
\left[
\begin{array}{c}
1 \\
\frac13
\end{array}
\right]
}
{
\theta
\left[
\begin{array}{c}
1 \\
\frac13
\end{array}
\right]
}
+
36
\left\{
\frac
{
\theta^{(2)}
\left[
\begin{array}{c}
1 \\
\frac23
\end{array}
\right]
}
{
\theta
\left[
\begin{array}{c}
1 \\
\frac23
\end{array}
\right]
}
\right\}^2  
+144
\left\{
\frac
{
\theta^{(1)}
\left[
\begin{array}{c}
1 \\
\frac23
\end{array}
\right]
}
{
\theta
\left[
\begin{array}{c}
1 \\
\frac23
\end{array}
\right]
}
\right\}^2  
\frac
{
\theta^{(2)}
\left[
\begin{array}{c}
1 \\
\frac23
\end{array}
\right]
}
{
\theta
\left[
\begin{array}{c}
1 \\
\frac23
\end{array}
\right]
}  \\
&+144
\frac
{
\theta^{(1)}
\left[
\begin{array}{c}
1 \\
\frac13
\end{array}
\right]
}
{
\theta
\left[
\begin{array}{c}
1 \\
\frac13
\end{array}
\right]
}
\frac
{
\theta^{(1)}
\left[
\begin{array}{c}
1 \\
\frac23
\end{array}
\right]
}
{
\theta
\left[
\begin{array}{c}
1 \\
\frac23
\end{array}
\right]
}
\frac
{
\theta^{(2)}
\left[
\begin{array}{c}
1 \\
\frac23
\end{array}
\right]
}
{
\theta
\left[
\begin{array}{c}
1 \\
\frac23
\end{array}
\right]
}  
+24
\frac
{
\theta^{(2)}
\left[
\begin{array}{c}
1 \\
\frac13
\end{array}
\right]
}
{
\theta
\left[
\begin{array}{c}
1 \\
\frac13
\end{array}
\right]
}
\frac
{
\theta^{(2)}
\left[
\begin{array}{c}
1 \\
\frac23
\end{array}
\right]
}
{
\theta
\left[
\begin{array}{c}
1 \\
\frac23
\end{array}
\right]
}  
+72
\left\{
\frac
{
\theta^{(1)}
\left[
\begin{array}{c}
1 \\
\frac23
\end{array}
\right]
}
{
\theta
\left[
\begin{array}{c}
1 \\
\frac23
\end{array}
\right]
}  
\right\}^2
\frac
{
\theta^{(2)}
\left[
\begin{array}{c}
1 \\
\frac13
\end{array}
\right]
}
{
\theta
\left[
\begin{array}{c}
1 \\
\frac13
\end{array}
\right]
}  \\
&+96
\frac
{
\theta^{(1)}
\left[
\begin{array}{c}
1 \\
\frac13
\end{array}
\right]
}
{
\theta
\left[
\begin{array}{c}
1 \\
\frac13
\end{array}
\right]
}  
\left\{
\frac
{
\theta^{(1)}
\left[
\begin{array}{c}
1 \\
\frac23
\end{array}
\right]
}
{
\theta
\left[
\begin{array}{c}
1 \\
\frac23
\end{array}
\right]
}  
\right\}^3
+
24
\left\{
\frac
{
\theta^{(1)}
\left[
\begin{array}{c}
1 \\
\frac23
\end{array}
\right]
}
{
\theta
\left[
\begin{array}{c}
1 \\
\frac23
\end{array}
\right]
}  
\right\}^4  
+10
\left\{
\frac
{
\theta^{(3)}
\left[
\begin{array}{c}
1 \\
1
\end{array}
\right]
}
{
\theta^{(1)}
\left[
\begin{array}{c}
1 \\
1
\end{array}
\right]
}
\right\}^2  \\
&
-10
\frac
{
\theta^{(3)}
\left[
\begin{array}{c}
1 \\
1
\end{array}
\right]
}
{
\theta^{(1)}
\left[
\begin{array}{c}
1 \\
1
\end{array}
\right]
}
\left[
\frac
{
\theta^{(2)}
\left[
\begin{array}{c}
1 \\
\frac13
\end{array}
\right]
}
{
\theta
\left[
\begin{array}{c}
1 \\
\frac13
\end{array}
\right]
}
+
8
\frac
{
\theta^{(1)}
\left[
\begin{array}{c}
1 \\
\frac13
\end{array}
\right]
}
{
\theta
\left[
\begin{array}{c}
1 \\
\frac13
\end{array}
\right]
}
\frac
{
\theta^{(1)}
\left[
\begin{array}{c}
1 \\
\frac23
\end{array}
\right]
}
{
\theta
\left[
\begin{array}{c}
1 \\
\frac23
\end{array}
\right]
}
+4
\frac
{
\theta^{(2)}
\left[
\begin{array}{c}
1 \\
\frac23
\end{array}
\right]
}
{
\theta
\left[
\begin{array}{c}
1 \\
\frac23
\end{array}
\right]
}
+
12
\left\{
\frac
{
\theta^{(1)}
\left[
\begin{array}{c}
1 \\
\frac23
\end{array}
\right]
}
{
\theta
\left[
\begin{array}{c}
1 \\
\frac23
\end{array}
\right]
}
\right\}^2
\right]
=0.
\end{align*}
}
\end{proposition}

\begin{proof}
Consider the following elliptic functions:
$
\theta^4
\left[
\begin{array}{c}
1 \\
\frac23
\end{array}
\right]
(z)
\theta
\left[
\begin{array}{c}
1 \\
\frac13
\end{array}
\right]
(z)
/
\theta^5
\left[
\begin{array}{c}
1 \\
1
\end{array}
\right]
(z).
$ 
\end{proof}

\subsection{The ODE for $P(q)$}

\begin{proof}
By Proposition \ref{prop:2nd-derivative-(1,1/3)-(1,2/3)-(1,0)}, we have 
\begin{equation}
\label{eqn:2nd-deri-X-(1,2/3)}
\frac{
\theta^{\prime \prime}
\left[
\begin{array}{c}
1 \\
\frac23
\end{array}
\right]  
}
{ 
\theta
\left[
\begin{array}{c}
1 \\
\frac23
\end{array}
\right]
}
=
\frac14 X^2-\frac32 XY+\frac14 Y^2+\frac13 Z,  \,\,
\frac{
\theta^{\prime \prime}
\left[
\begin{array}{c}
1 \\
0
\end{array}
\right]  
}
{ 
\theta
\left[
\begin{array}{c}
1 \\
0
\end{array}
\right]
}
=
-\frac34 X^2+\frac32 XY+\frac14 Y^2+\frac13 Z.
\end{equation}
By Proposition \ref{prop:3rd-derivative-(1,1/3)-(1,2/3)}, we obtain
\begin{equation}
\label{eqn:3rd-deri-X-(1,2/3)}
\frac{
\theta^{\prime \prime \prime}
\left[
\begin{array}{c}
1 \\
\frac23
\end{array}
\right]  
}
{ 
\theta
\left[
\begin{array}{c}
1 \\
\frac23
\end{array}
\right]
}
=
-\frac54 X^3+\frac92 X^2 Y-\frac94 X Y^2+ XZ.  
\end{equation}
\par
Since 
\begin{equation*}
\frac{
\theta^{\prime \prime \prime}
\left[
\begin{array}{c}
1 \\
\frac23
\end{array}
\right]  
}
{ 
\theta
\left[
\begin{array}{c}
1 \\
\frac23
\end{array}
\right]
}
=
4 \pi i
\frac{d}{d\tau}
\left\{
\frac{
\theta^{\prime }
\left[
\begin{array}{c}
1 \\
\frac23
\end{array}
\right]  
}
{ 
\theta
\left[
\begin{array}{c}
1 \\
\frac23
\end{array}
\right]
}
\right\}
+
\frac{
\theta^{\prime }
\left[
\begin{array}{c}
1 \\
\frac23
\end{array}
\right]  
}
{ 
\theta
\left[
\begin{array}{c}
1 \\
\frac23
\end{array}
\right]
}
\cdot
\frac{
\theta^{\prime \prime }
\left[
\begin{array}{c}
1 \\
\frac23
\end{array}
\right]  
}
{ 
\theta
\left[
\begin{array}{c}
1 \\
\frac23
\end{array}
\right]
},  
\end{equation*}
it follows that 
\begin{equation*}
4\pi i X^{\prime}=
-\frac32 X^3+6 X^2 Y-\frac52 X Y^2+\frac23 XZ, 
\end{equation*}
where ${ }^{\prime}=d/d\tau. $ 
Considering $d/d\tau=2 \pi i q d/dq,$ we obtain the ODE for $P(q)$. 
\end{proof}

\subsection{The ODE for $Q(q)$}

\begin{proof}
By Proposition \ref{prop:2nd-derivative-(1,1/3)-(1,2/3)-(1,0)}, we have 
\begin{equation}
\label{eqn:2nd-deri-(1,1/3)}
\frac{
\theta^{\prime \prime}
\left[
\begin{array}{c}
1 \\
\frac13
\end{array}
\right]  
}
{ 
\theta
\left[
\begin{array}{c}
1 \\
\frac13
\end{array}
\right]
}
=
-\frac12 (X-Y)^2+\frac13 Z=   -2W^2+\frac13Z. 
\end{equation}
By Proposition \ref{prop:3rd-derivative-(1,1/3)-(1,2/3)}, we obtain
\begin{equation}
\label{eqn:3rd-deri-X-(1,1/3)}
\frac{
\theta^{\prime \prime \prime}
\left[
\begin{array}{c}
1 \\
\frac13
\end{array}
\right]  
}
{ 
\theta
\left[
\begin{array}{c}
1 \\
\frac13
\end{array}
\right]
}
=
-X^3+3X^2 Y+\frac12 XZ-\frac12 YZ=-X^3+3 X^2Y+WZ. 
\end{equation}
\par
Since 
\begin{align*}
\frac{
\theta^{\prime \prime \prime}
\left[
\begin{array}{c}
1 \\
\frac13
\end{array}
\right]  
}
{ 
\theta
\left[
\begin{array}{c}
1 \\
\frac13
\end{array}
\right]
}
=&
4 \pi i
\frac{d}{d\tau}
\left\{
\frac{
\theta^{\prime }
\left[
\begin{array}{c}
1 \\
\frac13
\end{array}
\right]  
}
{ 
\theta
\left[
\begin{array}{c}
1 \\
\frac13
\end{array}
\right]
}
\right\}
+
\frac{
\theta^{\prime }
\left[
\begin{array}{c}
1 \\
\frac13
\end{array}
\right]  
}
{ 
\theta
\left[
\begin{array}{c}
1 \\
\frac23
\end{array}
\right]
}
\cdot
\frac{
\theta^{\prime \prime }
\left[
\begin{array}{c}
1 \\
\frac13
\end{array}
\right]  
}
{ 
\theta
\left[
\begin{array}{c}
1 \\
\frac13
\end{array}
\right]
}  \\
=&
-2\pi i Y^{\prime}
+\frac12
(
-2X^3+\frac{15}{2}X^2 Y-4XY^2+\frac12 Y^3+XZ-\frac13 YZ
)\\ 
=&4 \pi i W^{\prime}-2 W^3+\frac13 ZW,
\end{align*}
it follows that 
\begin{align}
4 \pi i Y^{\prime}=&
\frac32 X^2 Y-4 X Y^2+\frac12 Y^3+\frac23 Y Z,  \notag  \\
4\pi i W^{\prime}=&-X^3+3X^2Y+2W^3+\frac23ZW,   \label{eqn:ODE-W}
\end{align}
where ${ }^{\prime}=d/d\tau. $ 
Considering $d/d\tau=2 \pi i q d/dq,$ we obtain the ODE for $Q(q)$. 
\end{proof}

\subsection{The ODE for $R(q)$}

\begin{proof}
We first note that 
\begin{align*}
\frac{
\theta^{(4)}
\left[
\begin{array}{c}
1 \\
\frac23
\end{array}
\right]  
}
{ 
\theta
\left[
\begin{array}{c}
1 \\
\frac23
\end{array}
\right]
}
=&
4 \pi i
\frac{d}{d\tau}
\left\{
\frac{
\theta^{\prime \prime }
\left[
\begin{array}{c}
1 \\
\frac23
\end{array}
\right]  
}
{ 
\theta
\left[
\begin{array}{c}
1 \\
\frac23
\end{array}
\right]
}
\right\}
+
\left\{
\frac{
\theta^{\prime \prime}
\left[
\begin{array}{c}
1 \\
\frac23
\end{array}
\right]  
}
{ 
\theta
\left[
\begin{array}{c}
1 \\
\frac23
\end{array}
\right]
}
\right\}^2,   \\
=&
\frac{4\pi i}{3}\frac{dZ}{d\tau}
-\frac{11 X^4}{16}+\frac{9 X^3 Y}{4}-\frac{9 X^2 Y^2}{8}+\frac{X^2 Z}{2}+\frac{X Y^3}{4}-3 X Y Z+\frac{5 Y^4}{16}+\frac{Y^2 Z}{2}+\frac{Z^2}{9},
\end{align*}
\begin{align*}
\frac{
\theta^{(4)}
\left[
\begin{array}{c}
1 \\
\frac13
\end{array}
\right]  
}
{ 
\theta
\left[
\begin{array}{c}
1 \\
\frac13
\end{array}
\right]
}
=&
4 \pi i
\frac{d}{d\tau}
\left\{
\frac{
\theta^{\prime \prime }
\left[
\begin{array}{c}
1 \\
\frac13
\end{array}
\right]  
}
{ 
\theta
\left[
\begin{array}{c}
1 \\
\frac13
\end{array}
\right]
}
\right\}
+
\left\{
\frac{
\theta^{\prime \prime}
\left[
\begin{array}{c}
1 \\
\frac13
\end{array}
\right]  
}
{ 
\theta
\left[
\begin{array}{c}
1 \\
\frac13
\end{array}
\right]
}
\right\}^2,   \\
=&
\frac{4\pi i}{3}\frac{dZ}{d\tau}
+
\frac{7 X^4}{4}-7 X^3 Y+\frac{9 X^2 Y^2}{2}-X^2 Z+X Y^3+2 X Y Z-\frac{Y^4}{4}-Y^2 Z+\frac{Z^2}{9},
\end{align*}
and
\begin{align*}
\frac{
\theta^{(5)}
\left[
\begin{array}{c}
1 \\
1
\end{array}
\right]  
}
{ 
\theta^{(1)}
\left[
\begin{array}{c}
1 \\
1
\end{array}
\right]
}
=&
4 \pi i
\frac{d}{d\tau}
\left\{
\frac{
\theta^{(3) }
\left[
\begin{array}{c}
1 \\
1
\end{array}
\right]  
}
{ 
\theta^{(1)}
\left[
\begin{array}{c}
1 \\
1
\end{array}
\right]
}
\right\}
+
\left\{
\frac{
\theta^{(3)}
\left[
\begin{array}{c}
1 \\
1
\end{array}
\right]  
}
{ 
\theta^{(1)}
\left[
\begin{array}{c}
1 \\
1
\end{array}
\right]
}
\right\}^2   
=
4\pi i \frac{dZ}{d\tau}+Z^2.
\end{align*}
By Proposition \ref{prop:4th-derivative-(1,1/3)-(1,2/3)},
we have 
\begin{equation*}
4\pi i Z^{\prime}=
-\frac{27 X^4}{8}+\frac{27 X^3 Y}{2}-\frac{45 X^2 Y^2}{4}+\frac{3 X Y^3}{2}-\frac{3 Y^4}{8}+\frac{2 Z^2}{3},
\end{equation*}
\begin{equation*}
\frac{
\theta^{(4)}
\left[
\begin{array}{c}
1 \\
\frac23
\end{array}
\right]  
}
{ 
\theta
\left[
\begin{array}{c}
1 \\
\frac23
\end{array}
\right]
}
=
-\frac{29 X^4}{16}+\frac{27 X^3 Y}{4}-\frac{39 X^2 Y^2}{8}+\frac{X^2 Z}{2}+\frac{3 X Y^3}{4}-3 X Y Z+\frac{3 Y^4}{16}+\frac{Y^2 Z}{2}+\frac{Z^2}{3},
\end{equation*}
and
\begin{equation*}
\frac{
\theta^{(4)}
\left[
\begin{array}{c}
1 \\
\frac13
\end{array}
\right]  
}
{ 
\theta
\left[
\begin{array}{c}
1 \\
\frac13
\end{array}
\right]
}
=
\frac{5 X^4}{8}-\frac{5 X^3 Y}{2}+\frac{3 X^2 Y^2}{4}-X^2 Z+\frac{3 X Y^3}{2}+2 X Y Z-\frac{3 Y^4}{8}-Y^2 Z+\frac{Z^2}{3},
\end{equation*}
where ${ }^{\prime}=d/d\tau.$ Considering $d/d\tau=2 \pi i q d/dq,$ we obtain the ODE for $R(q)$. 
\end{proof}

\subsection{Note on $E_4(q)$ and $E_6(q)$ $\cdots$ (1)}

\begin{theorem}
\label{thm:E_2-E_4-E_6-q-(1,2/3)}
{\it
For $q\in\mathbb{C}$ with $|q|<1,$ we have 
\begin{align*}
E_4(q)=&\frac{ 729 P^4-972 P^3 Q+270 P^2 Q^2-12 P Q^3+Q^4 }{16}
=
\frac{(3 P-Q) \left(243 P^3-243 P^2 Q+9 P Q^2-Q^3\right)}{16}, 
\\
E_6(q)=&
\frac{-19683 P^6+39366 P^5 Q-24057 P^4 Q^2+4860 P^3 Q^3-405 P^2 Q^4-18 P Q^5+Q^6}{64}  \\
=&
-\frac{\left(27 P^2-18 P Q-Q^2\right) \left(729 P^4-972 P^3 Q+270 P^2 Q^2-36 P Q^3+Q^4\right)}{64},  \\
q(q;q)^{24}_{\infty} 
=&
\frac{1}{4096} P^2 Q^6 (P-Q) (9 P-Q)^3,    \\
j(q)
=&
\frac{(3 P-Q)^3 \left(243 P^3-243 P^2 Q+9 P Q^2-Q^3\right)^3}{P^2 Q^6 (P-Q) (9 P-Q)^3} 
=
\frac{(3 f-1)^3 \left(243 f^3-243 f^2+9 f-1\right)^3}{ f^2  (f-1)  (9 f-1)^3},
\end{align*}
where 
\begin{equation*}
f(q)=P(q)/Q(q)=
\frac
{
(q^2; q^2)_{\infty}^4 (q^3; q^3)_{\infty}^8
}
{
(q;q)_{\infty}^8
(q^6; q^6)_{\infty}^4
}. 
\end{equation*}
}
\end{theorem}

\begin{proof}
For the proof, we use 
\begin{equation*}
\wp=\wp(z;1,\tau)
=
\frac13
\frac
{
\theta^{\prime  \prime  \prime}
\left[
\begin{array}{c}
1 \\
1
\end{array}
\right]
}
{
\theta^{\prime  }
\left[
\begin{array}{c}
1 \\
1
\end{array}
\right]
}
-
\frac{d^2}{dz^2}
\log
\theta
\left[
\begin{array}{c}
1 \\
1
\end{array}
\right](z). 
\end{equation*} 
\par
We first note that 
\begin{equation*}
\wp^{\prime \prime}=6\wp^2-g_2, \,\,
g_2(1,\tau)=\frac{4\pi^4}{3}E_4(q).
\end{equation*}
Substituting $z=-1/6,$ we have 
\begin{equation*}
g_2=\frac{27 X^4}{4}-27 X^3 Y+\frac{45 X^2 Y^2}{2}-3 X Y^3+\frac{3 Y^4}{4}
=
\frac{3}{4} (X-Y) \left(9 X^3-27 X^2 Y+3 X Y^2-Y^3\right).
\end{equation*}
\par
We next note that 
\begin{equation*}
(\wp^{\prime})^2=4\wp^3-g_2\wp-g_3, \,\,
g_3(1,\tau)=\frac{8\pi^6}{27}E_6(q).
\end{equation*}
Substituting $z=-1/6,$ we obtain 
\begin{align*}
g_3=&-\frac{27 X^6}{8}+\frac{81 X^5 Y}{4}-\frac{297 X^4 Y^2}{8}+\frac{45 X^3 Y^3}{2}-\frac{45 X^2 Y^4}{8}-\frac{3 X Y^5}{4}+\frac{Y^6}{8}     \\
=&
-\frac{1}{8} \left(3 X^2-6 X Y-Y^2\right) \left(9 X^4-36 X^3 Y+30 X^2 Y^2-12 X Y^3+Y^4\right).
\end{align*}
\end{proof}

\subsection{Note on $E_2(q)$}

\begin{theorem}
\label{thm:E_2(q)-level-6-(1,2/3)}
{\it
For $q\in\mathbb{C}$ with $|q|<1,$ we have 
\begin{align*}
P^2=&
\frac16 
\left(
-E_2(q)-E_2(q^3)+8 E_2(q^6)
\right),  \,\,
PQ=\frac16
\left(
E_2(q)-4 E_2(q^2)-3 E_2(q^3)+12 E_2(q^6)
\right), \\
Q^2=&\frac12
\left(
E_2(q)-8E_2(q^2)+9 E_2(q^3)
\right),
\end{align*}
which implies that 
\begin{align*}
E_2(q^2)
=&
\frac{1}{16}
\left(
27 P^2-18 PQ-Q^2+8 R
\right), \,\,
E_2(q^3)
=\frac16
\left(
9 P^2-6 P Q + Q^2 + 2 R
\right),  \\
E_2(q^6)=&
\frac{1}{48}
\left(
45 P^2-6 P Q+Q^2+8 R
\right). 
\end{align*}
}
\end{theorem}

\begin{proof}
The formula for $P^2$ can be obtained by applying Lemma \ref{lem:first-order-deri-square} to the characteristic $(\epsilon, \epsilon^{\prime})=(1,2/3).$ 
The formulas for $PQ$ and $Q^2$ follow from equations (\ref{eqn:2nd-deri-(1,1/3)}) and (\ref{eqn:2nd-deri-(1,0)}) in Proposition \ref{prop:2nd-derivative-(1,1/3)-(1,2/3)-(1,0)}. 
\end{proof}

\subsubsection*{Remark}
The formula for $P^2$ was proved in \cite{Matsuda-2}.

\subsection{Note on $E_4(q)$ and $E_6(q)$ $\cdots$ (2)}

\begin{theorem}
\label{thm:E_2-E_4-E_6-q^2-(1,2/3)}
{\it
For $q\in\mathbb{C}$ with $|q|<1,$ we have 
\begin{align*}
E_4(q^2)=&\frac{1}{256} (3 P+Q) \left(243 P^3-405 P^2 Q+225 P Q^2+Q^3\right),    \\
E_6(q^2)=&-\frac{\left(27 P^2-18 P Q-Q^2\right) \left(729 P^4-972 P^3 Q+270 P^2 Q^2-540 P Q^3+Q^4\right)}{4096}, 
\end{align*}
and
\begin{align*}
& q^2(q^2;q^2)^{24}_{\infty} 
=
\frac{ P Q^3 (P-Q)^2 (9 P-Q)^6}{2^{24}},  \\
&j(q^2)
=
\frac{(3 P+Q)^3 \left(243 P^3-405 P^2 Q+225 P Q^2+Q^3\right)^3}{P Q^3 (P-Q)^2 (9 P-Q)^6}  
=\frac{(3 f+1)^3 \left(243 f^3-405 f^2+225 f+1\right)^3}{ f   (f-1)^2 (9 f-1)^6}, 
\end{align*}
where 
\begin{equation*}
f(q)=P(q)/Q(q)=
\frac
{
(q^2; q^2)_{\infty}^4 (q^3; q^3)_{\infty}^8
}
{
(q;q)_{\infty}^8
(q^6; q^6)_{\infty}^4
}. 
\end{equation*}
}
\end{theorem}

\begin{proof}
Changing $q\rightarrow q^2$ in Ramanujan's system (\ref{eqn:Ramanujan-ODE}) implies 
\begin{equation*}
q\frac{d}{dq} E_2(q^2)=\frac{ \left\{ E_2(q^2)  \right\}^2-  E_4(q^2)  }{6}, \,\,
\mathrm{and} \,\, 
q\frac{d}{dq} E_4(q^2)=\frac{ 2 \left(   E_2(q^2)E_4(q^2)-E_6(q^2)  \right)  }{3}. 
\end{equation*}
From Theorem \ref{thm:E_2(q)-level-6-(1,2/3)}, we have 
\begin{equation*}
q\frac{d}{dq} E_2(q^2)=
\frac{9 P^2 R}{32}-\frac{P Q^3}{8}-\frac{3 P Q R}{16}-\frac{Q^2 R}{96}+\frac{R^2}{24}. 
\end{equation*}
Therefore, it follows that 
\begin{equation*}
E_4(q^2) =\left\{ E_2(q^2)  \right\}^2-6q\frac{d}{dq} E_2(q^2)
=\frac{1}{256} (3 P+Q) \left(243 P^3-405 P^2 Q+225 P Q^2+Q^3\right). 
\end{equation*}
\par
In the same way, we can obtain 
\begin{align*}
E_6(q^2)=&E_2(q^2)E_4(q^2)-\frac32 q\frac{d}{dq} E_4(q^2)  \\
=&-\frac{\left(27 P^2-18 P Q-Q^2\right) \left(729 P^4-972 P^3 Q+270 P^2 Q^2-540 P Q^3+Q^4\right)}{4096}. 
\end{align*}
\end{proof}

\begin{theorem}
\label{thm:E_2-E_4-E_6-q^3-(1,2/3)}
{\it
For $q\in\mathbb{C}$ with $|q|<1,$ 
we have 
\begin{align*}
E_4(q^3)=&\frac{1}{16} (3 P-Q) \left(3 P^3-3 P^2 Q+9 P Q^2-Q^3\right),  \\
E_6(q^3)=&-\frac{1}{64} \left(3 P^2+6 P Q-Q^2\right) \left(9 P^4-36 P^3 Q+30 P^2 Q^2-12 P Q^3+Q^4\right),
\end{align*}
and
\begin{align*}
q^3(q^3;q^3)^{24}_{\infty} 
=&
\frac{P^6 Q^2 (P-Q)^3 (9 P-Q)}{4096},  \\
j(q^3)
=&
\frac{(3 P-Q)^3 \left(3 P^3-3 P^2 Q+9 P Q^2-Q^3\right)^3}{P^6 Q^2 (P-Q)^3 (9 P-Q)} 
=\frac{(3 f-1)^3 \left(3 f^3-3 f^2+9 f-1\right)^3}{  f^6 (f-1)^3  (9 f-1)}, 
\end{align*}
where 
\begin{equation*}
f(q)=P(q)/Q(q)=
\frac
{
(q^2; q^2)_{\infty}^4 (q^3; q^3)_{\infty}^8
}
{
(q;q)_{\infty}^8
(q^6; q^6)_{\infty}^4
}. 
\end{equation*}
}
\end{theorem}

\begin{proof}
The theorem can be proved in the same as Theorem \ref{thm:E_2-E_4-E_6-q^2-(1,2/3)}. 
\end{proof}

\begin{theorem}
\label{thm:E_2-E_4-E_6-q^6-(1,2/3)}
{\it
For $q\in\mathbb{C}$ with $|q|<1,$ we have 
\begin{align*}
E_4(q^6)=&\frac{1}{256} (3 P+Q) \left(3 P^3+75 P^2 Q-15 P Q^2+Q^3\right),    \\
E_6(q^6)=&-\frac{\left(3 P^2+6 P Q-Q^2\right) \left(9 P^4-540 P^3 Q+30 P^2 Q^2-12 P Q^3+Q^4\right)}{4096}, 
\end{align*}
and
\begin{align*}
&q^6(q^6;q^6)^{24}_{\infty} 
=
\frac{ P^3 Q (P-Q)^6 (9 P-Q)^2}{2^{24}},  \\
&j(q^6)
=
\frac{(3 P+Q)^3 \left(3 P^3+75 P^2 Q-15 P Q^2+Q^3\right)^3}{P^3 Q (P-Q)^6 (9 P-Q)^2}  
=\frac{(3 f+1)^3 \left(3 f^3+75 f^2-15 f+1\right)^3}{ f^3  (f-1)^6  (9 f-1)^2},
\end{align*}
where 
\begin{equation*}
f(q)=P(q)/Q(q)=
\frac
{
(q^2; q^2)_{\infty}^4 (q^3; q^3)_{\infty}^8
}
{
(q;q)_{\infty}^8
(q^6; q^6)_{\infty}^4
}. 
\end{equation*}
}
\end{theorem}

\begin{proof}
The theorem can be proved in the same as Theorem \ref{thm:E_2-E_4-E_6-q^2-(1,2/3)}. 
\end{proof}

\subsection{Note on the cubic theta functions}

\begin{theorem}
\label{thm:cubic-q-(1,2/3)}
{\it
For $q\in\mathbb{C}$ with $|q|<1,$ we have 
\begin{equation*}
a(q)=\frac32 P-\frac12 Q, \,\,
b^3(q)=\frac18(9 PQ^2-Q^3), \,\,c^3(q)=\frac{27}{8}(P^3-P^2Q). 
\end{equation*}
}
\end{theorem}

\begin{proof}
The formula for $a(q)$ can be obtained by direct calculation. 
By Theorem \ref{thm-level6-E_2(q)-(1,2/3)}, we have 
\begin{equation*}
q\frac{d}{d q} a(q)=
\frac{27 P^3}{32}-\frac{27 P^2 Q}{32}-\frac{3 P Q^2}{32}+\frac{P R}{8}+\frac{Q^3}{96}-\frac{Q R}{24}. 
\end{equation*}
Moreover, by equation (\ref{eqn:Matsda-ODE-1}), we obtain 
\begin{align*}
q\frac{d}{d q} a(q)=&\frac{1}{12} 
\left\{
3a^3(q)+a(q)E_2(q)-4b^3(q) 
\right\}  \\
=&
\frac{27 P^3}{32}-\frac{27 P^2 Q}{32}+\frac{9 P Q^2}{32}+\frac{P R}{8}-\frac{Q^3}{32}-\frac{Q R}{24}
-\frac13 b^3(q),
\end{align*}
which yields the formula of $b^3(q)$. 
The formula of $c^3(q)$ can be obtained by the cubic identity, $a^3(q)=b^3(q)+c^3(q).$ 
\end{proof}

\begin{theorem}
\label{thm:cubic-q^2-(1,2/3)}
{\it
For $q\in\mathbb{C}$ with $|q|<1,$ we have 
\begin{equation*}
a(q^2)=\frac34 P+\frac14 Q, \,\,
b^3(q^2)=\frac{1}{64}(81 PQ^2-18 P Q^2+Q^3), \,\,c^3(q^2)=\frac{27}{64}(P^3-2P^2Q+P Q^2). 
\end{equation*}
}
\end{theorem}

\begin{proof}
By Theorem \ref{thm-level6-E_2(q)-(1,2/3)}, we have 
\begin{equation*}
q\frac{d}{dq} a(q^2)=\frac{27 P^3}{64}-\frac{45 P^2 Q}{64}+\frac{13 P Q^2}{64}+\frac{P R}{16}-\frac{Q^3}{192}+\frac{Q R}{48}. 
\end{equation*}
\par
Changing $q\rightarrow q^2$ in equation (\ref{eqn:Matsda-ODE-1}) yields  
\begin{equation*}
q\frac{d}{dq} a(q^2)
=
\frac16
\left\{
3 a^3(q^2)+a(q^2)E_2(q^2)-4b^3(q^2)
\right\}.
\end{equation*}
Theorem \ref{thm:E_2(q)-level-6-(1,2/3)} implies that  
\begin{align*}
b^3(q^2)=&\frac14 \left\{3 a^3(q^2)+a(q^2)E_2(q^2)-6 q\frac{d}{dq} a(q^2)\right\}   \\
=&\frac{1}{64}(81 PQ^2-18 P Q^2+Q^3).
\end{align*}
The formula of $c^3(q^2)$ can be obtained by the cubic identity, $a^3(q)=b^3(q)+c^3(q).$ 
\end{proof}

\begin{theorem}
\label{thm:cubic-ODE-(l,m,n)}
{\it
For $q\in\mathbb{C}$ with $|q|<1,$ set 
$l=a(q), \,\,m=a(q^2)$ and $n=E_2(q).$ Then, we have 
\begin{align*}
q\frac{d}{dq} l=&\frac{l^3+6 l^2 m-8 m^3+l n}{12}, \,\,
 q\frac{d}{dq} m=\frac{2 l^3+l^2 m-4 l m^2+m n }{12}, \\
q\frac{d}{d q} n=&
\frac{-5 l^4-12 l^3 m+16 l m^3+n^2}{12}. 
\end{align*}
}
\end{theorem}

\begin{proof}
Theorems \ref{thm:cubic-q-(1,2/3)} and \ref{thm:cubic-q^2-(1,2/3)} yield 
\begin{equation*}
l=\frac32 P-\frac12 Q, \,\,m=\frac34 P+\frac14 Q \,\,\mathrm{and} \,\,
P=\frac13 l+\frac23 m, \,\,Q=-l+2m. 
\end{equation*} 
The theorem follows from Theorem \ref{thm-level6-E_2(q)-(1,2/3)}. 
\end{proof}

\subsection{Application to number theory}
Throughout this subsection, we set 
$
l=a(q), \,\,m=a(q^2). 
$

\begin{theorem}
\label{thm:two-sums-x^2+xy+y^2}
{\it
For $q\in\mathbb{C}$ with $|q|<1,$ 
we have 
\begin{align}
a^2(q)=&1+12\sum_{n=1}^{\infty} (\sigma(n)-3\sigma(n/3))q^n,   \label{eqn:a(q)-1-1}  \\
a(q)a(q^2)=&1+6\sum_{n=1}^{\infty} (\sigma(n)-2\sigma(n/2)+3\sigma(n/3)-6\sigma(n/6))q^n.   \label{eqn:a(q)-1-2}
\end{align}
}
\end{theorem}

\begin{proof}
The formula (\ref{eqn:a(q)-1-1}) follows from equation (\ref{eqn:2nd-deri-(1,1/3)}). 
For the detailed proof, see \cite{Matsuda-1}. 
\par
Theorems \ref{thm:E_2(q)-level-6-(1,2/3)}, \ref{thm:cubic-q-(1,2/3)}, and \ref{thm:cubic-q^2-(1,2/3)} imply that 
\begin{align*}
a(q) a(q^2)=&\left(\frac32P-\frac12 Q\right)  \left(\frac34 P+\frac14 Q\right)  
=\frac14\left\{
-E_2(q)+2E_2(q^2)-3E_2(q^3)+6E_2(q^6)
\right\},
\end{align*}
which proves the theorem. 
\end{proof}

\begin{theorem}
\label{thm:three-sums-x^2+xy+y^2}
{\it
For $q\in\mathbb{C}$ with $|q|<1,$ 
we have 
\begin{align}
a^3(q)=&1
-9
\sum_{n=1}^{\infty} \left( \sum_{d|n} d^2\left( \frac{d}{3}  \right)  \right)q^n
+27
\sum_{n=1}^{\infty} \left( \sum_{d|n} d^2\left( \frac{n/d}{3}  \right)  \right)  q^n,    \label{eqn:a(q)-1-1-1} \\
a^2(q)a(q^2)=&
1
+3
\sum_{n=1}^{\infty} \left( \sum_{d|n} d^2\left( \frac{d}{3}  \right)  \right)q^n
+9
\sum_{n=1}^{\infty} \left( \sum_{d|n} d^2\left( \frac{n/d}{3}  \right)  \right)  q^n   \notag  \\
&-12
\sum_{n=1}^{\infty} \left( \sum_{d|n} d^2\left( \frac{d}{3}  \right)  \right)q^{2n}
+36
\sum_{n=1}^{\infty} \left( \sum_{d|n} d^2\left( \frac{n/d}{3}  \right)  \right)  q^{2n},      \label{eqn:a(q)-1-1-2}   \\
a(q)a^2(q^2)=&
1
-3
\sum_{n=1}^{\infty} \left( \sum_{d|n} d^2\left( \frac{d}{3}  \right)  \right)q^n
+9
\sum_{n=1}^{\infty} \left( \sum_{d|n} d^2\left( \frac{n/d}{3}  \right)  \right)  q^n   \notag  \\
&-6
\sum_{n=1}^{\infty} \left( \sum_{d|n} d^2\left( \frac{d}{3}  \right)  \right)q^{2n}
-18
\sum_{n=1}^{\infty} \left( \sum_{d|n} d^2\left( \frac{n/d}{3}  \right)  \right)  q^{2n}.      \label{eqn:a(q)-1-2-2} 
\end{align}
}
\end{theorem}

\begin{proof}
For the proof of equation (\ref{eqn:a(q)-1-1-1}), see \cite{Matsuda-1}. 
By Theorems \ref{thm:cubic-q-(1,2/3)} and \ref{thm:cubic-q^2-(1,2/3)}, 
we note that 
$P=(l+2m)/3$ and $Q=-l+2m,$ 
which yields 
\begin{equation*}
b^3(q)=\frac{l^3}{2}-\frac{3 l^2 m}{2}+2 m^3, \,\,
b^3(q^2)=-\frac{l^3}{4}+\frac{3 l m^2}{4}+\frac{m^3}{2}. 
\end{equation*}
Therefore, it follows that 
\begin{align}
a^2(q)a(q^2)=&-\frac13 b^3(q)+\frac13 c^3(q)+\frac43 a^3(q^2),  \label{eqn:a(q)-1-1-2-b-c}  \\
a(q)a^2(q^2)=&\frac23 b^3(q^2)-\frac23 c^3(q^2)+\frac13 a^3(q),  \label{eqn:a(q)-1-2-2-b-c}
\end{align}
which proves the theorem. 
\end{proof}

\begin{theorem}
\label{thm:four-sums-x^2+xy+y^2}
{\it
For $q\in\mathbb{C}$ with $|q|<1,$ 
we have 
\begin{align}
a^4(q)=&1+24\sum_{n=1}^{\infty}(\sigma_3(n)+9\sigma_3(n/3))q^n,    \label{eqn:a(q)-1-1-1-1} \\
a^3(q)a(q^2)=&
1
+
6\sum_{n=1}^{\infty} (3\sigma_3(n)-8\sigma_3 (n/2)-27 \sigma_3  (n/3)+72 \sigma_3(n/6)) q^n,     \label{eqn:a(q)-1-1-1-2}   \\
a(q)a^3(q^2)=&
1
+
6\sum_{n=1}^{\infty} (\sigma_3(n)-6 \sigma_3 (n/2)-9 \sigma_3  (n/3)+54 \sigma_3  (n/6)) q^n.     \label{eqn:a(q)-1-2-2-2} 
\end{align}
}
\end{theorem}

\begin{proof}
By Theorems \ref{thm:cubic-q-(1,2/3)} and \ref{thm:cubic-q^2-(1,2/3)}, 
we note that 
$P=(l+2m)/3$ and $Q=-l+2m.$ 
Theorems \ref{thm:E_2-E_4-E_6-q-(1,2/3)}, \ref{thm:E_2-E_4-E_6-q^2-(1,2/3)}, and \ref{thm:E_2-E_4-E_6-q^3-(1,2/3)}, 
we have 
\begin{align}
E_4(q)=&5 l^4+12 l^3 m-16 l m^3,   \label{eqn:E_4-q-l,m} \\
E_4(q^2)=&2 l^3 m-6 l m^3+5 m^4,   \label{eqn:E_4-q^2-l,m} \\
E_4(q^3)=&\frac{5 l^4}{9}-\frac{4 l^3 m}{3}+\frac{16 l m^3}{9}.   \label{eqn:E_4-q^3-l,m}
\end{align}
Equations (\ref{eqn:E_4-q-l,m}) and (\ref{eqn:E_4-q^3-l,m}) yield 
\begin{equation*}
a^4(q)=\frac{1}{10}\left\{ E_4(q)+9E_4(q^3)   \right\}, \,\,\mathrm{and}  \,\,
a^4(q^2)=\frac{1}{10}\left\{ E_4(q^2)+9E_4(q^6)   \right\}.
\end{equation*}
Equations (\ref{eqn:E_4-q-l,m}) and (\ref{eqn:E_4-q^2-l,m}) imply that 
\begin{equation*}
12 l^3m-16 l m^3=\frac12 E_4(q)-\frac92 E_4(q^3), \,\,
2l^3m-6 l m^3=\frac12 E_4(q^2)-\frac92E_4(q^6),
\end{equation*}
which proves the theorem. 
\end{proof}

\section{ODEs satisfied by modular forms of level six (2)}
\label{sec:proof:level-six-E_2(q^6)}
The main result of this section is as follows:

\begin{theorem}
\label{thm-level6-E_2(q^6)-(2/3,1)}
{\it
For $q\in\mathbb{C}$ with $|q|<1,$ set
\begin{align*}
P(q)=&\frac{ (q^3; q^3)_{\infty} (q^2; q^2)_{\infty}^6  }{ (q; q)_{\infty}^3 (q^6; q^6)_{\infty}^2 }=1+ 3 \sum_{n=1}^{\infty} (d_{1,6}(n)-d_{5,6}(n)) q^{n},   \\
Q(q)=&q \frac{(q; q)_{\infty}  (q^6; q^6)_{\infty}^6   }{  (q^2; q^2)_{\infty}^2 (q^3; q^3)_{\infty}^3   }=   \sum_{n=1}^{\infty}(d_{1,6}(n)-d_{5,6}(n))q^{n}-2  \sum_{n=1}^{\infty}  (d_{1,3}(n)-d_{2,3}(n)) q^{2n},   \\
R(q)=&E_2(q^6)=1-24\sum_{n=1}^{\infty} \sigma_1(n) q^{6n}. 
\end{align*}
Then, we have 
\begin{align*}
q\frac{d}{dq}P=&\frac{-P^3+12P^2Q-15PQ^2+PR}{2},    \,\,
q\frac{d}{dq}Q=\frac{P^2Q-8PQ^2+3Q^3+QR}{2} \\
q\frac{d}{dq}R=&\frac{-P^4+12P^3Q-30P^2Q^2+12PQ^3-9Q^4+R^2}{2}. 
\end{align*}
}
\end{theorem}

\begin{proof}
The theorem can be proved in the same way as Theorem \ref{thm-level6-E_2(q)-(1,2/3)}. 
We first set 
\begin{equation*}
X=
\frac
{
\theta^{\prime}
\left[
\begin{array}{c}
\frac23 \\
1
\end{array}
\right]
}
{
\theta
\left[
\begin{array}{c}
\frac23 \\
1
\end{array}
\right]
},  \,\,
Y=
\exp\left( \frac{4\pi i}{3} \right)
\frac
{
\theta^{\prime}
\left[
\begin{array}{c}
1 \\
1
\end{array}
\right]
\theta
\left[
\begin{array}{c}
\frac23 \\
1
\end{array}
\right]
}
{
\theta
\left[
\begin{array}{c}
0 \\
1
\end{array}
\right]
\theta
\left[
\begin{array}{c}
\frac13 \\
1
\end{array}
\right]
},  \,\,
Z=
\frac
{
\theta^{\prime  \prime  \prime}
\left[
\begin{array}{c}
1 \\
1
\end{array}
\right]
}
{
\theta^{\prime  }
\left[
\begin{array}{c}
1 \\
1
\end{array}
\right]
}, \,\,
W=
\frac
{
\theta^{\prime}
\left[
\begin{array}{c}
\frac13 \\
1
\end{array}
\right]
}
{
\theta
\left[
\begin{array}{c}
\frac13 \\
1
\end{array}
\right]
}
=\frac{X-Y}{2}.
\end{equation*}
In particular, we note that 
\begin{equation*}
X=\frac{2\pi i}{3} P(y), 
\,\,
Y=
2
\pi iQ(y), \,\,
Z=
-\pi^2 R(y^6),
\end{equation*}
where $y=\exp(2 \pi i \tau/6)$ and $q=\exp(2\pi i \tau).$
\par
Applying the residue theorem to the elliptic functions, 
\begin{equation*}
\frac
{
\theta^3
\left[
\begin{array}{c}
\frac13 \\
1
\end{array}
\right]
(z)
}
{
\theta^3
\left[
\begin{array}{c}
1 \\
1
\end{array}
\right]
(z)
}, \quad
\frac
{
\theta^2
\left[
\begin{array}{c}
\frac23 \\
1
\end{array}
\right]
(z)
\theta
\left[
\begin{array}{c}
-\frac13 \\
1
\end{array}
\right]
(z)
}
{
\theta^3
\left[
\begin{array}{c}
1 \\
1
\end{array}
\right]
(z)
}, \,\,
\frac
{
\theta
\left[
\begin{array}{c}
\frac13 \\
1
\end{array}
\right]
(z)
\theta
\left[
\begin{array}{c}
\frac23 \\
1
\end{array}
\right]
(z)
\theta
\left[
\begin{array}{c}
0 \\
1
\end{array}
\right]
(z)
}
{
\theta^3
\left[
\begin{array}{c}
1 \\
1
\end{array}
\right]
(z)
},
\end{equation*}
and 
\begin{equation*}
\frac
{
\theta^3
\left[
\begin{array}{c}
\frac23\\
1
\end{array}
\right](z)
\theta
\left[
\begin{array}{c}
0 \\
1
\end{array}
\right](z)
}
{
\theta^4
\left[
\begin{array}{c}
1 \\
1
\end{array}
\right](z)
}, \quad
\frac
{
\theta^2
\left[
\begin{array}{c}
\frac13 \\
1
\end{array}
\right](z)
\theta^2
\left[
\begin{array}{c}
\frac23 \\
1
\end{array}
\right](z)
}
{
\theta^4
\left[
\begin{array}{c}
1 \\
1
\end{array}
\right](z)
},
\end{equation*}
we have 
\begin{equation*}
4\pi i X^{\prime}=
-\frac32 X^3+6 X^2 Y-\frac52 X Y^2+\frac23 XZ,
\,\,\mathrm{and}\,\,
4 \pi i Y^{\prime}=
\frac32 X^2 Y-4 X Y^2+\frac12 Y^3+\frac23 Y Z,
\end{equation*}
where ${ }^{\prime}=d/d\tau.$ 
\par
Applying the residue theorem to 
$
\theta^4
\left[
\begin{array}{c}
\frac23 \\
1
\end{array}
\right]
(z)
\theta
\left[
\begin{array}{c}
\frac13 \\
1
\end{array}
\right]
(z)
/
\theta^5
\left[
\begin{array}{c}
1 \\
1
\end{array}
\right]
(z),
$ 
we have 
\begin{equation*}
4\pi i Z^{\prime}=
-\frac{27 X^4}{8}+\frac{27 X^3 Y}{2}-\frac{45 X^2 Y^2}{4}+\frac{3 X Y^3}{2}-\frac{3 Y^4}{8}+\frac{2 Z^2}{3}.
\end{equation*}
Changing $\tau\longrightarrow 6 \tau,$ we obtain the theorem. 
\end{proof}

\subsubsection*{Remark}
Throughout this section, 
for $q\in\mathbb{C}$ with $|q|<1,$ we set
\begin{align*}
P(q)=&\frac{ (q^3; q^3)_{\infty} (q^2; q^2)_{\infty}^6  }{ (q; q)_{\infty}^3 (q^6; q^6)_{\infty}^2 }=1+ 3 \sum_{n=1}^{\infty} (d_{1,6}(n)-d_{5,6}(n)) q^{n},   \\
Q(q)=&q \frac{(q; q)_{\infty}  (q^6; q^6)_{\infty}^6   }{  (q^2; q^2)_{\infty}^2 (q^3; q^3)_{\infty}^3   }=   \sum_{n=1}^{\infty}(d_{1,6}(n)-d_{5,6}(n))q^{n}-2  \sum_{n=1}^{\infty}  (d_{1,3}(n)-d_{2,3}(n)) q^{2n},   \\
R(q)=&E_2(q^6)=1-24\sum_{n=1}^{\infty} \sigma_1(n) q^{6n}. 
\end{align*}

\subsection{Note on $E_4(q)$ and $E_6(q)$ $\cdots$ (1)}

\begin{theorem}
\label{thm-E4-E6-level-6-q^6}
{\it
For $q\in\mathbb{C}$ with $|q|<1,$ we have 
\begin{align*}
E_4(q^6)=&(P-3 Q) \left(P^3-9 P^2 Q+3 P Q^2-3 Q^3\right),  \\
E_6(q^6)=&
\left(P^2-6 P Q-3 Q^2\right) \left(P^4-12 P^3 Q+30 P^2 Q^2-36 P Q^3+9 Q^4\right),
\end{align*}
and
\begin{align*}
&q^6(q^6;q^6)_{\infty}^{24}=
P^2Q^6(P-9 Q)(P-Q)^3, \\
&j(q^6)=
\frac{(P-3 Q)^3 \left(P^3-9 P^2 Q+3 P Q^2-3 Q^3\right)^3}{P^2 Q^6 (P-9 Q) (P-Q)^3}  
=
\frac{(3 g-1)^3 \left(3 g^3-3 g^2+9 g-1\right)^3}{ g^6 (g-1)^3  (9 g-1)},
\end{align*}
where 
\begin{equation*}
g(q)=Q(q)/P(q)
=
q
\frac
{
(q;q)_{\infty}^4 (q^6;q^6)_{\infty}^8
}
{
(q^2;q^2)_{\infty}^8 (q^3;q^3)_{\infty}^4
}. 
\end{equation*}
}
\end{theorem}

\begin{proof}
The theorem can be proved in the same way as Theorem \ref{thm:E_2-E_4-E_6-q-(1,2/3)}. 
Substitute $z=-\tau/6$ in 
\begin{equation*}
\wp=\wp(z;1,\tau)
=
\frac13
\frac
{
\theta^{\prime  \prime  \prime}
\left[
\begin{array}{c}
1 \\
1
\end{array}
\right]
}
{
\theta^{\prime  }
\left[
\begin{array}{c}
1 \\
1
\end{array}
\right]
}
-
\frac{d^2}{dz^2}
\log
\theta
\left[
\begin{array}{c}
1 \\
1
\end{array}
\right](z),
\end{equation*} 
and 
\begin{equation*}
(\wp^{\prime})^2=4\wp^3-g_2\wp-g_3, \,\,  
\wp^{\prime \prime}(z)=6\wp^2(z)-\frac12g_2.
\end{equation*}
\end{proof}

\subsection{Note on $E_2(q)$}

\begin{theorem}
\label{thm:E_2(q)-level-6-(2/3,1)}
{\it
For $q\in\mathbb{C}$ with $|q|<1,$ we have 
\begin{align*}
P^2=&\frac18 
\left(
-2 E_2(q)+E_2(q^2)+9 E_2(q^6)
\right), \,\,
PQ=\frac{1}{24}
\left(
-E_2(q)+E_2(q^2)+3 E_2(q^3)-3 E_2(q^6)
\right), \\
Q^2=&\frac{1}{24}
\left(
-E_2(q^2)+2 E_2(q^3)-E_2(q^6)
\right),
\end{align*}
which implies that 
\begin{align*}
E_2(q)=&-5 P^2+6 P Q-9 Q^2+6 R, \,\,
E_2(q^2)=-2 P^2+12 P Q-18 Q^2+3 R, \\
E_2(q^3)=&-P^2+6 PQ+3 Q^2+2 R. 
\end{align*}
}
\end{theorem}

\begin{proof}
The theorem can be proved in the same way as Theorem \ref{thm:E_2(q)-level-6-(1,2/3)}. 
\end{proof}

\subsection{Note on $E_4(q)$ and $E_6(q)$ $\cdots$ (2)}

\begin{theorem}
\label{thm:E_2-E_4-E_6-q-(2/3,1)}
{\it
For $q\in\mathbb{C}$ with $|q|<1,$ we have 
\begin{align*}
E_4(q)=&(P+3 Q) \left(P^3+225 P^2 Q-405 P Q^2+243 Q^3\right),    \\
E_6(q)=&\left(P^2+18 P Q-27 Q^2\right) \left(P^4-540 P^3 Q+270 P^2 Q^2-972 P Q^3+729 Q^4\right), 
\end{align*}
and
\begin{align*}
& q(q;q)^{24}_{\infty} 
=
P^3 Q (P-9 Q)^6 (P-Q)^2,  \\
&j(q)
=\frac{(P+3 Q)^3 \left(P^3+225 P^2 Q-405 P Q^2+243 Q^3\right)^3}{P^3 Q (P-9 Q)^6 (P-Q)^2}  
=\frac{(3 g+1)^3 \left(243 g^3-405 g^2+225 g+1\right)^3}{ g(g-1)^2 (9 g-1)^6  }, 
\end{align*}
where 
\begin{equation*}
g(q)=Q(q)/P(q)
=
q
\frac
{
(q;q)_{\infty}^4 (q^6;q^6)_{\infty}^8
}
{
(q^2;q^2)_{\infty}^8 (q^3;q^3)_{\infty}^4
}. 
\end{equation*}
}
\end{theorem}

\begin{proof}
The theorem can be proved in the same as Theorem \ref{thm:E_2-E_4-E_6-q^2-(1,2/3)}. 
\end{proof}

\begin{theorem}
\label{thm:E_2-E_4-E_6-q^2-(2/3,1)}
{\it
For $q\in\mathbb{C}$ with $|q|<1,$ we have 
\begin{align*}
E_4(q^2)=&(P-3 Q) \left(P^3-9 P^2 Q+243 P Q^2-243 Q^3\right),    \\
E_6(q^2)=&\left(P^2+18 P Q-27 Q^2\right) \left(P^4-36 P^3 Q+270 P^2 Q^2-972 P Q^3+729 Q^4\right), 
\end{align*}
and
\begin{align*}
&q^2(q^2;q^2)^{24}_{\infty} 
=
P^6 Q^2 (P-9 Q)^3 (P-Q),  \\
&j(q^2)
=\frac{(P-3 Q)^3 \left(P^3-9 P^2 Q+243 P Q^2-243 Q^3\right)^3}{P^6 Q^2 (P-9 Q)^3 (P-Q)}  
=&\frac{(3 g-1)^3 \left(243 g^3-243 g^2+9 g-1\right)^3}{(g-1) g^2 (9 g-1)^3}, 
\end{align*}
where 
\begin{equation*}
g(q)=Q(q)/P(q)
=
q
\frac
{
(q;q)_{\infty}^4 (q^6;q^6)_{\infty}^8
}
{
(q^2;q^2)_{\infty}^8 (q^3;q^3)_{\infty}^4
}. 
\end{equation*}
}
\end{theorem}

\begin{proof}
The theorem can be proved in the same as Theorem \ref{thm:E_2-E_4-E_6-q^2-(1,2/3)}. 
\end{proof}

\begin{theorem}
\label{thm:E_2-E_4-E_6-q^3-(2/3,1)}
{\it
For $q\in\mathbb{C}$ with $|q|<1,$ we have 
\begin{align*}
E_4(q^3)=&(P+3 Q) \left(P^3-15 P^2 Q+75 P Q^2+3 Q^3\right),    \\
E_6(q^3)=&\left(P^2-6 P Q-3 Q^2\right) \left(P^4-12 P^3 Q+30 P^2 Q^2-540 P Q^3+9 Q^4\right), 
\end{align*}
and
\begin{align*}
&q^3(q^3;q^3)^{24}_{\infty} 
=
 P Q^3 (P-9 Q)^2 (P-Q)^6,  \\
&j(q^3)
=\frac{(P+3 Q)^3 \left(P^3-15 P^2 Q+75 P Q^2+3 Q^3\right)^3}{P Q^3 (P-9 Q)^2 (P-Q)^6}  
=\frac{(3 g+1)^3 \left(3 g^3+75 g^2-15 g+1\right)^3}{  g^3 (g-1)^6 (9 g-1)^2 }, 
\end{align*}
where 
\begin{equation*}
g(q)=Q(q)/P(q)
=
q
\frac
{
(q;q)_{\infty}^4 (q^6;q^6)_{\infty}^8
}
{
(q^2;q^2)_{\infty}^8 (q^3;q^3)_{\infty}^4
}. 
\end{equation*}
}
\end{theorem}

\begin{proof}
The theorem can be proved in the same as Theorem \ref{thm:E_2-E_4-E_6-q^2-(1,2/3)}. 
\end{proof}

\subsection{Note on the cubic theta functions}

\begin{theorem}
\label{thm:cubic-q-(2/3,1)}
{\it
For $q\in\mathbb{C}$ with $|q|<1,$ we have 
\begin{equation*}
a(q)=P+3Q, \,\,
b^3(q)=P^3-18 P^2 Q+81 P Q^2, \,\,
c^3(q)=27 P^2 Q-54 P Q^2+27 Q^3. 
\end{equation*}
}
\end{theorem}

\begin{proof}
The theorem can be proved in the same way as Theorem \ref{thm:cubic-q-(1,2/3)}. 
\end{proof}

\begin{theorem}
\label{thm:cubic-q^2-(2/3,1)}
{\it
For $q\in\mathbb{C}$ with $|q|<1,$ we have 
\begin{equation*}
a(q^2)=P-3Q, \,\,
b^3(q^2)=P^3-9P^2 Q, \,\,
c^3(q^2)=27 P Q^2-27 Q^3. 
\end{equation*}
}
\end{theorem}

\begin{proof}
The theorem can be proved in the same way as Theorem \ref{thm:cubic-q^2-(1,2/3)}. 
\end{proof}

\subsection{Note on Theorems \ref{thm-level6-E_2(q)-(1,2/3)} and \ref{thm-level6-E_2(q^6)-(2/3,1)} }
Theorems \ref{thm:cubic-q-(1,2/3)}, \ref{thm:cubic-q^2-(1,2/3)},  \ref{thm:cubic-q-(2/3,1)}, and \ref{thm:cubic-q^2-(2/3,1)} show that 
the differential fields of Theorems \ref{thm-level6-E_2(q)-(1,2/3)} and \ref{thm-level6-E_2(q^6)-(2/3,1)} are the same and can be expressed by 
$\mathbb{C}(a(q), a(q^2), E_2(q)).$ 
\par
Then, the differential field, $\mathbb{C}(a(q), a(q^2), E_2(q))$ contains the following differential subfields: 
\begin{align*}
&\mathbb{C}\left( E_2(q), E_4(q), E_6(q)    \right),  \mathbb{C}\left( E_2(q^2), E_4(q^2), E_6(q^2)    \right), \,\,  
\mathbb{C}\left( E_2(q^3), E_4(q^3), E_6(q^3)    \right), \\           
&\mathbb{C}\left( E_2(q^6), E_4(q^6), E_6(q^6)    \right), \,\,        \mathbb{C}\left( a(q), E_2(q), b^3(q)  \right), \,\,
\mathbb{C}\left( a(q^2), E_2(q^2), b^3(q^2)  \right).
\end{align*}

\section{ODEs satisfied by modular forms of level six (3)}
\label{sec:proof:level-six-E_2(q)-split}
The main result of this section is as follows:

\begin{theorem}
\label{thm-level6-E_2(q)-(0,1/3)-(0,2/3)}
{\it
For $q\in\mathbb{C}$ with $|q|<1,$ set
\begin{align*}
P(q)=&a(q)=1+6\sum_{n=1}^{\infty}( d_{1,3}(n)-d_{2,3}(n) ) q^n, \,\,
Q(q)= q^{\frac12} \frac{ (q^{\frac12}; q^{\frac12})_{\infty}  (q^3;q^3)_{\infty}^6 }{  (q;q)_{\infty}^2  (q^{\frac32}; q^{\frac32})_{\infty}^3 }=\sum_{n=1}^{\infty}(d_{1,3}^{*}(n)-d_{2,3}^{*}(n)) q^{\frac{n}{2}}, \\
R(q)=&q^{\frac12} \frac{   (q;q)_{\infty}   (q^{\frac32}; q^{\frac32})_{\infty}^3    (q^6;q^6)_{\infty}^3    }{   (q^{\frac12}; q^{\frac12})_{\infty}   (q^2;q^2)_{\infty}    (q^3;q^3)_{\infty}^3 }
=\sum_{n=1}^{\infty}(d_{1,6}^{*}(n)+d_{2,6}^{*}(n)    - d_{4,6}^{*}(n)-d_{5,6}^{*}(n)  ) q^{\frac{n}{2}}, \\
S(q)=&E_2(q)=1-24\sum_{n=1}^{\infty} \sigma_1(n) q^n. 
\end{align*}
Then, we have 
\begin{align*}
q\frac{d}{dq} P=&\frac{-P^3+108 P Q^2+216 Q^3+P S  }{12   }=\frac{ -P^3+108 P R^2-216 R^3+PS  }{ 12 }, \\
q\frac{d}{dq} Q=&\frac{ 5 P^2 Q-6 P Q^2-36Q^3+QS }{ 12 },  \,\,
q\frac{d}{dq} R=\frac{ 5 P^2 R+6 P R^2-36 R^3+ RS }{12  },  \\
q\frac{d}{dq} S=&\frac{ -P^4-216 P^2 Q^2-432 P Q^3+S^2 }{ 12  }=\frac{ -P^4-216P^2 R^2+432 P R^3+S^2 }{  12},
\end{align*}
and 
\begin{equation}
\label{eqn:(P,Q,R)-relation}
P(Q-R)+2(Q^2-QR+R^2)=0.
\end{equation}
}
\end{theorem}

\begin{proof}
The proof consists of Subsections 8.1-8.5.
\end{proof}

\subsubsection*{Remark}
Throughout this section, we set
\begin{align*}
P(q)=&a(q)=1+6\sum_{n=1}^{\infty}( d_{1,3}(n)-d_{2,3}(n) ) q^n, \,\,
Q(q)= q^{\frac12} \frac{ (q^{\frac12}; q^{\frac12})_{\infty}  (q^3;q^3)_{\infty}^6 }{  (q;q)_{\infty}^2  (q^{\frac32}; q^{\frac32})_{\infty}^3 }=\sum_{n=1}^{\infty}(d_{1,3}^{*}(n)-d_{2,3}^{*}(n)) q^{\frac{n}{2}}, \\
R(q)=&q^{\frac12} \frac{   (q;q)_{\infty}   (q^{\frac32}; q^{\frac32})_{\infty}^3    (q^6;q^6)_{\infty}^3    }{   (q^{\frac12}; q^{\frac12})_{\infty}   (q^2;q^2)_{\infty}    (q^3;q^3)_{\infty}^3 }
=\sum_{n=1}^{\infty}(d_{1,6}^{*}(n)+d_{2,6}^{*}(n)    - d_{4,6}^{*}(n)-d_{5,6}^{*}(n)  ) q^{\frac{n}{2}}, \\
S(q)=&E_2(q)=1-24\sum_{n=1}^{\infty} \sigma_1(n) q^n, \,\,\mathrm{for} \,\,q\in\mathbb{C}\,\, \mathrm{with}\,\, |q|<1.
\end{align*}
Moreover, we define
\begin{align*}
T(q)=&\frac{(q ;q)_{\infty}^6 (q^{\frac32};q^{\frac32})_{\infty}  }{(q^{\frac12};q^{\frac12})_{\infty}^3 (q^3;q^3)_{\infty}^2  } =
1+3\sum_{n=1}^{\infty} (d_{1,3}(n)-d_{2,3}(n)) q^{\frac{n}{2}}+3 \sum_{n=1}^{\infty} ( d_{1,3}(n)-d_{2,3}(n) )q^n=P+3Q,  \\
U(q)=&\frac{(q^{\frac12} ; q^{\frac12} )_{\infty}^3  (q^2; q^2)_{\infty}^3   (q^3; q^3)_{\infty}   }{  (q; q)_{\infty}^3 (q^{\frac32}; q^{\frac32})_{\infty}  (q^6; q^6)_{\infty}   }=
1-3  \sum_{n=1}^{\infty}(d_{1,6}^{*}(n)+d_{2,6}^{*}(n)  -d_{4,6}^{*}(n)-d_{5,6}^{*}(n))q^{\frac{n}{2}}    \\
&\hspace{60mm}+6  \sum_{n=1}^{\infty}  (d_{1,3}(n)-d_{2,3}(n)) q^{n}=P-3 R. 
\end{align*}

\begin{corollary}
{\it 
We have 
\begin{align*}
q\frac{d}{d q}
Q=&
\frac{ 27 Q^3-36 Q^2 T+5 Q T^2 +Q S }{12},   \,\,
q\frac{d}{d q}
T=
\frac{-27 Q^2 T+24 Q T^2-T^3 +S T }{12},   \\
q\frac{d}{d q}
S=&
\frac{-729 Q^4+972 Q^3 T-270 Q^2 T^2+12 Q T^3-T^4  +S^2}{12},  
\end{align*}
and
\begin{align*}
q\frac{d}{dq}
R
=&
\frac{27 R^3+36 R^2 U+5 R U^2  +R S }{12},  \,\,
q\frac{d}{dq}
U
= 
\frac{-27 R^2 U-24 R U^2-U^3+S U}{12},  \\
q\frac{d}{dq}
S
=&
\frac{-729 R^4-972 R^3 U-270 R^2 U^2-12 R U^3-U^4 +S^2 }{12}. 
\end{align*}
}
\end{corollary}

\subsection{Notations}
In this section, we set 
\begin{equation*}
X=
\frac{
\theta^{\prime}
\left[
\begin{array}{c}
1 \\
\frac13
\end{array}
\right]  
}
{ 
\theta
\left[
\begin{array}{c}
1 \\
\frac13
\end{array}
\right]
}, \,\,
Y=
\frac{
\theta^{\prime}
\left[
\begin{array}{c}
0 \\
\frac13
\end{array}
\right]  
}
{ 
\theta
\left[
\begin{array}{c}
0 \\
\frac13
\end{array}
\right]
}, \,\,
Z=
\frac{
\theta^{\prime}
\left[
\begin{array}{c}
0 \\
\frac23
\end{array}
\right]  
}
{ 
\theta
\left[
\begin{array}{c}
0 \\
\frac23
\end{array}
\right]
}, \,\,
W
=
\frac
{
\theta^{\prime  \prime  \prime}
\left[
\begin{array}{c}
1 \\
1
\end{array}
\right]
}
{
\theta^{\prime  }
\left[
\begin{array}{c}
1 \\
1
\end{array}
\right]
}, 
\end{equation*}
implying that 
\begin{equation*}
X=
-\frac{\pi}{\sqrt{3}}P(q), \,\,
Y=
-2\sqrt{3} \pi Q(q), \,\,
Z=
-2\sqrt{3} \pi R(q), \,\,
W=
-\pi^2 S(q)=-\pi^2 E_2(q),  \,\,q=\exp(2 \pi i \tau). 
\end{equation*}

\subsection{Relations between theta derivatives      }

\begin{proposition}
\label{prop:2nd-derivative-(0,1/3)-(0,2/3)}
{\it
For every $\tau\in\mathbb{H}^2,$ we have 
\begin{equation}
\label{eqn:2nd-deri-(1,1/3)-again}
3
\frac
{
\theta^{\prime \prime}
\left[
\begin{array}{c}
1 \\
\frac13
\end{array}
\right]
}
{
\theta
\left[
\begin{array}{c}
1 \\
\frac13
\end{array}
\right]
}
-
\frac
{
\theta^{\prime \prime \prime}
\left[
\begin{array}{c}
1 \\
1
\end{array}
\right]
}
{
\theta^{\prime}
\left[
\begin{array}{c}
1 \\
1
\end{array}
\right]
}
+
6
\left\{
\frac
{
\theta^{\prime}
\left[
\begin{array}{c}
1 \\
\frac13
\end{array}
\right]
}
{
\theta
\left[
\begin{array}{c}
1 \\
\frac13
\end{array}
\right]
}
\right\}^2=0,  
\end{equation}
\begin{equation}
\label{eqn:2nd-deri-(0,1/3)}
\frac
{
\theta^{\prime \prime}
\left[
\begin{array}{c}
1 \\
\frac13
\end{array}
\right]
}
{
\theta
\left[
\begin{array}{c}
1 \\
\frac13
\end{array}
\right]
}
+
2
\frac
{
\theta^{\prime \prime}
\left[
\begin{array}{c}
0 \\
\frac13
\end{array}
\right]
}
{
\theta
\left[
\begin{array}{c}
0 \\
\frac13
\end{array}
\right]
}
-
\frac
{
\theta^{\prime \prime \prime}
\left[
\begin{array}{c}
1 \\
1
\end{array}
\right]
}
{
\theta^{\prime}
\left[
\begin{array}{c}
1\\
1
\end{array}
\right]
}
+
4
\frac
{
\theta^{\prime}
\left[
\begin{array}{c}
0\\
\frac13
\end{array}
\right]
}
{
\theta
\left[
\begin{array}{c}
0\\
\frac13
\end{array}
\right]
}
\frac
{
\theta^{\prime}
\left[
\begin{array}{c}
1\\
\frac13
\end{array}
\right]
}
{
\theta
\left[
\begin{array}{c}
1\\
\frac13
\end{array}
\right]
}
+
2
\left\{
\frac
{
\theta^{\prime}
\left[
\begin{array}{c}
0\\
\frac13
\end{array}
\right]
}
{
\theta
\left[
\begin{array}{c}
0\\
\frac13
\end{array}
\right]
}
\right\}^2
=0,
\end{equation}
\begin{equation}
\label{eqn:2nd-deri-(0,2/3)}
\frac
{
\theta^{\prime \prime}
\left[
\begin{array}{c}
1 \\
\frac13
\end{array}
\right]
}
{
\theta
\left[
\begin{array}{c}
1 \\
\frac13
\end{array}
\right]
}
+
2
\frac
{
\theta^{\prime \prime}
\left[
\begin{array}{c}
0 \\
\frac23
\end{array}
\right]
}
{
\theta
\left[
\begin{array}{c}
0 \\
\frac23
\end{array}
\right]
}
-
\frac
{
\theta^{\prime \prime \prime}
\left[
\begin{array}{c}
1 \\
1
\end{array}
\right]
}
{
\theta^{\prime}
\left[
\begin{array}{c}
1\\
1
\end{array}
\right]
}
-
4
\frac
{
\theta^{\prime}
\left[
\begin{array}{c}
1\\
\frac13
\end{array}
\right]
}
{
\theta
\left[
\begin{array}{c}
1\\
\frac13
\end{array}
\right]
}
\frac
{
\theta^{\prime}
\left[
\begin{array}{c}
0\\
\frac23
\end{array}
\right]
}
{
\theta
\left[
\begin{array}{c}
0\\
\frac23
\end{array}
\right]
}
+
2
\left\{
\frac
{
\theta^{\prime}
\left[
\begin{array}{c}
0\\
\frac23
\end{array}
\right]
}
{
\theta
\left[
\begin{array}{c}
0\\
\frac23
\end{array}
\right]
}
\right\}^2
=0,
\end{equation}
\begin{equation}
\label{eqn:2nd-deri-(0,1)--(0,1/3)}
\frac
{
\theta^{\prime \prime}
\left[
\begin{array}{c}
1 \\
\frac13
\end{array}
\right]
}
{
\theta
\left[
\begin{array}{c}
1 \\
\frac13
\end{array}
\right]
}
+
\frac
{
\theta^{\prime \prime}
\left[
\begin{array}{c}
0 \\
\frac13
\end{array}
\right]
}
{
\theta
\left[
\begin{array}{c}
0 \\
\frac13
\end{array}
\right]
}
+
\frac
{
\theta^{\prime \prime}
\left[
\begin{array}{c}
0 \\
1
\end{array}
\right]
}
{
\theta
\left[
\begin{array}{c}
0 \\
1
\end{array}
\right]
}
-
\frac
{
\theta^{\prime \prime \prime}
\left[
\begin{array}{c}
1 \\
1
\end{array}
\right]
}
{
\theta^{\prime}
\left[
\begin{array}{c}
1\\
1
\end{array}
\right]
}
-
2
\frac
{
\theta^{\prime}
\left[
\begin{array}{c}
1\\
\frac13
\end{array}
\right]
}
{
\theta
\left[
\begin{array}{c}
1\\
\frac13
\end{array}
\right]
}
\frac
{
\theta^{\prime}
\left[
\begin{array}{c}
0\\
\frac13
\end{array}
\right]
}
{
\theta
\left[
\begin{array}{c}
0\\
\frac13
\end{array}
\right]
}
=0,
\end{equation}
\begin{equation}
\label{eqn:2nd-deri-(0,0)--(0,2/3)}
\frac
{
\theta^{\prime \prime}
\left[
\begin{array}{c}
1 \\
\frac13
\end{array}
\right]
}
{
\theta
\left[
\begin{array}{c}
1 \\
\frac13
\end{array}
\right]
}
+
\frac
{
\theta^{\prime \prime}
\left[
\begin{array}{c}
0 \\
\frac23
\end{array}
\right]
}
{
\theta
\left[
\begin{array}{c}
0 \\
\frac23
\end{array}
\right]
}
+
\frac
{
\theta^{\prime \prime}
\left[
\begin{array}{c}
0 \\
0
\end{array}
\right]
}
{
\theta
\left[
\begin{array}{c}
0 \\
0
\end{array}
\right]
}
-
\frac
{
\theta^{\prime \prime \prime}
\left[
\begin{array}{c}
1 \\
1
\end{array}
\right]
}
{
\theta^{\prime}
\left[
\begin{array}{c}
1\\
1
\end{array}
\right]
}
+
2
\frac
{
\theta^{\prime}
\left[
\begin{array}{c}
1\\
\frac13
\end{array}
\right]
}
{
\theta
\left[
\begin{array}{c}
1\\
\frac13
\end{array}
\right]
}
\frac
{
\theta^{\prime}
\left[
\begin{array}{c}
0\\
\frac23
\end{array}
\right]
}
{
\theta
\left[
\begin{array}{c}
0\\
\frac23
\end{array}
\right]
}
=0,
\end{equation}
\begin{equation}
\label{eqn:2nd-deri-(0,1/3)--(0,2/3)}
\frac
{
\theta^{\prime \prime}
\left[
\begin{array}{c}
0 \\
\frac13
\end{array}
\right]
}
{
\theta
\left[
\begin{array}{c}
0 \\
\frac13
\end{array}
\right]
}
+
\frac
{
\theta^{\prime \prime}
\left[
\begin{array}{c}
0 \\
\frac23
\end{array}
\right]
}
{
\theta
\left[
\begin{array}{c}
0 \\
\frac23
\end{array}
\right]
}
+
\frac
{
\theta^{\prime \prime}
\left[
\begin{array}{c}
1 \\
0
\end{array}
\right]
}
{
\theta
\left[
\begin{array}{c}
1 \\
0
\end{array}
\right]
}
-
\frac
{
\theta^{\prime \prime \prime}
\left[
\begin{array}{c}
1 \\
1
\end{array}
\right]
}
{
\theta^{\prime}
\left[
\begin{array}{c}
1\\
1
\end{array}
\right]
}
+
2
\frac
{
\theta^{\prime}
\left[
\begin{array}{c}
0\\
\frac13
\end{array}
\right]
}
{
\theta
\left[
\begin{array}{c}
0\\
\frac13
\end{array}
\right]
}
\frac
{
\theta^{\prime}
\left[
\begin{array}{c}
0\\
\frac23
\end{array}
\right]
}
{
\theta
\left[
\begin{array}{c}
0\\
\frac23
\end{array}
\right]
}
=0.
\end{equation}
}
\end{proposition}

\begin{proof}
Consider the following elliptic functions: 
\begin{equation*}
\frac
{
\theta^3
\left[
\begin{array}{c}
1 \\
\frac13
\end{array}
\right]
(z)
}
{
\theta^3
\left[
\begin{array}{c}
1 \\
1
\end{array}
\right]
(z)
},   \,\,
\frac
{
\theta^2
\left[
\begin{array}{c}
0 \\
\frac13
\end{array}
\right]
(z)
\theta
\left[
\begin{array}{c}
1 \\
\frac13
\end{array}
\right]
(z)
}
{
\theta^3
\left[
\begin{array}{c}
1 \\
1
\end{array}
\right]
(z)
}, \quad
\frac
{
\theta^2
\left[
\begin{array}{c}
0 \\
\frac23
\end{array}
\right]
(z)
\theta
\left[
\begin{array}{c}
1 \\
-\frac13
\end{array}
\right]
(z)
}
{
\theta^3
\left[
\begin{array}{c}
1 \\
1
\end{array}
\right]
(z)
}
\end{equation*}
\begin{equation*}
\frac
{
\theta
\left[
\begin{array}{c}
0 \\
\frac13
\end{array}
\right]
(z)
\theta
\left[
\begin{array}{c}
1 \\
-\frac13
\end{array}
\right]
(z)
\theta
\left[
\begin{array}{c}
0 \\
1
\end{array}
\right]
(z)
}
{
\theta^3
\left[
\begin{array}{c}
1 \\
1
\end{array}
\right]
(z)
}, \quad
\frac
{
\theta
\left[
\begin{array}{c}
0 \\
\frac23
\end{array}
\right]
(z)
\theta
\left[
\begin{array}{c}
1 \\
\frac13
\end{array}
\right]
(z)
\theta
\left[
\begin{array}{c}
0 \\
0
\end{array}
\right]
(z)
}
{
\theta^3
\left[
\begin{array}{c}
1 \\
1
\end{array}
\right]
(z)
},
\end{equation*}
and
\begin{equation*}
\frac
{
\theta
\left[
\begin{array}{c}
0 \\
\frac13
\end{array}
\right]
(z)
\theta
\left[
\begin{array}{c}
0 \\
\frac23
\end{array}
\right]
(z)
\theta
\left[
\begin{array}{c}
1 \\
0
\end{array}
\right]
(z)
}
{
\theta^3
\left[
\begin{array}{c}
1 \\
1
\end{array}
\right]
(z)
}.
\end{equation*}
\end{proof}

\begin{proposition}
\label{prop:3rd-derivative-(0,1/3)-(0,2/3)-(1,1/3)}
{\it
For every $\tau\in\mathbb{H}^2,$ we have 
\begin{equation*}
\label{eqn:3rd-order-(0,1/3)}
\frac
{
\theta^{(3)}
\left[
\begin{array}{c}
0 \\
\frac13
\end{array}
\right]
}
{
\theta
\left[
\begin{array}{c}
0 \\
\frac13
\end{array}
\right]
}
+
6
\frac
{
\theta^{(2)}
\left[
\begin{array}{c}
0 \\
\frac13
\end{array}
\right]
}
{
\theta
\left[
\begin{array}{c}
0 \\
\frac13
\end{array}
\right]
}
\frac
{
\theta^{(1)}
\left[
\begin{array}{c}
0 \\
\frac13
\end{array}
\right]
}
{
\theta
\left[
\begin{array}{c}
0 \\
\frac13
\end{array}
\right]
}
+
2
\left\{
\frac
{
\theta^{(1)}
\left[
\begin{array}{c}
0 \\
\frac13
\end{array}
\right]
}
{
\theta
\left[
\begin{array}{c}
0 \\
\frac13
\end{array}
\right]
}
\right\}^3
+
3
\frac
{
\theta^{(1)}
\left[
\begin{array}{c}
0 \\
\frac13
\end{array}
\right]
}
{
\theta
\left[
\begin{array}{c}
0 \\
\frac13
\end{array}
\right]
}
\frac
{
\theta^{(2)}
\left[
\begin{array}{c}
0 \\
1
\end{array}
\right]
}
{
\theta
\left[
\begin{array}{c}
0 \\
1
\end{array}
\right]
}
-
4
\frac
{
\theta^{(1)}
\left[
\begin{array}{c}
0 \\
\frac13
\end{array}
\right]
}
{
\theta
\left[
\begin{array}{c}
0 \\
\frac13
\end{array}
\right]
}
\frac
{
\theta^{(3)}
\left[
\begin{array}{c}
1 \\
1
\end{array}
\right]
}
{
\theta^{(1)}
\left[
\begin{array}{c}
1 \\
1
\end{array}
\right]
}
=0,
\end{equation*}
\begin{equation*}
\label{eqn:3rd-order-(0,2/3)}
\frac
{
\theta^{(3)}
\left[
\begin{array}{c}
0 \\
\frac23
\end{array}
\right]
}
{
\theta
\left[
\begin{array}{c}
0 \\
\frac23
\end{array}
\right]
}
+
6
\frac
{
\theta^{(2)}
\left[
\begin{array}{c}
0 \\
\frac23
\end{array}
\right]
}
{
\theta
\left[
\begin{array}{c}
0 \\
\frac23
\end{array}
\right]
}
\frac
{
\theta^{(1)}
\left[
\begin{array}{c}
0 \\
\frac13
\end{array}
\right]
}
{
\theta
\left[
\begin{array}{c}
0 \\
\frac13
\end{array}
\right]
}
+
2
\left\{
\frac
{
\theta^{(1)}
\left[
\begin{array}{c}
0 \\
\frac23
\end{array}
\right]
}
{
\theta
\left[
\begin{array}{c}
0 \\
\frac23
\end{array}
\right]
}
\right\}^3
+
3
\frac
{
\theta^{(1)}
\left[
\begin{array}{c}
0 \\
\frac23
\end{array}
\right]
}
{
\theta
\left[
\begin{array}{c}
0 \\
\frac23
\end{array}
\right]
}
\frac
{
\theta^{(2)}
\left[
\begin{array}{c}
0 \\
0
\end{array}
\right]
}
{
\theta
\left[
\begin{array}{c}
0 \\
0
\end{array}
\right]
}
-
4
\frac
{
\theta^{(1)}
\left[
\begin{array}{c}
0 \\
\frac23
\end{array}
\right]
}
{
\theta
\left[
\begin{array}{c}
0 \\
\frac23
\end{array}
\right]
}
\frac
{
\theta^{(3)}
\left[
\begin{array}{c}
1 \\
1
\end{array}
\right]
}
{
\theta^{(1)}
\left[
\begin{array}{c}
1 \\
1
\end{array}
\right]
}
=0,
\end{equation*}
\begin{align*}
&
\frac
{
\theta^{(3)}
\left[
\begin{array}{c}
0 \\
\frac13
\end{array}
\right]
}
{
\theta
\left[
\begin{array}{c}
0 \\
\frac13
\end{array}
\right]
}
-
\frac
{
\theta^{(3)}
\left[
\begin{array}{c}
1 \\
\frac13
\end{array}
\right]
}
{
\theta
\left[
\begin{array}{c}
1 \\
\frac13
\end{array}
\right]
}
+
3
\frac
{
\theta^{(2)}
\left[
\begin{array}{c}
0 \\
\frac13
\end{array}
\right]
}
{
\theta
\left[
\begin{array}{c}
0 \\
\frac13
\end{array}
\right]
}
\frac
{
\theta^{(1)}
\left[
\begin{array}{c}
0 \\
\frac13
\end{array}
\right]
}
{
\theta
\left[
\begin{array}{c}
0 \\
\frac13
\end{array}
\right]
}
-
6
\frac
{
\theta^{(2)}
\left[
\begin{array}{c}
0 \\
\frac13
\end{array}
\right]
}
{
\theta
\left[
\begin{array}{c}
0 \\
\frac13
\end{array}
\right]
}
\frac
{
\theta^{(1)}
\left[
\begin{array}{c}
1 \\
\frac13
\end{array}
\right]
}
{
\theta
\left[
\begin{array}{c}
1 \\
\frac13
\end{array}
\right]
}
-
6
\frac
{
\theta^{(1)}
\left[
\begin{array}{c}
1 \\
\frac13
\end{array}
\right]
}
{
\theta
\left[
\begin{array}{c}
1 \\
\frac13
\end{array}
\right]
}
\left\{
\frac
{
\theta^{(1)}
\left[
\begin{array}{c}
0 \\
\frac13
\end{array}
\right]
}
{
\theta
\left[
\begin{array}{c}
0 \\
\frac13
\end{array}
\right]
}
\right\}^2  \notag \\
&
+6
\frac
{
\theta^{(2)}
\left[
\begin{array}{c}
1 \\
\frac13
\end{array}
\right]
}
{
\theta
\left[
\begin{array}{c}
1 \\
\frac13
\end{array}
\right]
}
\frac
{
\theta^{(1)}
\left[
\begin{array}{c}
0 \\
\frac13
\end{array}
\right]
}
{
\theta
\left[
\begin{array}{c}
0 \\
\frac13
\end{array}
\right]
}
+
6
\frac
{
\theta^{(1)}
\left[
\begin{array}{c}
0 \\
\frac13
\end{array}
\right]
}
{
\theta
\left[
\begin{array}{c}
0 \\
\frac13
\end{array}
\right]
}
\left\{
\frac
{
\theta^{(1)}
\left[
\begin{array}{c}
1 \\
\frac13
\end{array}
\right]
}
{
\theta
\left[
\begin{array}{c}
1 \\
\frac13
\end{array}
\right]
}
\right\}^2
-
3
\frac
{
\theta^{(2)}
\left[
\begin{array}{c}
1 \\
\frac13
\end{array}
\right]
}
{
\theta
\left[
\begin{array}{c}
1 \\
\frac13
\end{array}
\right]
}
\frac
{
\theta^{(1)}
\left[
\begin{array}{c}
1 \\
\frac13
\end{array}
\right]
}
{
\theta
\left[
\begin{array}{c}
1 \\
\frac13
\end{array}
\right]
}  \\
&
-
4
\frac
{
\theta^{(3)}
\left[
\begin{array}{c}
1 \\
1
\end{array}
\right]
}
{
\theta^{(1)}
\left[
\begin{array}{c}
1 \\
1
\end{array}
\right]
}
\left\{
-
\frac
{
\theta^{(1)}
\left[
\begin{array}{c}
1 \\
\frac13
\end{array}
\right]
}
{
\theta
\left[
\begin{array}{c}
1 \\
\frac13
\end{array}
\right]
} 
+
\frac
{
\theta^{(1)}
\left[
\begin{array}{c}
0 \\
\frac13
\end{array}
\right]
}
{
\theta
\left[
\begin{array}{c}
0 \\
\frac13
\end{array}
\right]
} 
\right\}
=0,
\end{align*}
and
\begin{align*}
&
\frac
{
\theta^{(3)}
\left[
\begin{array}{c}
0 \\
\frac23
\end{array}
\right]
}
{
\theta
\left[
\begin{array}{c}
0 \\
\frac23
\end{array}
\right]
}
+
\frac
{
\theta^{(3)}
\left[
\begin{array}{c}
1 \\
\frac13
\end{array}
\right]
}
{
\theta
\left[
\begin{array}{c}
1 \\
\frac13
\end{array}
\right]
}
+
3
\frac
{
\theta^{(2)}
\left[
\begin{array}{c}
0 \\
\frac23
\end{array}
\right]
}
{
\theta
\left[
\begin{array}{c}
0 \\
\frac23
\end{array}
\right]
}
\frac
{
\theta^{(1)}
\left[
\begin{array}{c}
0 \\
\frac13
\end{array}
\right]
}
{
\theta
\left[
\begin{array}{c}
0 \\
\frac13
\end{array}
\right]
}
+
6
\frac
{
\theta^{(2)}
\left[
\begin{array}{c}
0 \\
\frac23
\end{array}
\right]
}
{
\theta
\left[
\begin{array}{c}
0 \\
\frac23
\end{array}
\right]
}
\frac
{
\theta^{(1)}
\left[
\begin{array}{c}
1 \\
\frac13
\end{array}
\right]
}
{
\theta
\left[
\begin{array}{c}
1 \\
\frac13
\end{array}
\right]
}
+
6
\frac
{
\theta^{(1)}
\left[
\begin{array}{c}
1 \\
\frac13
\end{array}
\right]
}
{
\theta
\left[
\begin{array}{c}
1 \\
\frac13
\end{array}
\right]
}
\left\{
\frac
{
\theta^{(1)}
\left[
\begin{array}{c}
0 \\
\frac23
\end{array}
\right]
}
{
\theta
\left[
\begin{array}{c}
0 \\
\frac23
\end{array}
\right]
}
\right\}^2  \notag \\
&
+6
\frac
{
\theta^{(2)}
\left[
\begin{array}{c}
1 \\
\frac13
\end{array}
\right]
}
{
\theta
\left[
\begin{array}{c}
1 \\
\frac13
\end{array}
\right]
}
\frac
{
\theta^{(1)}
\left[
\begin{array}{c}
0 \\
\frac23
\end{array}
\right]
}
{
\theta
\left[
\begin{array}{c}
0 \\
\frac23
\end{array}
\right]
}
+
6
\frac
{
\theta^{(1)}
\left[
\begin{array}{c}
0 \\
\frac23
\end{array}
\right]
}
{
\theta
\left[
\begin{array}{c}
0 \\
\frac23
\end{array}
\right]
}
\left\{
\frac
{
\theta^{(1)}
\left[
\begin{array}{c}
1 \\
\frac13
\end{array}
\right]
}
{
\theta
\left[
\begin{array}{c}
1 \\
\frac13
\end{array}
\right]
}
\right\}^2
+
3
\frac
{
\theta^{(2)}
\left[
\begin{array}{c}
1 \\
\frac13
\end{array}
\right]
}
{
\theta
\left[
\begin{array}{c}
1 \\
\frac13
\end{array}
\right]
}
\frac
{
\theta^{(1)}
\left[
\begin{array}{c}
1 \\
\frac13
\end{array}
\right]
}
{
\theta
\left[
\begin{array}{c}
1 \\
\frac13
\end{array}
\right]
}  \\
&
-
4
\frac
{
\theta^{(3)}
\left[
\begin{array}{c}
1 \\
1
\end{array}
\right]
}
{
\theta^{(1)}
\left[
\begin{array}{c}
1 \\
1
\end{array}
\right]
}
\left\{
\frac
{
\theta^{(1)}
\left[
\begin{array}{c}
1 \\
\frac13
\end{array}
\right]
}
{
\theta
\left[
\begin{array}{c}
1 \\
\frac13
\end{array}
\right]
} 
+
\frac
{
\theta^{(1)}
\left[
\begin{array}{c}
0 \\
\frac23
\end{array}
\right]
}
{
\theta
\left[
\begin{array}{c}
0 \\
\frac23
\end{array}
\right]
} 
\right\}
=0.
\end{align*}
}
\end{proposition}

\begin{proof}
Consider the following elliptic functions:
\begin{equation*}
\frac
{
\theta^3
\left[
\begin{array}{c}
0 \\
\frac13
\end{array}
\right](z)
\theta
\left[
\begin{array}{c}
0 \\
1
\end{array}
\right](z)
}
{
\theta^4
\left[
\begin{array}{c}
1 \\
1
\end{array}
\right](z)
}, \quad
\frac
{
\theta^3
\left[
\begin{array}{c}
0 \\
\frac23
\end{array}
\right](z)
\theta
\left[
\begin{array}{c}
0 \\
0
\end{array}
\right](z)
}
{
\theta^4
\left[
\begin{array}{c}
1 \\
1
\end{array}
\right](z)
},
\,\,
\frac
{
\theta^2
\left[
\begin{array}{c}
1\\
-\frac13
\end{array}
\right](z)
\theta^2
\left[
\begin{array}{c}
0\\
\frac13
\end{array}
\right](z)
}
{
\theta^4
\left[
\begin{array}{c}
1 \\
1
\end{array}
\right](z)
}, 
\quad
\frac
{
\theta^2
\left[
\begin{array}{c}
1\\
\frac13
\end{array}
\right](z)
\theta^2
\left[
\begin{array}{c}
0\\
\frac23
\end{array}
\right](z)
}
{
\theta^4
\left[
\begin{array}{c}
1 \\
1
\end{array}
\right](z)
}. 
\end{equation*}
\end{proof}

\begin{proposition}
\label{prop:4th-derivative-(0,1/3)-(0,2/3)-(1,1/3)}
{\it 
For every $\tau\in\mathbb{H}^2,$ we have 
\begin{align*}
&
4
\frac
{
\theta^{(4)}
\left[
\begin{array}{c}
0 \\
\frac13
\end{array}
\right]
}
{
\theta
\left[
\begin{array}{c}
0 \\
\frac13
\end{array}
\right]
}
+
\frac
{
\theta^{(4)}
\left[
\begin{array}{c}
1 \\
\frac13
\end{array}
\right]
}
{
\theta
\left[
\begin{array}{c}
1 \\
\frac13
\end{array}
\right]
}
-
\frac
{
\theta^{(5)}
\left[
\begin{array}{c}
1 \\
1
\end{array}
\right]
}
{
\theta^{(1)}
\left[
\begin{array}{c}
1 \\
1
\end{array}
\right]
}   
+
48
\frac
{
\theta^{(1)}
\left[
\begin{array}{c}
0 \\
\frac13
\end{array}
\right]
}
{
\theta
\left[
\begin{array}{c}
0 \\
\frac13
\end{array}
\right]
}
\frac
{
\theta^{(3)}
\left[
\begin{array}{c}
0 \\
\frac13
\end{array}
\right]
}
{
\theta
\left[
\begin{array}{c}
0 \\
\frac13
\end{array}
\right]
}  
-
16
\frac
{
\theta^{(1)}
\left[
\begin{array}{c}
1 \\
\frac13
\end{array}
\right]
}
{
\theta
\left[
\begin{array}{c}
1 \\
\frac13
\end{array}
\right]
}
\frac
{
\theta^{(3)}
\left[
\begin{array}{c}
0 \\
\frac13
\end{array}
\right]
}
{
\theta
\left[
\begin{array}{c}
0 \\
\frac13
\end{array}
\right]
}   \\
&-
16
\frac
{
\theta^{(1)}
\left[
\begin{array}{c}
0 \\
\frac13
\end{array}
\right]
}
{
\theta
\left[
\begin{array}{c}
0 \\
\frac13
\end{array}
\right]
}
\frac
{
\theta^{(3)}
\left[
\begin{array}{c}
1 \\
\frac13
\end{array}
\right]
}
{
\theta
\left[
\begin{array}{c}
1 \\
\frac13
\end{array}
\right]
}
+
36
\left\{
\frac
{
\theta^{(2)}
\left[
\begin{array}{c}
0 \\
\frac13
\end{array}
\right]
}
{
\theta
\left[
\begin{array}{c}
0 \\
\frac13
\end{array}
\right]
}
\right\}^2  
+144
\left\{
\frac
{
\theta^{(1)}
\left[
\begin{array}{c}
0 \\
\frac13
\end{array}
\right]
}
{
\theta
\left[
\begin{array}{c}
0 \\
\frac13
\end{array}
\right]
}
\right\}^2  
\frac
{
\theta^{(2)}
\left[
\begin{array}{c}
0 \\
\frac13
\end{array}
\right]
}
{
\theta
\left[
\begin{array}{c}
0 \\
\frac13
\end{array}
\right]
}  \\
&-144
\frac
{
\theta^{(1)}
\left[
\begin{array}{c}
1 \\
\frac13
\end{array}
\right]
}
{
\theta
\left[
\begin{array}{c}
1 \\
\frac13
\end{array}
\right]
}
\frac
{
\theta^{(1)}
\left[
\begin{array}{c}
0 \\
\frac13
\end{array}
\right]
}
{
\theta
\left[
\begin{array}{c}
0 \\
\frac13
\end{array}
\right]
}
\frac
{
\theta^{(2)}
\left[
\begin{array}{c}
0 \\
\frac13
\end{array}
\right]
}
{
\theta
\left[
\begin{array}{c}
0 \\
\frac13
\end{array}
\right]
}  
+24
\frac
{
\theta^{(2)}
\left[
\begin{array}{c}
1 \\
\frac13
\end{array}
\right]
}
{
\theta
\left[
\begin{array}{c}
1 \\
\frac13
\end{array}
\right]
}
\frac
{
\theta^{(2)}
\left[
\begin{array}{c}
0 \\
\frac13
\end{array}
\right]
}
{
\theta
\left[
\begin{array}{c}
0 \\
\frac13
\end{array}
\right]
}  
+72
\left\{
\frac
{
\theta^{(1)}
\left[
\begin{array}{c}
0 \\
\frac13
\end{array}
\right]
}
{
\theta
\left[
\begin{array}{c}
0 \\
\frac13
\end{array}
\right]
}  
\right\}^2
\frac
{
\theta^{(2)}
\left[
\begin{array}{c}
1 \\
\frac13
\end{array}
\right]
}
{
\theta
\left[
\begin{array}{c}
1 \\
\frac13
\end{array}
\right]
}  \\
&-96
\frac
{
\theta^{(1)}
\left[
\begin{array}{c}
1 \\
\frac13
\end{array}
\right]
}
{
\theta
\left[
\begin{array}{c}
1 \\
\frac13
\end{array}
\right]
}  
\left\{
\frac
{
\theta^{(1)}
\left[
\begin{array}{c}
0 \\
\frac13
\end{array}
\right]
}
{
\theta
\left[
\begin{array}{c}
0 \\
\frac13
\end{array}
\right]
}  
\right\}^3
+
24
\left\{
\frac
{
\theta^{(1)}
\left[
\begin{array}{c}
0 \\
\frac13
\end{array}
\right]
}
{
\theta
\left[
\begin{array}{c}
0 \\
\frac13
\end{array}
\right]
}  
\right\}^4  
+10
\left\{
\frac
{
\theta^{(3)}
\left[
\begin{array}{c}
1 \\
1
\end{array}
\right]
}
{
\theta^{(1)}
\left[
\begin{array}{c}
1 \\
1
\end{array}
\right]
}
\right\}^2  \\
&
-10
\frac
{
\theta^{(3)}
\left[
\begin{array}{c}
1 \\
1
\end{array}
\right]
}
{
\theta^{(1)}
\left[
\begin{array}{c}
1 \\
1
\end{array}
\right]
}
\left[
\frac
{
\theta^{(2)}
\left[
\begin{array}{c}
1 \\
\frac13
\end{array}
\right]
}
{
\theta
\left[
\begin{array}{c}
1 \\
\frac13
\end{array}
\right]
}
-
8
\frac
{
\theta^{(1)}
\left[
\begin{array}{c}
1 \\
\frac13
\end{array}
\right]
}
{
\theta
\left[
\begin{array}{c}
1 \\
\frac13
\end{array}
\right]
}
\frac
{
\theta^{(1)}
\left[
\begin{array}{c}
0 \\
\frac13
\end{array}
\right]
}
{
\theta
\left[
\begin{array}{c}
0 \\
\frac13
\end{array}
\right]
}
+4
\frac
{
\theta^{(2)}
\left[
\begin{array}{c}
0 \\
\frac13
\end{array}
\right]
}
{
\theta
\left[
\begin{array}{c}
0 \\
\frac13
\end{array}
\right]
}
+
12
\left\{
\frac
{
\theta^{(1)}
\left[
\begin{array}{c}
0 \\
\frac13
\end{array}
\right]
}
{
\theta
\left[
\begin{array}{c}
0 \\
\frac13
\end{array}
\right]
}
\right\}^2
\right]
=0,
\end{align*}
and
\begin{align*}
&
4
\frac
{
\theta^{(4)}
\left[
\begin{array}{c}
0 \\
\frac23
\end{array}
\right]
}
{
\theta
\left[
\begin{array}{c}
0 \\
\frac23
\end{array}
\right]
}
+
\frac
{
\theta^{(4)}
\left[
\begin{array}{c}
1 \\
\frac13
\end{array}
\right]
}
{
\theta
\left[
\begin{array}{c}
1 \\
\frac13
\end{array}
\right]
}
-
\frac
{
\theta^{(5)}
\left[
\begin{array}{c}
1 \\
1
\end{array}
\right]
}
{
\theta^{(1)}
\left[
\begin{array}{c}
1 \\
1
\end{array}
\right]
}   
+
48
\frac
{
\theta^{(1)}
\left[
\begin{array}{c}
0 \\
\frac23
\end{array}
\right]
}
{
\theta
\left[
\begin{array}{c}
0 \\
\frac23
\end{array}
\right]
}
\frac
{
\theta^{(3)}
\left[
\begin{array}{c}
0 \\
\frac23
\end{array}
\right]
}
{
\theta
\left[
\begin{array}{c}
0 \\
\frac23
\end{array}
\right]
}  
+
16
\frac
{
\theta^{(1)}
\left[
\begin{array}{c}
1 \\
\frac13
\end{array}
\right]
}
{
\theta
\left[
\begin{array}{c}
1 \\
\frac13
\end{array}
\right]
}
\frac
{
\theta^{(3)}
\left[
\begin{array}{c}
0 \\
\frac23
\end{array}
\right]
}
{
\theta
\left[
\begin{array}{c}
0 \\
\frac23
\end{array}
\right]
}   \\
&+
16
\frac
{
\theta^{(1)}
\left[
\begin{array}{c}
0 \\
\frac23
\end{array}
\right]
}
{
\theta
\left[
\begin{array}{c}
0 \\
\frac23
\end{array}
\right]
}
\frac
{
\theta^{(3)}
\left[
\begin{array}{c}
1 \\
\frac13
\end{array}
\right]
}
{
\theta
\left[
\begin{array}{c}
1 \\
\frac13
\end{array}
\right]
}
+
36
\left\{
\frac
{
\theta^{(2)}
\left[
\begin{array}{c}
0 \\
\frac23
\end{array}
\right]
}
{
\theta
\left[
\begin{array}{c}
0 \\
\frac23
\end{array}
\right]
}
\right\}^2  
+144
\left\{
\frac
{
\theta^{(1)}
\left[
\begin{array}{c}
0 \\
\frac23
\end{array}
\right]
}
{
\theta
\left[
\begin{array}{c}
0 \\
\frac23
\end{array}
\right]
}
\right\}^2  
\frac
{
\theta^{(2)}
\left[
\begin{array}{c}
0 \\
\frac23
\end{array}
\right]
}
{
\theta
\left[
\begin{array}{c}
0 \\
\frac23
\end{array}
\right]
}  \\
&+144
\frac
{
\theta^{(1)}
\left[
\begin{array}{c}
1 \\
\frac13
\end{array}
\right]
}
{
\theta
\left[
\begin{array}{c}
1 \\
\frac13
\end{array}
\right]
}
\frac
{
\theta^{(1)}
\left[
\begin{array}{c}
0 \\
\frac23
\end{array}
\right]
}
{
\theta
\left[
\begin{array}{c}
0 \\
\frac23
\end{array}
\right]
}
\frac
{
\theta^{(2)}
\left[
\begin{array}{c}
0 \\
\frac23
\end{array}
\right]
}
{
\theta
\left[
\begin{array}{c}
0 \\
\frac23
\end{array}
\right]
}  
+24
\frac
{
\theta^{(2)}
\left[
\begin{array}{c}
1 \\
\frac13
\end{array}
\right]
}
{
\theta
\left[
\begin{array}{c}
1 \\
\frac13
\end{array}
\right]
}
\frac
{
\theta^{(2)}
\left[
\begin{array}{c}
0 \\
\frac23
\end{array}
\right]
}
{
\theta
\left[
\begin{array}{c}
0 \\
\frac23
\end{array}
\right]
}  
+72
\left\{
\frac
{
\theta^{(1)}
\left[
\begin{array}{c}
0 \\
\frac23
\end{array}
\right]
}
{
\theta
\left[
\begin{array}{c}
0 \\
\frac23
\end{array}
\right]
}  
\right\}^2
\frac
{
\theta^{(2)}
\left[
\begin{array}{c}
1 \\
\frac13
\end{array}
\right]
}
{
\theta
\left[
\begin{array}{c}
1 \\
\frac13
\end{array}
\right]
}  \\
&+96
\frac
{
\theta^{(1)}
\left[
\begin{array}{c}
1 \\
\frac13
\end{array}
\right]
}
{
\theta
\left[
\begin{array}{c}
1 \\
\frac13
\end{array}
\right]
}  
\left\{
\frac
{
\theta^{(1)}
\left[
\begin{array}{c}
0 \\
\frac23
\end{array}
\right]
}
{
\theta
\left[
\begin{array}{c}
0 \\
\frac23
\end{array}
\right]
}  
\right\}^3
+
24
\left\{
\frac
{
\theta^{(1)}
\left[
\begin{array}{c}
0 \\
\frac23
\end{array}
\right]
}
{
\theta
\left[
\begin{array}{c}
0 \\
\frac23
\end{array}
\right]
}  
\right\}^4  
+10
\left\{
\frac
{
\theta^{(3)}
\left[
\begin{array}{c}
1 \\
1
\end{array}
\right]
}
{
\theta^{(1)}
\left[
\begin{array}{c}
1 \\
1
\end{array}
\right]
}
\right\}^2  \\
&
-10
\frac
{
\theta^{(3)}
\left[
\begin{array}{c}
1 \\
1
\end{array}
\right]
}
{
\theta^{(1)}
\left[
\begin{array}{c}
1 \\
1
\end{array}
\right]
}
\left[
\frac
{
\theta^{(2)}
\left[
\begin{array}{c}
1 \\
\frac13
\end{array}
\right]
}
{
\theta
\left[
\begin{array}{c}
1 \\
\frac13
\end{array}
\right]
}
+
8
\frac
{
\theta^{(1)}
\left[
\begin{array}{c}
1 \\
\frac13
\end{array}
\right]
}
{
\theta
\left[
\begin{array}{c}
1 \\
\frac13
\end{array}
\right]
}
\frac
{
\theta^{(1)}
\left[
\begin{array}{c}
0 \\
\frac23
\end{array}
\right]
}
{
\theta
\left[
\begin{array}{c}
0 \\
\frac23
\end{array}
\right]
}
+4
\frac
{
\theta^{(2)}
\left[
\begin{array}{c}
0 \\
\frac23
\end{array}
\right]
}
{
\theta
\left[
\begin{array}{c}
0 \\
\frac23
\end{array}
\right]
}
+
12
\left\{
\frac
{
\theta^{(1)}
\left[
\begin{array}{c}
0 \\
\frac23
\end{array}
\right]
}
{
\theta
\left[
\begin{array}{c}
0 \\
\frac23
\end{array}
\right]
}
\right\}^2
\right]
=0.
\end{align*}
}
\end{proposition}

\begin{proof}
Consider the following elliptic functions:
\begin{equation*}
\frac
{
\theta^4
\left[
\begin{array}{c}
0 \\
\frac13
\end{array}
\right](z)
\theta
\left[
\begin{array}{c}
1 \\
-\frac13
\end{array}
\right](z)
}
{
\theta^5
\left[
\begin{array}{c}
1 \\
1
\end{array}
\right](z)
}, \quad
\frac
{
\theta^4
\left[
\begin{array}{c}
0 \\
\frac23
\end{array}
\right](z)
\theta
\left[
\begin{array}{c}
1 \\
\frac13
\end{array}
\right](z)
}
{
\theta^5
\left[
\begin{array}{c}
1 \\
1
\end{array}
\right](z)
}.
\end{equation*}
\end{proof}

\subsection{The ODEs for $Q(q)$ and $R(q)$        }

\begin{proof}
By Proposition \ref{prop:2nd-derivative-(0,1/3)-(0,2/3)}, 
we have 
\begin{equation*}
\frac
{
\theta^{\prime \prime}
\left[
\begin{array}{c}
0 \\
\frac13
\end{array}
\right]
}
{
\theta
\left[
\begin{array}{c}
0 \\
\frac13
\end{array}
\right]
}
=X^2-2XY-Y^2+\frac13 W, \,\,
\frac
{
\theta^{\prime \prime}
\left[
\begin{array}{c}
0 \\
\frac23
\end{array}
\right]
}
{
\theta
\left[
\begin{array}{c}
0 \\
\frac23
\end{array}
\right]
}
=X^2+2XZ-Z^2+\frac13 W, \,\,
\end{equation*}
\begin{equation*}
\frac
{
\theta^{\prime \prime}
\left[
\begin{array}{c}
0 \\
1
\end{array}
\right]
}
{
\theta
\left[
\begin{array}{c}
0 \\
1
\end{array}
\right]
}
=
X^2+4XY+Y^2+\frac13 W,  \,\,
\frac
{
\theta^{\prime \prime}
\left[
\begin{array}{c}
0 \\
0
\end{array}
\right]
}
{
\theta
\left[
\begin{array}{c}
0 \\
0
\end{array}
\right]
}
=
X^2-4 XZ +Z^2+\frac13 W,
\end{equation*}
and
\begin{equation*}
\frac
{
\theta^{\prime \prime}
\left[
\begin{array}{c}
1 \\
\frac13
\end{array}
\right]
}
{
\theta
\left[
\begin{array}{c}
1 \\
\frac13
\end{array}
\right]
}
=-2 X^2+\frac13 W.\,\,
\end{equation*}
\par
By Proposition \ref{prop:3rd-derivative-(0,1/3)-(0,2/3)-(1,1/3)}, 
we obtain 
\begin{equation*}
\frac
{
\theta^{\prime \prime \prime}
\left[
\begin{array}{c}
0 \\
\frac13
\end{array}
\right]
}
{
\theta
\left[
\begin{array}{c}
0 \\
\frac13
\end{array}
\right]
}
=
-9 X^2 Y+Y^3+YW, \,\,
\frac
{
\theta^{\prime \prime \prime}
\left[
\begin{array}{c}
0 \\
\frac23
\end{array}
\right]
}
{
\theta
\left[
\begin{array}{c}
0 \\
\frac23
\end{array}
\right]
}
=
-9X^2Z+Z^3+ZW.
\end{equation*}
\par
Since 
\begin{equation*}
\frac{
\theta^{\prime \prime \prime}
\left[
\begin{array}{c}
0 \\
\frac{j}{3}
\end{array}
\right]  
}
{ 
\theta
\left[
\begin{array}{c}
0 \\
\frac{j}{3}
\end{array}
\right]
}
=
4 \pi i
\frac{d}{d\tau}
\left\{
\frac{
\theta^{\prime }
\left[
\begin{array}{c}
0 \\
\frac{j}{3}
\end{array}
\right]  
}
{ 
\theta
\left[
\begin{array}{c}
0 \\
\frac{j}{3}
\end{array}
\right]
}
\right\}
+
\frac{
\theta^{\prime }
\left[
\begin{array}{c}
0 \\
\frac{j}{3}
\end{array}
\right]  
}
{ 
\theta
\left[
\begin{array}{c}
0 \\
\frac{j}{3}
\end{array}
\right]
}
\cdot
\frac{
\theta^{\prime \prime }
\left[
\begin{array}{c}
0\\
\frac{j}{3}
\end{array}
\right]  
}
{ 
\theta
\left[
\begin{array}{c}
0 \\
\frac{j}{3}
\end{array}
\right]
},  \,\,(j=1,2),
\end{equation*}
it follows that 
\begin{equation*}
4 \pi i Y^{\prime}
=-10 X^2 Y+2 X Y^2+2 Y^3+\frac23 YW, \,\,
4 \pi i Z^{\prime}
=-10 X^2 Z-2 X Z^2+2 Z^3+\frac23 ZW,
\end{equation*}
where ${ }^{\prime}=d/d\tau.$ Considering $d/d\tau=2 \pi i q d/dq,$ we have the ODEs for $Q(q)$ and $R(q).$ 
\end{proof}

\subsection{The ODE for $P(q)$}

\begin{proof}
\par
By Proposition \ref{prop:3rd-derivative-(0,1/3)-(0,2/3)-(1,1/3)}, 
we have
\begin{equation*}
\frac
{
\theta^{\prime \prime \prime}
\left[
\begin{array}{c}
1 \\
\frac13
\end{array}
\right]
}
{
\theta
\left[
\begin{array}{c}
1 \\
\frac13
\end{array}
\right]
}
=-6XY^2-2 Y^3+X W=-6X Z^2+2 Z^3+ XW,
\end{equation*}
which implies that 
\begin{equation}
\label{eqn:relation-(X,Y,Z)}
3 X(Y-Z)+Y^2-YZ +Z^2=0.
\end{equation}
The relation (\ref{eqn:(P,Q,R)-relation}) follows from this relation. 
\par
Since 
\begin{equation*}
\frac{
\theta^{\prime \prime \prime}
\left[
\begin{array}{c}
1 \\
\frac13
\end{array}
\right]  
}
{ 
\theta
\left[
\begin{array}{c}
1 \\
\frac13
\end{array}
\right]
}
=
4 \pi i
\frac{d}{d\tau}
\left\{
\frac{
\theta^{\prime }
\left[
\begin{array}{c}
1 \\
\frac13
\end{array}
\right]  
}
{ 
\theta
\left[
\begin{array}{c}
1 \\
\frac13
\end{array}
\right]
}
\right\}
+
\frac{
\theta^{\prime }
\left[
\begin{array}{c}
1 \\
\frac13
\end{array}
\right]  
}
{ 
\theta
\left[
\begin{array}{c}
1 \\
\frac13
\end{array}
\right]
}
\cdot
\frac{
\theta^{\prime \prime }
\left[
\begin{array}{c}
1 \\
\frac13
\end{array}
\right]  
}
{ 
\theta
\left[
\begin{array}{c}
1 \\
\frac13
\end{array}
\right]
}, 
\end{equation*}
it follows that 
\begin{equation*}
4 \pi i X^{\prime}=2X^3-6X Y^2-2 Y^3+\frac23 XW=2X^3-6XZ^2+2Z^3+\frac23 XW, \,\,(\,\,{ }^{\prime}=d/d\tau).
\end{equation*}
\end{proof}

\subsection{The ODE for $S(q)$}

\begin{proof}
We first note that 
\begin{equation*}
\frac{
\theta^{(4)}
\left[
\begin{array}{c}
0 \\
\frac{j}{3}
\end{array}
\right]  
}
{ 
\theta
\left[
\begin{array}{c}
0 \\
\frac{j}{3}
\end{array}
\right]
}
=
4 \pi i
\frac{d}{d\tau}
\left\{
\frac{
\theta^{\prime \prime }
\left[
\begin{array}{c}
0 \\
\frac{j}{3}
\end{array}
\right]  
}
{ 
\theta
\left[
\begin{array}{c}
0 \\
\frac{j}{3}
\end{array}
\right]
}
\right\}
+
\left\{
\frac{
\theta^{\prime \prime  }
\left[
\begin{array}{c}
0 \\
\frac{j}{3}
\end{array}
\right]  
}
{ 
\theta
\left[
\begin{array}{c}
0 \\
\frac{j}{3}
\end{array}
\right]
}
\right\}^2,   \,\,(j=1,2),
\end{equation*}
which implies that 
\begin{equation*}
\frac{
\theta^{(4)}
\left[
\begin{array}{c}
0 \\
\frac13
\end{array}
\right]  
}
{ 
\theta
\left[
\begin{array}{c}
0 \\
\frac13
\end{array}
\right]
}
=
\frac{4\pi i}{3}W^{\prime}
+\frac{W^2}{9}+2 X^2W -4  X YW-2  Y^2W+5 X^4+12 X^3 Y+6 X^2 Y^2+4 X Y^3+Y^4
\end{equation*}
and 
\begin{equation*}
\frac{
\theta^{(4)}
\left[
\begin{array}{c}
0 \\
\frac23
\end{array}
\right]  
}
{ 
\theta
\left[
\begin{array}{c}
0 \\
\frac23
\end{array}
\right]
}
=
\frac{4\pi i}{3}W^{\prime}
+\frac{W^2}{9}+2 X^2W +4X Z W -2  Z^2W+5 X^4-12 X^3 Z+6 X^2 Z^2-4 X Z^3+Z^4, \,\,(\,{ }^{\prime}=d/d\tau).
\end{equation*}
In the same way, we obtain 
\begin{align*}
\frac{
\theta^{(4)}
\left[
\begin{array}{c}
1 \\
\frac13
\end{array}
\right]  
}
{ 
\theta
\left[
\begin{array}{c}
1 \\
\frac13
\end{array}
\right]
}
=&
\frac{4\pi i}{3}W^{\prime}
+\frac{W^2}{9}-4 X^2W -4 X^4+24 X^2 Y^2+8 X Y^3  \\
=&
\frac{4\pi i}{3}W^{\prime}
+\frac{W^2}{9}-4 X^2W -4 X^4+24 X^2 Z^2-8 X Z^3
\end{align*}
and 
\begin{equation*}
\frac{
\theta^{(5)}
\left[
\begin{array}{c}
1 \\
1
\end{array}
\right]  
}
{ 
\theta^{(1)}
\left[
\begin{array}{c}
1 \\
1
\end{array}
\right]
}=4 \pi i W^{\prime}+W^2.
\end{equation*}
\par
By Proposition \ref{prop:4th-derivative-(0,1/3)-(0,2/3)-(1,1/3)}, 
we obtain 
\begin{equation*}
4 \pi i W^{\prime}=\frac23 W^2-6X^4-36 X^2 Y^2-12 X Y^3
=\frac23 W^2-6X^4-36X^2 Z^2+12 XZ^3,
\end{equation*}
\begin{equation*}
\frac{
\theta^{(4)}
\left[
\begin{array}{c}
1 \\
\frac13
\end{array}
\right]  
}
{ 
\theta
\left[
\begin{array}{c}
1 \\
\frac13
\end{array}
\right]
}
=\frac{W^2}{3}-4 X^2W -6 X^4+12 X^2 Y^2+4 X Y^3
=\frac{W^2}{3}-4 X^2W -6 X^4+12 X^2 Z^2-4 X Z^3.
\end{equation*}
\end{proof}

\subsection{Note on $E_4(q)$ and $E_6(q)$ $\cdots$ (1)}

\begin{theorem}
\label{thm-E4-E6-level-6-(0,1/3)-(0,2/3)}
{\it
For $q\in\mathbb{C}$ with $|q|<1,$ we have 
\begin{align*}
E_4(q)=&
P^4+216 P^2 Q^2+432 P Q^3=P^4+216P^2R^2-432PR^3, \\
E_6(q)=&
P^6-540 P^4 Q^2-1080 P^3 Q^3-5832 P^2 Q^4-23328 P Q^5-23328 Q^6  \\
=&P^6-540 P^4 R^2+1080 P^3 R^3-5832 P^2 R^4+23328 P R^5-23328 R^6, 
\end{align*}
and
\begin{align*}
&q(q;q)_{\infty}^{24}
=
Q^2 (P-6 Q)^3 (P+2 Q) (P+3 Q)^6
=
R^2 (P-3 R)^6 (P-2 R) (P+6 R)^3,  \\
&j(q)=
\frac{P^3 \left(P^3+216 P Q^2+432 Q^3\right)^3}{Q^2 (P-6 Q)^3 (P+2 Q) (P+3 Q)^6}
=
\frac{P^3 \left(P^3+216 P R^2-432 R^3\right)^3}{R^2 (P-3 R)^6 (P-2 R) (P+6 R)^3}  \\
=&
\frac{(3 Q-T)^3 \left(243 Q^3-243 Q^2 T+9 Q T^2-T^3\right)^3}{Q^2 T^6 (Q-T) (9 Q-T)^3}  
=
\frac{(3 R+U)^3 \left(243 R^3+243 R^2 U+9 R U^2+U^3\right)^3}{R^2 U^6 (R+U) (9 R+U)^3}  \\
=&
\frac{(3 t-1)^3 \left(243 t^3-243 t^2+9 t-1\right)^3}{ t^2(t-1)  (9 t-1)^3}
=
\frac{(3 u+1)^3 \left(243 u^3+243 u^2+9 u+1\right)^3}{u^2 (u+1) (9 u+1)^3},
\end{align*}
where 
\begin{equation*}
t(q)=\frac{Q(q)}{T(q)}=q^{\frac12} 
\frac{
(q^{\frac12} ;q^{\frac12}  )_{\infty}^4   (q^3; q^3)_{\infty}^8
}
{
(q; q)_{\infty}^8 (q^{\frac32}; q^{\frac32})_{\infty}^4
}, \,\,
u(q)=\frac{R(q)}{U(q)}
=
q^{\frac12} 
\frac{
(q ;q )_{\infty}^4   (q^{\frac32}; q^{\frac32})_{\infty}^4   (q^6; q^6)_{\infty}^4
}
{
(q^{\frac12}; q^{\frac12})_{\infty}^4        (q^2 ;q^2 )_{\infty}^4  (q^3 ;q^3 )_{\infty}^4  
}. 
\end{equation*}
}
\end{theorem}

\begin{proof}
The theorem can be proved in the same way as Theorem \ref{thm:E_2-E_4-E_6-q-(1,2/3)}. 
Substitute $z=-1/3$ in 
\begin{equation*}
\wp=\wp(z;1,\tau)
=
\frac13
\frac
{
\theta^{\prime  \prime  \prime}
\left[
\begin{array}{c}
1 \\
1
\end{array}
\right]
}
{
\theta^{\prime  }
\left[
\begin{array}{c}
1 \\
1
\end{array}
\right]
}
-
\frac{d^2}{dz^2}
\log
\theta
\left[
\begin{array}{c}
1 \\
1
\end{array}
\right](z),
\end{equation*} 
and 
\begin{equation*}
(\wp^{\prime})^2=4\wp^3-g_2\wp-g_3.
\end{equation*}
\end{proof}

\begin{corollary}
{\it
For $q\in\mathbb{C}$ with $|q|<1,$ we have 
\begin{equation*}
E_4(q)=
\frac{16 \left(Q^2-Q R+R^2\right) \left(Q^6-3 Q^5 R+60 Q^4 R^2-115 Q^3 R^3+60 Q^2 R^4-3 Q R^5+R^6\right)}{(Q-R)^4},
\end{equation*}
\begin{align*}
E_6(q)=&
\frac{32 \left(Q^4+10 Q^3 R-12 Q^2 R^2+4 Q R^3-2 R^4\right) \left(2 Q^4-4 Q^3 R+12 Q^2 R^2-10 Q R^3-R^4\right) }{(Q-R)^6} \times \\
&\times \left(Q^4-14 Q^3 R+24 Q^2 R^2-14 Q R^3+R^4\right),
\end{align*}
\begin{equation*}
q(q;q)_{\infty}^{24}=
\frac{16 Q^2 R^2 (Q-2 R)^6 (2 Q-R)^6 (Q+R)^6}{(Q-R)^{10}},
\end{equation*}
and
\begin{align*}
j(q)
=&
\frac{256 \left(Q^2-Q R+R^2\right)^3 \left(Q^6-3 Q^5 R+60 Q^4 R^2-115 Q^3 R^3+60 Q^2 R^4-3 Q R^5+R^6\right)^3}{Q^2 R^2 (Q-R)^2(Q-2 R)^6  (2 Q-R)^6 (Q+R)^6}. 
\end{align*}
}
\end{corollary}

\begin{proof}
The relation (\ref{eqn:(P,Q,R)-relation}) implies that 
\begin{equation*}
P=-
\frac
{
2Q^2-2QR+2R^2
}
{
Q-R
}. 
\end{equation*}
The corollary follows from this relation. 
\end{proof}

\subsection{Note on $E_2(q)$} 

\begin{theorem}
\label{thm:E_2-(1,1/3)-(0,1/3)-(0,2/3)}
{\it
For $q\in\mathbb{C}$ with $|q|<1,$ we have 
\begin{align*}
P^2=&\frac12 \left(-E_2(q)+3 E_2(q^3) \right),  \\
Q^2=&\frac{1}{24}  \left(   -E_2(q)+2 E_2(q^{\frac32})-E_2(q^3)   \right),  \\
PQ=&\frac{1}{24} \left( -E_2(q^{\frac12})+4 E_2(q)-3 E_2(q^{\frac32})  \right), \\
R^2=&\frac{1}{24}  \left( -E_2(q)-2E_2(q^{\frac32})+11 E_2(q^3)-8E_2(q^6)        \right),   \\
PR=&\frac{1}{24}\left( -E_2(q^{\frac12})+2 E_2(q)-3 E_2(q^{\frac32})-4E_2(q^2)+18 E_2(q^3)-12 E_2(q^6) \right),  \\
QR=&\frac{1}{24} \left( -E_2(q)+2 E_2(q^2)+E_2(q^3)-2E_2(q^6)           \right),
\end{align*}
which implies that 
\begin{align*}
E_2(q^{\frac12})=&-P^2-24 PQ-36 Q^2+2 S, \\
E_2(q^{\frac32})=&\frac13
\left(
P^2+36 Q^2+2 S
\right),  \\
E_2(q^2)=&\frac12
\left(
P^2-6Q^2+24 Q R-6 R^2+S
\right),  \\
E_2(q^3)=&\frac13 \left(2 P^2+S \right),  \\
E_2(q^6)=&
\frac16(5P^2-18 Q^2-18 R^2+S). \\
\end{align*}
}
\end{theorem}

\begin{proof}
The formulas for $P^2, Q^2$ and $R^2$ can be obtained by applying Lemma \ref{lem:first-order-deri-square} to the characteristics 
$(\epsilon, \epsilon^{\prime})=(1,1/3), (0,1/3)$, and $(0,2/3).$ 
The other formulas follows from equations (\ref{eqn:2nd-deri-(0,1)--(0,1/3)}), (\ref{eqn:2nd-deri-(0,0)--(0,2/3)}), and (\ref{eqn:2nd-deri-(0,1/3)--(0,2/3)}) 
in Proposition \ref{prop:2nd-derivative-(0,1/3)-(0,2/3)}. 
\end{proof}

\begin{corollary}
\label{coro-4squares-triangular-E_2-(1)}
{\it
For $q\in\mathbb{C}$ with $|q|<1,$ 
we have 
\begin{align*}
\varphi^4(q^{\frac12})=&
P^2+8 P Q+8 Q^2+16 Q R-4 R^2
=
-\frac{4 Q (Q-2 R)^3}{(Q-R)^2}, \\
16q^{\frac12} \psi^4(q)=&
16 P Q+28 Q^2-16 Q R+4 R^2
=
-\frac{4 (Q+R)^3}{Q-R},  \\
\varphi^4(-q^{\frac12})=&
P^2-8 P Q-20 Q^2+32 Q R-8 R^2
=
\frac{4 R (2 Q-R)^3}{(Q-R)^2},
\end{align*}
where 
$
\varphi(q)=\sum_{n\in\mathbb{Z}} q^{n^2}$ and $\psi(q)=\sum_{n=0}^{\infty} q^{\frac{n(n+1)}{2}}.
$
}
\end{corollary}

\begin{proof}
From Matsuda \cite{Matsuda-2}, we recall 
\begin{align}
\varphi^4(q^{\frac12})=&\frac13\left(
-E_2(q^{\frac12})+4E_2(q^2)
\right),  \label{eqn:4squares-q}    \\
16q^{\frac12} \psi^4(q)=&
\frac13
\left(
-2E_2(q^{\frac12})+6E_2(q)-4E_2(q^2)
\right),     \label{eqn:4triangular-q}    \\
\varphi^4(-q^{\frac12})=&
\frac13
\left(
E_2(q^{\frac12})-6E_2(q)+8E_2(q^2)
\right),   \label{eqn:4squares-(-q)} 
\end{align}
which proves the corollary. 
\end{proof}

\begin{corollary}
\label{coro-4squares-triangular-E_2-(2)}
{\it
For $q\in\mathbb{C}$ with $|q|<1,$ 
we have 
\begin{align*}
\varphi^4(q^{\frac32})=&\frac13\left(
-E_2(q^{\frac32})+4E_2(q^6)
\right)
=
P^2-8 Q^2-4 R^2
=
-\frac{4 Q^3 (Q-2 R)}{(Q-R)^2}, \\
16q^{\frac32} \psi^4(q^3)=&
\frac13
\left(
-2E_2(q^{\frac32})+6E_2(q^3)-4E_2(q^6)
\right)
=
-4 Q^2+4 R^2,   \\
\varphi^4(-q^{\frac32})=&
\frac13
\left(
E_2(q^{\frac32})-6E_2(q^3)+8E_2(q^6)
\right)
=
P^2-4 Q^2-8 R^2
=
\frac{4 R^3 (2 Q-R)}{(Q-R)^2}. 
\end{align*}
}
\end{corollary}

\subsubsection*{Remark}
The formulas for $P^2, Q^2$ and $R^2$ were proved in \cite{Matsuda-2}.

\begin{theorem}
\label{thm:system-(A,B,C,D,E,F)}
{\it
For $q\in\mathbb{C}$ with $|q|<1,$ set 
\begin{equation*}
A=E_2(q^{\frac12}), \, B=E_2(q), \,C=E_2(q^2), \,D=E_2(q^{\frac32}), \, E=E_2(q^3), \, F=E_2(q^6).
\end{equation*}
Then, we have
\begin{align*}
q\frac{d}{dq}A
=&
\frac{-2 A^2 -3 B^2 +16 C^2    +9 A B-24 B C+4 A C }{18},     \\ 
q\frac{d}{dq}B
=&
\frac{-A^2   -9 B^2  -16 C^2    +6 A B  +24 B C-4 A C }{36},     \\
q\frac{d}{dq}C
=&
\frac{A^2  -3 B^2  -32 C^2-    6 A B  +36 B C+4 A C }{72},  \\
q\frac{d}{dq}D
=&
\frac{7 A^2+24 B^2-128 C^2-72 A B-32 A C+192 B C-81 D^2+90 A D}{432},  \\ 
q\frac{d}{dq}E
=&
\frac{2 A^2+9 B^2+32 C^2-12 A B-48 B C+8 A C-81 E^2+90 BE }{216},  \\
q\frac{d}{dq}F
=&
\frac{-A^2+3 B^2+14 C^2+6 A B-4 A C-36 B C-162 F^2+180 C F}{216},
\end{align*}
and
\begin{align}
&A^2-8 A B+6 A D+4 B^2+24 B E+9 D^2-72 D E+36 E^2=0, \label{eqn:relation-(1)-with-cubic-simple}  \\
&B^2-8 B C+6 B E+4 C^2+24 C F+9 E^2-72 E F+36 F^2=0,  \label{eqn:relation-(2)-with-cubic-simple}  \\
&A B-10 B C+2 A C+9 D E-90 E F+18 D F-3 A D+9 B E-6 C F+7 C^2+63 F^2=0, 
\label{eqn:relation-(3)-with-cubic-simple}  \\
&3 C^2+27 F^2+A D-3 A E+2 A F-6 B C-3 B D+21 B E-12 B F  \notag \\
&\hspace{60mm}+2 C D-12 C E+34 C F-54 E F=0,  
\label{eqn:relation-(4)-with-cubic-simple} \\
&C^2+D^2+3 E^2+F^2  -B C+2 B E-B F+C E-2 C F-6 D E+4 D F-3 E F=0. \label{eqn:relation-(5)-with-cubic-simple} 
\end{align}
}
\end{theorem}

\begin{proof}
The system of ODEs and proofs of the relations (\ref{eqn:relation-(1)-with-cubic-simple}),  (\ref{eqn:relation-(2)-with-cubic-simple}), and 
 (\ref{eqn:relation-(3)-with-cubic-simple}) can be found in Matsuda \cite{Matsuda-3}. 
Note that the relations (\ref{eqn:relation-(1)-with-cubic-simple}),  (\ref{eqn:relation-(2)-with-cubic-simple}) and 
(\ref{eqn:relation-(3)-with-cubic-simple}) are equivalent to 
\begin{align}
&7 A^2+24 B^2-128 C^2-72 A B+192 B C-32 A C \notag \\
&\hspace{15mm} +63 D^2+216 E^2-1152 F^2+90 A D-648 D E+1728 E F-288 D F=0,  \label{eqn:relation-(1)-with-cubic}  \\
&2 A^2+9 B^2+32 C^2-12 A B-48 B C+8 A C \notag \\
&\hspace{15mm} +18 D^2 +81 E^2 +288 F^2 +90 B E-108 D E-432 E F+72 D F=0,  \label{eqn:relation-(2)-with-cubic}  \\
&A^2-3 B^2 -14 C^2   -6 A B+36 B C+4 A C   \notag  \\
&\hspace{15mm} 
 +9 D^2-27 E^2-126 F^2 -180 C F-54  DE+324 E F+36 D F
=0.  
\label{eqn:relation-(3)-with-cubic} 
\end{align}
\par
Because $P^2R^2=(PR)^2$ and $Q^2R^2=(QR)^2$, Theorem \ref{thm:E_2-(1,1/3)-(0,1/3)-(0,2/3)} implies that 
\begin{align}
&A^2-8 B^2+16 C^2 +9 D^2-72 E^2+144 F^2-4 A B-16 B C+8 A C    -36 D E+72 D F-144 E F    \notag \\
&+6 A D-36 A E+24 A F   -36 B D+240 B E-144 B F+24 C D-144 C E+96 C F =0,  \label{eqn:relation-(4)-with-cubic}  \\
&C^2+D^2+3 E^2+F^2  -B C+2 B E-B F+C E-2 C F-6 D E+4 D F-3 E F=0.
\end{align}
\par
Considering 
$
-5/21\times \mathrm{Eq.}(\ref{eqn:relation-(1)-with-cubic})
+8/21\times  \mathrm{Eq.}(\ref{eqn:relation-(2)-with-cubic})
+40/21\times \mathrm{Eq.}(\ref{eqn:relation-(3)-with-cubic})
$ 
implies 
\begin{align*}
&A^2-8 B^2+16 C^2+9 D^2 -72 E^2 +144 F^2  \\
&+\frac{8 A B}{7}+\frac{128 A C}{7}-\frac{150 A D}{7}+\frac{32 B C}{7}+\frac{240 B E}{7}-\frac{2400 C F}{7}+\frac{72 D E}{7}+\frac{1152 D F}{7}+\frac{288 E F}{7}=0.
\end{align*}
Considering Eq. (\ref{eqn:relation-(4)-with-cubic}) shows 
\begin{align}
&3 A B+12 B C+6 A C  +27 D E+108 E F+54 D F  -16 A D+21 A E-14 A F  \notag \\
&\hspace{15mm}+21 B D-120 B E+84 B F-14 C D+84 C E-256 C F  =0.  \label{eqn:relation-(5)-with-cubic} 
\end{align}
$\mathrm{Eq.}(\ref{eqn:relation-(5)-with-cubic})-3 \times   \mathrm{Eq.}(\ref{eqn:relation-(3)-with-cubic-simple}) $ 
implies the relation (\ref{eqn:relation-(4)-with-cubic-simple}). 
\end{proof}

\subsubsection*{Remark}
Halphen's differential field is expressed by 
\begin{equation*}
\mathbb{C}
\left(
E_2(q^{\frac12}), E_2(q), E_2(q^2)
\right)
=
\mathbb{C}
\left(
\varphi^4(q^{\frac12}),  16q^{\frac12} \psi^4(q),   E_2(q)
\right)
=
\mathbb{C}
\left(
\varphi^4(q^{\frac12}),  \varphi^4(-q^{\frac12}),   E_2(q)
\right).
\end{equation*}

\subsection{Note on the cubic theta functions}

\begin{theorem}
\label{thm:cubic-level-6-q-(0,1/3)-(0,2/3)}
{\it
For $q\in\mathbb{C}$ with $|q|<1,$ we have 
\begin{align*}
a(q)=&-\frac{2Q^2-2QR+2 R^2   }{Q-R  },  \\
b^3(q)=&P^3-27PQ^2-54Q^3=P^3-27PR^2+54R^3=
-\frac{2 (Q-2 R)^2 (2 Q-R)^2 (Q+R)^2}{(Q-R)^3}, \\
c^3(q)=&27PQ^2+54Q^3=27PR^2-54R^3
=-\frac{54 Q^2 R^2}{Q-R}.
\end{align*}
}
\end{theorem}

\begin{proof}
The theorem can be proved in the same way as Theorem \ref{thm:cubic-q-(1,2/3)}. 
\end{proof}

\begin{theorem}
\label{thm:cubic-level-6-q^{1/2}-(0,1/3)-(0,2/3)}
{\it
For $q\in\mathbb{C}$ with $|q|<1,$ we have 
\begin{align*}
a(q^{\frac12})=&P+6 Q=\frac{2 \left(2 Q^2-2 Q R-R^2\right)}{Q-R},  \\
b^3(q^{\frac12})=&P^3-9 P^2 Q+108 Q^3=P^3-9 P^2 Q+54 P Q^2-54 P R^2+216 Q^3+108 R^3  \\
=&\frac{4 (Q-2 R) (2 Q-R)^4 (Q+R)}{(Q-R)^3} \\
c^3(q^{\frac12})=&27 P^2 Q+108 P Q^2+108 Q^3 =27 P^2 Q+54 P Q^2+54 P R^2-108 R^3
=\frac{108 Q R^4}{(Q-R)^2}.
\end{align*}
}
\end{theorem}

\begin{proof}
The theorem can be proved in the same way as Theorem \ref{thm:cubic-q^2-(1,2/3)}. 
\end{proof}

\begin{theorem}
\label{thm:cubic-level-6-q^2-(0,1/3)-(0,2/3)}
{\it
For $q\in\mathbb{C}$ with $|q|<1,$ we have 
\begin{align*}
a(q^2)=&P+3Q-3R=\frac{Q^2-4 Q R+R^2}{Q-R},  \\
b^3(q^2)
=&\frac{ \left(4 P^3+21 P^2 Q-21 P^2 R+33 P Q^2-150 P Q R+87 P R^2+18 Q^3-198 Q^2 R+198 Q R^2-126 R^3\right)}{4}  \\
=&\frac{\left(4 P^3+21 P^2 Q-21 P^2 R+87 P Q^2-150 P Q R+33 P R^2+126 Q^3-198 Q^2 R+198 Q R^2-18 R^3\right)}{4}  \\
=&\frac{(Q-2 R) (2 Q-R) (Q+R)^4}{2 (Q-R)^3}, \\
c^3(q^2)
=&\frac{3\left(5 P^2 Q-5 P^2 R+25 P Q^2-22 P Q R+7 P R^2+30 Q^3-42 Q^2 R+42 Q R^2+6 R^3\right)  }{4}  \\
=&\frac{3\left(5 P^2 Q-5 P^2 R+7 P Q^2-22 P Q R+25 P R^2-6 Q^3-42 Q^2 R+42 Q R^2-30 R^3\right)}{4}  \\
=&-\frac{27}{2} Q R (Q-R). 
\end{align*}
}
\end{theorem}

\begin{proof}
The theorem can be proved in the same way as Theorem \ref{thm:cubic-q^2-(1,2/3)}. 
\end{proof}

\begin{theorem}
\label{thm:cubic-ODE-(k, l,m,n)}
{\it
For $q\in\mathbb{C}$ with $|q|<1,$ set 
\begin{equation*}
k=a(q^{\frac12}), \,\,l=a(q), \,\,m=a(q^2), \,\, \text{and} \,\,n=E_2(q).
\end{equation*}
Then, we have 
\begin{align*}
q\frac{d}{dq} k=&\frac{2 k^2 l+k l^2-4 l^3 +k n}{12} , \,\,
 q\frac{d}{dq} m=\frac{ 2 l^3+l^2 m-4 l m^2+m n }{12},   \\
q\frac{d}{d q} l=&
\frac{  k^3-3 k l^2+l^3+l n  }{12} 
=\frac{ l^3+6 l^2 m-8 m^3 +l n}{12},  \\
q\frac{d}{dq}n=&
\frac{ -2 k^3 l+6 k l^3-5 l^4+n^2 }{12} =\frac{ -5 l^4-12 l^3 m+16 l m^3+n^2  }{12}. 
\end{align*}
and 
\begin{equation}
\label{eqn:relation-(k,l,m)}
 k^2 - 3 l^2 - 2 k m + 4 m^2=0.
\end{equation}
}
\end{theorem}

\begin{proof}
Theorems \ref{thm:cubic-level-6-q^{1/2}-(0,1/3)-(0,2/3)} and \ref{thm:cubic-level-6-q^2-(0,1/3)-(0,2/3)} yield 
\begin{equation*}
k=P+6 Q, \,\,l=P, \,\,m=P+3 Q-3 R, \,\,\mathrm{and} \,\,
P=l, \,\,Q=\frac16 k-\frac16 l, \,\,R=\frac16 k-\frac16 l-\frac13 m. 
\end{equation*}
The relation (\ref{eqn:relation-(k,l,m)}) follows from the relation (\ref{eqn:(P,Q,R)-relation}). 
\par
We prove the second formula of  $q d \,l/dq.$ 
The others can be proved in the same way.  
Theorem \ref{thm-level6-E_2(q)-(0,1/3)-(0,2/3)} implies that 
\begin{align*}
q\frac{d}{d q} l=&q\frac{d}{d q} P
=\frac{ -P^3+108 P R^2-216 R^3+PS  }{ 12 } \\
=&
\frac{-k^3+6 k^2 m+3 k l^2-12 k m^2+l^3-6 l^2 m+l n+8 m^3 }{12}  \\
&\quad 
+
\left(\frac{k}{12}-\frac{m}{3}\right) \left(k^2-2 k m-3 l^2+4 m^2\right)  \\
=&
\frac{ l^3+6 l^2 m-8 m^3 +l n}{12}. 
\end{align*}
\end{proof}

\subsubsection*{Remark}
Theorem \ref{thm:E_2-(1,1/3)-(0,1/3)-(0,2/3)} yields 
\begin{equation*}
E_2(q^{\frac12})=-k^2-2 k l+2 l^2+2 n, \,\,
E_2(q^2)=\frac{  k^2-2 k m+6 l m-2 m^2+3 n }{6}
=\frac{  l^2+2 l m-2 m^2+n }{2}. 
\end{equation*}
Changing $q\rightarrow q^{\frac12}, $ 
we can add $a(q^{\frac12},), a(q^{\frac14}), a(q^{\frac18}), \ldots$ to the differential fields $\mathbb{C}(k,l,n),$ and $\mathbb{C}(l,m,n)$; 
\begin{align*}
\mathbb{C}(a(q^{\frac12}), a(q), E_2(q)) &\subset \mathbb{C}(a(q^{\frac14}), a(q^{\frac12}), a(q), E_2(q)) \subset
\mathbb{C}(a(q^{\frac18}),   a(q^{\frac14}), a(q^{\frac12}), a(q), E_2(q))\subset \cdots,  \\
\mathbb{C}(a(q), a(q^2), E_2(q)) &\subset \mathbb{C}(a(q^{\frac12}), a(q), a(q^2), E_2(q)) \subset
\mathbb{C}( a(q^{\frac14}), a(q^{\frac12}), a(q), a(q^2), E_2(q))\subset \cdots.  
\end{align*}
\par
Moreover, 
changing $q\rightarrow q^2,$ 
we can add $a(q^2), a(q^4), a(q^8), \ldots$ to the differential fields $\mathbb{C}(k,l,n),$ and $\mathbb{C}(l,m,n)$; 
\begin{align*}
\mathbb{C}(a(q^{\frac12}), a(q), E_2(q)) &\subset \mathbb{C}(a(q^{\frac12}), a(q), a(q^2), E_2(q)) \subset
\mathbb{C}(a(q^{\frac12}), a(q), a(q^2), a(q^4), E_2(q))\subset \cdots,  \\
\mathbb{C}(a(q), a(q^2), E_2(q)) &\subset \mathbb{C}(a(q), a(q^2), a(q^4), E_2(q)) \subset \mathbb{C}(a(q), a(q^2), a(q^4), a(q^8), E_2(q))
\subset \cdots. 
\end{align*}

\subsection{Note on $E_4(q)$ and $E_6(q)$ $\cdots$ (2)}
From Berndt \cite[pp. 126-128]{Berndt} and our papers \cite{Matsuda-1, Matsuda-3}, 
we recall 
\begin{align}
&E_4(q)=9 a^4(q)-8 a(q) b^3(q), \,\, E_6(q)=-27 a^6(q)+36 a^3(q) b^3(q)-8 b^6(q), \label{eqn:E_4-E_6-q-cubic}   \\
&E_4(q^3)=a^4(q)-\frac89 a(q) c^3(q), \,\,E_6(q^3)=a^6(q) -\frac43 a^3(q) c^3(q)+\frac{8}{27} c^6(q),   \label{eqn:E_4-E_6-q^3-cubic} 
\end{align}
and 
\begin{align}
&E_4(q^{\frac12})=G^2+14 GH+H^2, \,\, E_6(q^{\frac12})=(G+H)(G^2-34 GH+H^2), \label{eqn:E_4-E_6-q^{1/2}-theta}    \\
&E_4(q^2)=G^2-GH+\frac{1}{16} H^2, \,\,E_6(q^2)=\left(G-\frac12 H \right) \left( G^2-GH-\frac{1}{32} H^2  \right),   
\label{eqn:E_4-E_6-q^2-theta}
\end{align}
where $G=\varphi^4(q^{\frac12})$ and $H=16 q^{\frac12} \psi^4(q).$

\begin{theorem}
\label{thm:E_2-E_4-E_6-q^3-(0,1/3)-(0,2/3)}
{\it 
For $q\in\mathbb{C}$ with $|q|<1,$ we have 
\begin{align*}
E_2(q^3)=&\frac23 a^2(q)+\frac13 E_2(q)= \frac{8 \left(Q^2-Q R+R^2\right)^2}{3 (Q-R)^2}+\frac13 E_2(q),   \\
E_4(q^3)=&P^4-24 P^2 Q^2-48 P Q^3=P^4-24 P^2 R^2+48 P R^3  \\
=&
\frac{16 \left(Q^2-Q R+R^2\right) \left(Q^6-3 Q^5 R+5 Q^3 R^3-3 Q R^5+R^6\right)}{(Q-R)^4}, \\ 
E_6(q^3) 
=&P^6-36 P^4 Q^2-72 P^3 Q^3+216 P^2 Q^4+864 P Q^5+864 Q^6  \\
=&P^6-36 P^4 R^2+72 P^3 R^3+216 P^2 R^4-864 P R^5+864 R^6  \\ 
=&\frac{32 \left(Q^4-2 Q^3 R+4 Q R^3-2 R^4\right) \left(2 Q^4-4 Q^3 R+2 Q R^3-R^4\right) \left(Q^4-2 Q^3 R-2 Q R^3+R^4\right)}{(Q-R)^6},   
\end{align*}
and
\begin{align*}
&q^3(q^3;q^3)_{\infty}^{24}
=
Q^6 (P-6 Q) (P+2 Q)^3 (P+3 Q)^2
=
R^6 (P-3 R)^2 (P-2 R)^3 (P+6 R)  \\
=&
\frac{16 Q^6 R^6 (Q-2 R)^2 (2 Q-R)^2 (Q+R)^2}{(Q-R)^6}, \\
&j(q^3)
=
\frac{P^3 \left(P^3-24 P Q^2-48 Q^3\right)^3}{Q^6 (P-6 Q) (P+2 Q)^3 (P+3 Q)^2}
=
\frac{P^3 \left(P^3-24 P R^2+48 R^3\right)^3}{R^6 (P-3 R)^2 (P-2 R)^3 (P+6 R)}  \\
=&
\frac{(3 Q-T)^3 \left(3 Q^3-3 Q^2 T+9 Q T^2-T^3\right)^3}{Q^6 T^2 (Q-T)^3 (9 Q-T)} 
=
\frac{(3 R+U)^3 \left(3 R^3+3 R^2 U+9 R U^2+U^3\right)^3}{R^6 U^2 (R+U)^3 (9 R+U)} \\
=&
\frac{(3 t-1)^3 \left(3 t^3-3 t^2+9 t-1\right)^3}{ t^6 (t-1)^3  (9 t-1)}
=
\frac{(3 u+1)^3 \left(3 u^3+3 u^2+9 u+1\right)^3}{u^6 (u+1)^3 (9 u+1)}  \\
=&
\frac{256 \left(Q^2-Q R+R^2\right)^3 \left(Q^6-3 Q^5 R+5 Q^3 R^3-3 Q R^5+R^6\right)^3}
{Q^6 R^6 (Q-2 R)^2(2 Q-R)^2 (Q-R)^6  (Q+R)^2},   \\
\end{align*}
where 
\begin{equation*}
t(q)=\frac{Q(q)}{T(q)}=q^{\frac12} 
\frac{
(q^{\frac12} ;q^{\frac12}  )_{\infty}^4   (q^3; q^3)_{\infty}^8
}
{
(q; q)_{\infty}^8 (q^{\frac32}; q^{\frac32})_{\infty}^4
}, \,\,
u(q)=\frac{R(q)}{U(q)}
=
q^{\frac12} 
\frac{
(q ;q )_{\infty}^4   (q^{\frac32}; q^{\frac32})_{\infty}^4   (q^6; q^6)_{\infty}^4
}
{
(q^{\frac12}; q^{\frac12})_{\infty}^4        (q^2 ;q^2 )_{\infty}^4  (q^3 ;q^3 )_{\infty}^4  
}. 
\end{equation*}
}
\end{theorem}

\begin{proof}
The theorem follows from equation (\ref{eqn:E_4-E_6-q^3-cubic}).  
\end{proof}

\begin{theorem}
\label{thm:E_2-E_4-E_6-q^{1/2}-(0,1/3)-(0,2/3)}
{\it
For $q\in\mathbb{C}$ with $|q|<1,$ 
we have 
\begin{align*}
&E_2(q^{\frac12})=-P^2-24 P Q-36 Q^2+2S=\frac{4 \left(2 Q^4-4 Q^3 R+12 Q^2 R^2-10 Q R^3-R^4\right)}{(Q-R)^2}+2E_2(q), \\
&E_4(q^{\frac12})
=(P+6 Q) \left(P^3+234 P^2 Q+972 P Q^2+1080 Q^3\right)  \\
&=(P+6 Q) \left(P^3+234 P^2 Q+540 P Q^2+432 P R^2+216 Q^3-864 R^3\right) \\
&=
\frac{16 \left(2 Q^2-2 Q R-R^2\right) \left(8 Q^6-24 Q^5 R+12 Q^4 R^2+16 Q^3 R^3+102 Q^2 R^4-114 Q R^5-R^6\right)}{(Q-R)^4}, \\
&E_6(q^{\frac12})=\left(P^2+24 P Q+36 Q^2\right) \left(P^4-528 P^3 Q-4536 P^2 Q^2-13824 P Q^3-14256 Q^4\right) \\
&=
\left(P^2+24 P Q+36 Q^2\right) \times \\
&\times \Big(  P^4-528 P^3 Q-4968 P^2 Q^2-10944 P Q^3  -6768 Q^4 -144 P^2 R^2+576 R^4\\
&\hspace{20mm} +576 P^2 Q R-4608 P Q^2 R+1152 P Q R^2+12096 Q^2 R^2  -11520 Q^3 R -4608 Q R^3  \Big) \\
&=
\frac{64 \left(2 Q^4-4 Q^3 R+12 Q^2 R^2-10 Q R^3-R^4\right) }{(Q-R)^6}\times  \\
&\times \left(32 Q^8-128 Q^7 R-64 Q^6 R^2+640 Q^5 R^3-400 Q^4 R^4-416 Q^3 R^5+68 Q^2 R^6+268 Q R^7-R^8\right),
\end{align*}
and
\begin{align*}
& q^{\frac12}(q^{\frac12};q^{\frac12} )_{\infty}^{24}= Q (P-6 Q)^6 (P+2 Q)^2 (P+3 Q)^3
=\frac{256Q R^4 (Q-2 R)^3 (2 Q-R)^{12} (Q+R)^3}{(Q-R)^{11}},  \\
&j(q^{\frac12})=\frac{(P+6 Q)^3 \left(P^3+234 P^2 Q+972 P Q^2+1080 Q^3\right)^3}{Q (P-6 Q)^6 (P+2 Q)^2 (P+3 Q)^3} \\
&\hspace{5mm}=\frac{16 \left(2 Q^2-2 Q R-R^2\right)^3 \left(8 Q^6-24 Q^5 R+12 Q^4 R^2+16 Q^3 R^3+102 Q^2 R^4-114 Q R^5-R^6\right)^3}{Q R^4 (Q-2 R)^3 (Q-R) (2 Q-R)^{12} (Q+R)^3}. 
\end{align*}
}
\end{theorem}

\begin{proof}
The formulas of $E_4(q^{\frac12})$ can be obtained from equation (\ref{eqn:E_4-E_6-q-cubic}).  
The formulas of $E_6(q^{\frac12})$ can be derived from equations (\ref{eqn:E_4-E_6-q-cubic}) and (\ref{eqn:E_4-E_6-q^{1/2}-theta}).
\end{proof}

\begin{theorem}
\label{thm:E_2-E_4-E_6-q^{3/2}-(0,1/3)-(0,2/3)}
{\it
For $q\in\mathbb{C}$ with $|q|<1,$ 
we have 
\begin{align*}
E_4(q^{\frac32})
=&
(P+6 Q) \left(P^3-6 P^2 Q+12 P Q^2+120 Q^3\right) \\
=&(P+6 Q) \left(P^3-6 P^2 Q+60 P Q^2-48 P R^2+216 Q^3+96 R^3\right)  \\
=&
\frac{16 \left(2 Q^2-2 Q R-R^2\right) \left(8 Q^6-24 Q^5 R+12 Q^4 R^2+16 Q^3 R^3-18 Q^2 R^4+6 Q R^5-R^6\right)}{(Q-R)^4},   \\
E_6(q^{\frac32})
=&\left(P^2-12 Q^2\right) \left(P^4-24 P^2 Q^2-576 P Q^3-1584 Q^4\right)  \\
=&
\left(P^2-12 Q^2\right) \left(P^4+120 P^2 Q^2-144 P^2 R^2-1008 Q^4+576 Q^2 R^2+576 R^4\right)  \\
=&
\frac{64 \left(2 Q^4-4 Q^3 R+2 Q R^3-R^4\right) }{(Q-R)^6} \times  \\
&\times \left(32 Q^8-128 Q^7 R+128 Q^6 R^2+64 Q^5 R^3-160 Q^4 R^4+64 Q^3 R^5-4 Q^2 R^6+4 Q R^7-R^8\right),
\end{align*}
and
\begin{align*}
&q^{\frac32}(q^{\frac32};q^{\frac32} )_{\infty}^{24}=Q^3 (P-6 Q)^2 (P+2 Q)^6 (P+3 Q)=
\frac{256Q^3 R^{12} (Q-2 R) (2 Q-R)^4 (Q+R)}{(Q-R)^9},    \\
&j(q^{\frac32})=\frac{(P+6 Q)^3 \left(P^3-6 P^2 Q+12 P Q^2+120 Q^3\right)^3}{Q^3 (P-6 Q)^2 (P+2 Q)^6 (P+3 Q)}  \\
&\hspace{10mm}=
\frac{16 \left(2 Q^2-2 Q R-R^2\right)^3 \left(8 Q^6-24 Q^5 R+12 Q^4 R^2+16 Q^3 R^3-18 Q^2 R^4+6 Q R^5-R^6\right)^3}{Q^3 R^{12} (Q-2 R)(2 Q-R)^4 (Q-R)^3  (Q+R)}. 
\end{align*}
}
\end{theorem}

\begin{proof}
The formulas of $E_4(q^{\frac32})$ can be obtained from equation (\ref{eqn:E_4-E_6-q^3-cubic}).  
The formulas of $E_6(q^{\frac32})$ can be derived from equations (\ref{eqn:E_4-E_6-q^3-cubic}) and (\ref{eqn:E_4-E_6-q^{1/2}-theta}).
\end{proof}

\begin{theorem}
\label{thm:E_2-E_4-E_6-q^2-(0,1/3)-(0,2/3)}
{\it
For $q\in\mathbb{C}$ with $|q|<1,$ we have 
\begin{align*}
E_2(q^2)=&
-\frac{Q^4-14 Q^3 R+24 Q^2 R^2-14 Q R^3+R^4}{(Q-R)^2}+\frac12 E_2(q), \\
E_4(q^2)=&
(P+3 Q-3 R) \times \\
\times& \left(P^3+39 P^2 Q-39 P^2 R+177 P Q^2-186 P Q R+69 P R^2+207 Q^3-333 Q^2 R+333 Q R^2+9 R^3\right)  \\
E_6(q^2)
=&\frac{1}{2} \left(P^2-6 Q^2+24 Q R-6 R^2\right) \times \\
&\times \Big(2 P^4-24 P^2 R^2 -504 P Q^3-369 Q^4 +1026 Q^2 R^2-504 Q R^3+63 R^4 -168 P^2 Q^2  -72 Q^3 R\\
&\hspace{20mm}  +96 P^2 Q R +288 P Q^2 R-72 P Q R^2  \Big)  \\
=&
\frac{\left(Q^4-14 Q^3 R+24 Q^2 R^2-14 Q R^3+R^4\right) }{(Q-R)^6} \times \\
&\times\Big(Q^8+260 Q^7 R-1916 Q^6 R^2+5564 Q^5 R^3-7850 Q^4 R^4  \\
&\hspace{70mm}+5564 Q^3 R^5-1916 Q^2 R^6+260 Q R^7+R^8  \Big), 
\end{align*}
and
\begin{align*}
& q^{2}(q^{2};q^{2} )_{\infty}^{24} =
-\frac{ Q R (Q-2 R)^3 (2 Q-R)^3 (Q+R)^{12}}{ 16(Q-R)^8}    \\
&j(q^{2})=
-\frac{16 \left(Q^2-4 Q R+R^2\right)^3 \left(Q^6-120 Q^5 R+483 Q^4 R^2-736 Q^3 R^3+483 Q^2 R^4-120 Q R^5+R^6\right)^3}{Q R (Q-2 R)^3 (2 Q-R)^3  (Q-R)^4  (Q+R)^{12}}. 
\end{align*}
}
\end{theorem}

\begin{proof}
The formulas of $E_4(q^2)$ can be obtained from equation (\ref{eqn:E_4-E_6-q-cubic}).  
The formulas of $E_6(q^2)$ can be derived from equation (\ref{eqn:E_4-E_6-q^2-theta}).
\end{proof}

\begin{theorem}
\label{thm:E_2-E_4-E_6-q^6-(0,1/3)-(0,2/3)}
{\it
For $q\in\mathbb{C}$ with $|q|<1,$ we have 
\begin{align*}
E_2(q^6)=&\frac{Q^4-2 Q^3 R+12 Q^2 R^2-2 Q R^3+R^4}{3 (Q-R)^2}+\frac16 E_2(q), \\
E_4(q^6)
=&P^4-12 P^2 Q^2-12 P^2 R^2+33 Q^4+78 Q^2 R^2+33 R^4 \\
=&
\frac{\left(Q^2-4 Q R+R^2\right) \left(Q^6+3 Q^4 R^2-16 Q^3 R^3+3 Q^2 R^4+R^6\right)}{(Q-R)^4}, \\
E_6(q^6)
=&\frac{1}{2} \left(P^2-6 Q^2-6 R^2\right) \left(2 P^4-24 P^2 Q^2-24 P^2 R^2+63 Q^4+162 Q^2 R^2+63 R^4\right)  \\
=&\frac{\left(Q^4-2 Q^3 R-2 Q R^3+R^4\right) }{(Q-R)^6}\times   \\
&\hspace{15mm} \times \left(Q^8-4 Q^7 R+4 Q^6 R^2+68 Q^5 R^3-170 Q^4 R^4+68 Q^3 R^5+4 Q^2 R^6-4 Q R^7+R^8\right),
\end{align*}
and 
\begin{align*}
&q^6(q^6;q^6)_{\infty}^{24}
=
-\frac{Q^3 R^3 (Q-2 R) (2 Q-R) (Q+R)^4}{16},  \\
&j(q^6)
=
-\frac{16 \left(Q^2-4 Q R+R^2\right)^3 \left(Q^6+3 Q^4 R^2-16 Q^3 R^3+3 Q^2 R^4+R^6\right)^3}{Q^3 R^3 (Q-2 R) (Q-R)^{12} (2 Q-R) (Q+R)^4}. 
\end{align*}
}
\end{theorem}

\begin{proof}
The formulas of $E_4(q^6)$ and $E_6(q^6)$ can be obtained from equation (\ref{eqn:E_4-E_6-q^2-theta}).  
\end{proof}

\subsection{Application to number theory}

\begin{theorem}
\label{thm:two-sums-x^2+xy+y^2-(1,4)}
{\it
For $q\in\mathbb{C}$ with $|q|<1,$ 
we have 
\begin{equation}
\label{eqn:a(q)-1-4}
a(q)a(q^4)=1+6\sum_{n=1}^{\infty} (\sigma(n)-3\sigma(n/2)+4\sigma(n/4)-3\sigma(n/3) +9\sigma(n/6)-12\sigma(n/12)        )q^n.   
\end{equation}
}
\end{theorem}

\begin{proof}
Equation (\ref{eqn:a(q)-1-1}) and the relation (\ref{eqn:relation-(k,l,m)}) yield 
\begin{equation*}
km=\frac12 k^2-\frac32 l^2+2 m^2
=\frac14
\left\{
-E_2(q^{\frac12})+3E_2(q^{\frac32})
+3E_2(q)-9 E_2(q^3)
-4 E_2(q^2)+12 E_2(q^6)
\right\}.
\end{equation*} 
The theorem can be obtained by changing $q\rightarrow q^2.$   
\end{proof}

\begin{theorem}
\label{thm:three-sums-x^2+xy+y^2-(1,2,4)}
{\it
For $q\in\mathbb{C}$ with $|q|<1,$ 
we have 
\begin{align*}
a(q)a(q^2) a(q^4)=&
1
+\frac32
\sum_{n=1}^{\infty} \left( \sum_{d|n} d^2\left( \frac{d}{3}  \right)  \right)q^n
+\frac92
\sum_{n=1}^{\infty} \left( \sum_{d|n} d^2\left( \frac{n/d}{3}  \right)  \right)  q^n   \notag  \\
&+\frac32
\sum_{n=1}^{\infty} \left( \sum_{d|n} d^2\left( \frac{d}{3}  \right)  \right)q^{2n}
-\frac92
\sum_{n=1}^{\infty} \left( \sum_{d|n} d^2\left( \frac{n/d}{3}  \right)  \right)  q^{2n}   \notag   \\
&-12
\sum_{n=1}^{\infty} \left( \sum_{d|n} d^2\left( \frac{d}{3}  \right)  \right)q^{4n}
-36
\sum_{n=1}^{\infty} \left( \sum_{d|n} d^2\left( \frac{n/d}{3}  \right)  \right)  q^{4n}.  
\end{align*}
}
\end{theorem}

\begin{proof}
The relation (\ref{eqn:relation-(k,l,m)}) and equations (\ref{eqn:a(q)-1-1-2-b-c}) and  (\ref{eqn:a(q)-1-2-2-b-c}) imply that 
\begin{equation*}
klm=\frac12 k^2l-\frac32 l^3+2 l m^2
=
-\frac16 b^3(q^{\frac12})+\frac16 c^3(q^{\frac12})-\frac16 a^3(q)+\frac43 b^3(q^2)-\frac43 c^3(q^2).
\end{equation*}
The theorem can be obtained by changing $q\rightarrow q^2.$   
\end{proof}

\subsection{Note on Theorem \ref{thm-level6-E_2(q)-(0,1/3)-(0,2/3)}  }
The differential field, $\mathbb{C}\left( P(q), Q(q), R(q), S(q)  \right) $ can be expressed by 
\begin{equation*}
\mathbb{C}\left( P(q), Q(q), R(q), S(q)  \right)= \mathbb{C}\left( Q(q), R(q), S(q)  \right)=
\mathbb{C}\left( a(q^{\frac12}), a(q),  a(q^2), E_2(q)  \right). 
\end{equation*}
This differential field contains the following differential subfields:
\begin{align*}
&\mathbb{C}\left( E_2(q), E_4(q), E_6(q) \right), 
\mathbb{C}\left( E_2(q^2), E_4(q^2), E_6(q^2) \right), \mathbb{C}\left( E_2(q^3), E_4(q^3), E_6(q^3) \right),  \\
&\mathbb{C}\left( E_2(q^{\frac12}), E_4(q^{\frac12}),  E_6(q^{\frac12}) \right),  \,\,  \mathbb{C}\left( E_2(q^{\frac32}), E_4(q^{\frac32}),  E_6(q^{\frac32}  )  \right), \,\,        \mathbb{C}\left( E_2(q^6), E_4(q^6), E_6(q^6) \right),  \\
&\mathbb{C}\left( a(q), E_2(q), b^3(q) \right), \,\,\mathbb{C}\left( a(q^{\frac12}), E_2(q^{\frac12}), b^3(q^{\frac12}) \right), \,\, 
\mathbb{C}\left( a(q^2), E_2(q^2), b^3(q^2) \right),  \\
&\mathbb{C}\left( P(q), Q(q), S(q)  \right), \,\,\mathbb{C}\left( P(q), R(q), S(q)  \right),   \mathbb{C}\left( a(q), a(q^2), E_2(q)  \right),   \\
&\mathbb{C}\left( E_2(q^{\frac12}), E_2(q), E_2(q^2) \right), \,\,  \mathbb{C}\left( E_2(q^{\frac32}), E_2(q^3), E_2(q^6) \right), 
\mathbb{C}\left( E_2(q^{\frac12}), E_2(q), E_2(q^2), E_2(q^{\frac32}) \right), \\
&\mathbb{C}\left( E_2(q^{\frac12}), E_2(q), E_2(q^2), E_2(q^3) \right), \,\,
\mathbb{C}\left( E_2(q^{\frac12}), E_2(q), E_2(q^2), E_2(q^6) \right).
\end{align*}
In particular, 
the differential field, 
$\mathbb{C}\left( a(q^{\frac12}), a(q),  a(q^2), E_2(q)  \right)$ is a two-dimensional extension of 
the differential field, 
$\mathbb{C}\left( a(q), a(q^2), E_2(q)  \right)$ in Theorems \ref{thm-level6-E_2(q)-(1,2/3)} and \ref{thm-level6-E_2(q^6)-(2/3,1)}.
\par
Note that 
\begin{equation*}
\mathbb{C}\left( P(q), Q(q), S(q)  \right)=\mathbb{C}(a(q^{\frac12}), a(q), E_2(q) ). 
\end{equation*}
In particular, changing $q\rightarrow q^2,$ the differential field, $\mathbb{C}\left( P(q), Q(q), S(q)  \right)$ 
can be transformed into the differential field in Theorems \ref{thm-level6-E_2(q)-(1,2/3)} and \ref{thm-level6-E_2(q^6)-(2/3,1)}.

\section{Riccati equations satisfied by modular functions of level six}
\label{sec:Riccati-level-6}

\subsection{Some theta constant identities}

\begin{proposition}
\label{prop:theta-relation-(1/3,2/3)-(1/3,2/3)}
{\it
For every $\tau\in\mathbb{H}^2,$ we have 
\begin{equation*}
\theta^3
\left[
\begin{array}{c}
1 \\
0
\end{array}
\right]
\theta
\left[
\begin{array}{c}
1 \\
\frac23
\end{array}
\right]
-
\theta^4
\left[
\begin{array}{c}
1 \\
\frac13
\end{array}
\right]
+
\theta^4
\left[
\begin{array}{c}
1 \\
\frac23
\end{array}
\right]=0,\hspace{5mm}
\theta^3
\left[
\begin{array}{c}
0 \\
1
\end{array}
\right]
\theta
\left[
\begin{array}{c}
\frac23 \\
1
\end{array}
\right]
+\zeta_3
\theta^4
\left[
\begin{array}{c}
\frac13 \\
1
\end{array}
\right]
-
\theta^4
\left[
\begin{array}{c}
\frac23 \\
1
\end{array}
\right]=0,
\end{equation*}
\begin{equation*}
\theta^3
\left[
\begin{array}{c}
0 \\
1
\end{array}
\right]
\theta
\left[
\begin{array}{c}
0 \\
\frac13
\end{array}
\right]
+
\theta^4
\left[
\begin{array}{c}
1 \\
\frac13
\end{array}
\right]
-
\theta^4
\left[
\begin{array}{c}
0 \\
\frac13
\end{array}
\right]=0,\hspace{5mm}
\theta^3
\left[
\begin{array}{c}
0 \\
0
\end{array}
\right]
\theta
\left[
\begin{array}{c}
0 \\
\frac23
\end{array}
\right]
-
\theta^4
\left[
\begin{array}{c}
1 \\
\frac13
\end{array}
\right]
-
\theta^4
\left[
\begin{array}{c}
0 \\
\frac23
\end{array}
\right]=0,
\end{equation*}
and
\begin{equation*}
\theta^3
\left[
\begin{array}{c}
1 \\
0
\end{array}
\right]
\theta
\left[
\begin{array}{c}
1 \\
\frac23
\end{array}
\right]
-
\theta^4
\left[
\begin{array}{c}
0 \\
\frac13
\end{array}
\right]
+
\theta^4
\left[
\begin{array}{c}
0 \\
\frac23
\end{array}
\right]=0,
\end{equation*}
where $\zeta_3=\exp(2\pi i/3).$
}
\end{proposition}

\begin{proof}
The proposition can be obtained by applying the residue theorem to the following elliptic functions:
\begin{equation*}
\frac
{
\theta^3
\left[
\begin{array}{c}
1 \\
1
\end{array}
\right](z)
}
{
\theta
\left[
\begin{array}{c}
1 \\
\frac13
\end{array}
\right](z)
\theta
\left[
\begin{array}{c}
1 \\
\frac23
\end{array}
\right](z)
\theta
\left[
\begin{array}{c}
1 \\
0
\end{array}
\right](z)
}, \hspace{5mm}
\frac
{
\theta^3
\left[
\begin{array}{c}
1 \\
1
\end{array}
\right](z)
}
{
\theta
\left[
\begin{array}{c}
\frac13\\
1
\end{array}
\right](z)
\theta
\left[
\begin{array}{c}
\frac23 \\
1
\end{array}
\right](z)
\theta
\left[
\begin{array}{c}
0 \\
1
\end{array}
\right](z)
},
\end{equation*}
\begin{equation*}
\frac
{
\theta^3
\left[
\begin{array}{c}
1 \\
1
\end{array}
\right](z)
}
{
\theta
\left[
\begin{array}{c}
1 \\
-\frac13
\end{array}
\right](z)
\theta
\left[
\begin{array}{c}
0 \\
\frac13
\end{array}
\right](z)
\theta
\left[
\begin{array}{c}
0 \\
1
\end{array}
\right](z)
}, \hspace{5mm}
\frac
{
\theta^3
\left[
\begin{array}{c}
1 \\
1
\end{array}
\right](z)
}
{
\theta
\left[
\begin{array}{c}
1\\
\frac13
\end{array}
\right](z)
\theta
\left[
\begin{array}{c}
0 \\
\frac23
\end{array}
\right](z)
\theta
\left[
\begin{array}{c}
0 \\
0
\end{array}
\right](z)
},
\end{equation*}
and
\begin{equation*}
\frac
{
\theta^3
\left[
\begin{array}{c}
1 \\
1
\end{array}
\right](z)
}
{
\theta
\left[
\begin{array}{c}
0 \\
\frac13
\end{array}
\right](z)
\theta
\left[
\begin{array}{c}
0 \\
\frac23
\end{array}
\right](z)
\theta
\left[
\begin{array}{c}
1 \\
0
\end{array}
\right](z)
}.
\end{equation*}
\end{proof}

\subsection{Some theta functional formulas}

\begin{proposition}
\label{prop:theta-function-identities-level-6}
{\it
For every $(z,\tau)\in\mathbb{C}\times\mathbb{H}^2,$ we have 
\begin{align}
&
\theta^2
\left[
\begin{array}{c}
1 \\
\frac23
\end{array}
\right]
\theta
\left[
\begin{array}{c}
1 \\
\frac13
\end{array}
\right](z,\tau)
\theta
\left[
\begin{array}{c}
1 \\
\frac53
\end{array}
\right](z,\tau)
-
\theta^2
\left[
\begin{array}{c}
1 \\
\frac13
\end{array}
\right]
\theta 
\left[
\begin{array}{c}
1 \\
\frac23
\end{array}
\right](z,\tau)
\theta
\left[
\begin{array}{c}
1 \\
\frac43
\end{array}
\right](z,\tau)      \notag \\
&
\hspace{30mm}
+\theta
\left[
\begin{array}{c}
1 \\
0
\end{array}
\right]
\theta
\left[
\begin{array}{c}
1 \\
\frac23
\end{array}
\right]
\theta^2
\left[
\begin{array}{c}
1 \\
1
\end{array}
\right](z,\tau)=0,  \label{eqn:two-theta-sum-2-deno-3-(0)}
\end{align}
\begin{align}
&
-\zeta_3
\theta^2
\left[
\begin{array}{c}
\frac23 \\
1
\end{array}
\right]
\theta
\left[
\begin{array}{c}
\frac13 \\
1
\end{array}
\right](z,\tau)
\theta
\left[
\begin{array}{c}
\frac53 \\
1
\end{array}
\right](z,\tau)
-
\theta^2
\left[
\begin{array}{c}
\frac13 \\
1
\end{array}
\right]
\theta 
\left[
\begin{array}{c}
\frac23 \\
1
\end{array}
\right](z,\tau)
\theta
\left[
\begin{array}{c}
\frac43 \\
1
\end{array}
\right](z,\tau)      \notag \\
&
\hspace{30mm}
+\theta
\left[
\begin{array}{c}
0 \\
1
\end{array}
\right]
\theta
\left[
\begin{array}{c}
\frac23 \\
1
\end{array}
\right]
\theta^2
\left[
\begin{array}{c}
1 \\
1
\end{array}
\right](z,\tau)=0,  \label{eqn:two-theta-sum-2-deno-3-(1)}
\end{align}
\begin{align}
&
\theta^2
\left[
\begin{array}{c}
0 \\
\frac13
\end{array}
\right]
\theta
\left[
\begin{array}{c}
1 \\
\frac13
\end{array}
\right](z,\tau)
\theta
\left[
\begin{array}{c}
1 \\
\frac53
\end{array}
\right](z,\tau)
+
\theta^2
\left[
\begin{array}{c}
1 \\
\frac13
\end{array}
\right]
\theta 
\left[
\begin{array}{c}
0 \\
\frac13
\end{array}
\right](z,\tau)
\theta
\left[
\begin{array}{c}
0 \\
\frac53
\end{array}
\right](z,\tau)      \notag \\
&
\hspace{30mm}
-
\theta
\left[
\begin{array}{c}
0 \\
1
\end{array}
\right]
\theta
\left[
\begin{array}{c}
0 \\
\frac13
\end{array}
\right]
\theta^2
\left[
\begin{array}{c}
1 \\
1
\end{array}
\right](z,\tau)=0,  \label{eqn:two-theta-sum-2-deno-3-(2)}
\end{align}
and 
\begin{align}
&
\theta^2
\left[
\begin{array}{c}
0 \\
\frac23
\end{array}
\right]
\theta
\left[
\begin{array}{c}
1 \\
\frac13
\end{array}
\right](z,\tau)
\theta
\left[
\begin{array}{c}
1 \\
\frac53
\end{array}
\right](z,\tau)
+
\theta^2
\left[
\begin{array}{c}
1 \\
\frac13
\end{array}
\right]
\theta 
\left[
\begin{array}{c}
0 \\
\frac23
\end{array}
\right](z,\tau)
\theta
\left[
\begin{array}{c}
0 \\
\frac43
\end{array}
\right](z,\tau)      \notag \\
&
\hspace{30mm}
-
\theta
\left[
\begin{array}{c}
0 \\
0
\end{array}
\right]
\theta
\left[
\begin{array}{c}
0 \\
\frac23
\end{array}
\right]
\theta^2
\left[
\begin{array}{c}
1 \\
1
\end{array}
\right](z,\tau)=0,  \label{eqn:two-theta-sum-2-deno-3-(3)}
\end{align}
where $\zeta_3=\exp(2\pi i/3).$ 
}
\end{proposition}

\begin{proof}
We prove equation (\ref{eqn:two-theta-sum-2-deno-3-(0)}). 
The others can be proved in the same way. 
We first note that 
$
\dim \mathcal{F}_{2}\left[
\begin{array}{c}
0 \\
0
\end{array}
\right]
=
2,
$ 
and 
\begin{equation*}
\theta
\left[
\begin{array}{c}
1 \\
\frac13
\end{array}
\right](z,\tau)
\theta
\left[
\begin{array}{c}
1 \\
\frac53
\end{array}
\right](z,\tau)
, \,
\theta 
\left[
\begin{array}{c}
1 \\
\frac23
\end{array}
\right](z,\tau)
\theta
\left[
\begin{array}{c}
1 \\
\frac43
\end{array}
\right](z,\tau), \,
\theta^2
\left[
\begin{array}{c}
1 \\
1
\end{array}
\right](z,\tau) \in
\mathcal{F}_{2}\left[
\begin{array}{c}
0 \\
0
\end{array}
\right].  
\end{equation*}
Therefore, there exist complex numbers $x_1, x_2$, and $ x_3$ not all zero such that 
 \begin{align*}
 &
 x_1
\theta
\left[
\begin{array}{c}
1 \\
\frac13
\end{array}
\right](z,\tau)
\theta
\left[
\begin{array}{c}
1 \\
\frac53
\end{array}
\right](z,\tau)
+x_2
 \theta
\left[
\begin{array}{c}
1\\
\frac23
\end{array}
\right](z,\tau)
\theta
\left[
\begin{array}{c}
1 \\
\frac43
\end{array}
\right](z,\tau)
+
x_3
\theta^2
\left[
\begin{array}{c}
1 \\
1
\end{array}
\right](z,\tau)=0. 
\end{align*}
Note that, in the fundamental parallelogram, 
the zeroes of 
$
\theta
\left[
\begin{array}{c}
1 \\
\frac13
\end{array}
\right](z),
$ 
$
\theta
\left[
\begin{array}{c}
1 \\
\frac23
\end{array}
\right](z),
$ 
and
$
\theta
\left[
\begin{array}{c}
1 \\
1
\end{array}
\right](z)
$ 
are $z=1/3,$ $1/6$, and $0$,
respectively.
Substituting $z=1/3, 1/6$ and $0,$ we have 
\begin{alignat*}{4}
&
&    
&
x_2
\theta
\left[
\begin{array}{c}
1\\
\frac23
\end{array}
\right]  
\theta
\left[
\begin{array}{c}
1\\
0
\end{array}
\right]  
&    
&+x_3
\theta^2
\left[
\begin{array}{c}
1 \\
\frac13
\end{array}
\right]  
&  
&=0,  \\
&-x_1
\theta
\left[
\begin{array}{c}
1\\
\frac23
\end{array}
\right]  
\theta
\left[
\begin{array}{c}
1\\
0
\end{array}
\right]  
&  
&
&
&+x_3
\theta^2
\left[
\begin{array}{c}
1 \\
\frac23
\end{array}
\right]  
&  
&=0,  \\
&
-
x_1
\theta^2
\left[
\begin{array}{c}
1 \\
\frac13
\end{array}
\right]  
&
-&x_2
\theta^2
\left[
\begin{array}{c}
1 \\
\frac23
\end{array}
\right] 
&
&
&
&=0. 
\end{alignat*}
Solving this system of equations, 
we have 
\begin{equation*}
(x_1,x_2,x_3)=\alpha
\left(
\theta^2
\left[
\begin{array}{c}
1 \\
\frac23
\end{array}
\right], 
-
\theta^2
\left[
\begin{array}{c}
1 \\
\frac13
\end{array}
\right], 
\theta
\left[
\begin{array}{c}
1 \\
\frac23
\end{array}
\right] 
\theta\left[
\begin{array}{c}
1\\
0
\end{array}
\right] 
\right) \,\,\,   \text{for some} \,\,\alpha\in\mathbb{C}\setminus\{0\},
\end{equation*}
which proves the proposition. 
\end{proof}

\subsection{Derivation of a Riccati equation (1)}

\begin{proposition}
\label{prop:2nd-order-theta(1,1/3)-(1,2/3)}
{\it
For every $\tau\in\mathbb{H}^2,$ we have 
\begin{equation*}
\frac
{
\theta^{\prime \prime}
\left[
\begin{array}{c}
1 \\
\frac23
\end{array}
\right] 
}
{
\theta
\left[
\begin{array}{c}
1 \\
\frac23
\end{array}
\right] 
}
-
\frac
{
\theta^{\prime \prime}
\left[
\begin{array}{c}
1 \\
\frac13
\end{array}
\right] 
}
{
\theta
\left[
\begin{array}{c}
1 \\
\frac13
\end{array}
\right] 
}
=
\frac{1}{12}
\frac
{
\left\{
\theta^{\prime}
\left[
\begin{array}{c}
1 \\
1
\end{array}
\right] 
\right\}^2
}
{
\theta^2
\left[
\begin{array}{c}
1 \\
0
\end{array}
\right] 
\theta^2
\left[
\begin{array}{c}
1 \\
\frac13
\end{array}
\right] 
\theta^6
\left[
\begin{array}{c}
1 \\
\frac23
\end{array}
\right] 
}
\left\{
\theta^{8}
\left[
\begin{array}{c}
1 \\
\frac13
\end{array}
\right] 
-10
\theta^4
\left[
\begin{array}{c}
1 \\
\frac13
\end{array}
\right] 
\theta^4
\left[
\begin{array}{c}
1 \\
\frac23
\end{array}
\right] 
+
9
\theta^{8}
\left[
\begin{array}{c}
1 \\
\frac23
\end{array}
\right] 
\right\}. 
\end{equation*}
}
\end{proposition}

\begin{proof}
The proposition can be obtained by 
considering the derivative formulas in Theorem \ref{thm:analogue-Jacobi-(1,2/3)-(2/3,1)} and Corollary \ref{coro:analogue-Jacobi-(1,1/3)-(1/3,1)}, and 
by comparing the coefficients of the term $z^2$ in equation (\ref{eqn:two-theta-sum-2-deno-3-(0)}). 
\end{proof}

\begin{theorem}
\label{thm:q-Riccati-diff-eq-modulus6-(1)}
{\it
For $q\in\mathbb{C}$ with $|q|<1,$ set 
\begin{equation*}
f(q)=
\frac
{
(q^2; q^2)_{\infty}^4 (q^3; q^3)_{\infty}^8
}
{
(q;q)_{\infty}^8
(q^6; q^6)_{\infty}^4
},    \,\,q=\exp(2\pi i \tau).
\end{equation*}
Then, $f$ satisfies the following Riccati equation:
\begin{equation}
\label{eqn:Riccati-q-diff-eq-modulus6-(1)}
q
\frac{d}{dq} f
=
\frac{1}{2^3}
\frac
{
(q; q)_{\infty}^4 (q^3; q^3)_{\infty}^4
}
{
(q^2; q^2)_{\infty}^2 (q^6; q^6)_{\infty}^2
}
\left(
9f^2-10 f+1
\right).
\end{equation}
}
\end{theorem}

\begin{proof}
The heat equation (\ref{eqn:heat}) and 
Proposition \ref{prop:2nd-order-theta(1,1/3)-(1,2/3)}
imply that 
\begin{align*}
&
4\pi i 
\frac
{
d
}
{
d\tau
}
\log 
\theta
\left[
\begin{array}{c}
1 \\
\frac23
\end{array}
\right] 
-
4\pi i 
\frac
{
d
}
{
d\tau
}
\log 
\theta
\left[
\begin{array}{c}
1 \\
\frac13
\end{array}
\right]   \\
=&
\frac
{
\left\{
\theta^{\prime}
\left[
\begin{array}{c}
1 \\
1
\end{array}
\right] 
\right\}^2
}
{
12
\theta^2
\left[
\begin{array}{c}
1 \\
0
\end{array}
\right] 
\theta^2
\left[
\begin{array}{c}
1 \\
\frac13
\end{array}
\right] 
\theta^6
\left[
\begin{array}{c}
1 \\
\frac23
\end{array}
\right] 
}
\left\{
\theta^{8}
\left[
\begin{array}{c}
1 \\
\frac13
\end{array}
\right] 
-10
\theta^4
\left[
\begin{array}{c}
1 \\
\frac13
\end{array}
\right] 
\theta^4
\left[
\begin{array}{c}
1 \\
\frac23
\end{array}
\right] 
+
9
\theta^{8}
\left[
\begin{array}{c}
1 \\
\frac23
\end{array}
\right] 
\right\}
\end{align*}
which shows that 
\begin{align*}
&
4\pi i 
\frac
{
d
}
{
d \tau
}
\log 
\theta^4
\left[
\begin{array}{c}
1 \\
\frac23
\end{array}
\right] 
-
4\pi i 
\frac
{
d
}
{
d\tau
}
\log 
\theta^4
\left[
\begin{array}{c}
1 \\
\frac13
\end{array}
\right]   \\
=&
\frac
{
\left\{
\theta^{\prime}
\left[
\begin{array}{c}
1 \\
1
\end{array}
\right] 
\right\}^2
}
{
3
\theta^2
\left[
\begin{array}{c}
1 \\
0
\end{array}
\right] 
\theta^2
\left[
\begin{array}{c}
1 \\
\frac13
\end{array}
\right] 
\theta^6
\left[
\begin{array}{c}
1 \\
\frac23
\end{array}
\right] 
}
\left\{
\theta^{8}
\left[
\begin{array}{c}
1 \\
\frac13
\end{array}
\right] 
-10
\theta^4
\left[
\begin{array}{c}
1 \\
\frac13
\end{array}
\right] 
\theta^4
\left[
\begin{array}{c}
1 \\
\frac23
\end{array}
\right] 
+
9
\theta^{8}
\left[
\begin{array}{c}
1 \\
\frac23
\end{array}
\right] 
\right\}
\end{align*}
\par
Jacobi's triple product identity (\ref{eqn:Jacobi-triple}) yields 
\begin{equation*}
\frac
{
\left\{
\theta^{\prime}
\left[
\begin{array}{c}
1 \\
1
\end{array}
\right] 
\right\}^2
\theta^2
\left[
\begin{array}{c}
1 \\
\frac13
\end{array}
\right] 
}
{
3
\theta^2
\left[
\begin{array}{c}
1 \\
0
\end{array}
\right] 
\theta^2
\left[
\begin{array}{c}
1 \\
\frac23
\end{array}
\right] 
}
=
\pi^2
\frac
{
\eta^4(\tau) \eta^4(3\tau)
}
{
\eta^2(2\tau) \eta^2(6\tau)
}, \,\,
\frac
{
\theta^4
\left[
\begin{array}{c}
1 \\
\frac13
\end{array}
\right] 
}
{
\theta^4
\left[
\begin{array}{c}
1 \\
\frac23
\end{array}
\right] 
}
=9
\frac
{
\eta^4(2\tau) \eta^8(3\tau)
}
{
\eta^8(\tau) \eta^4(6\tau)
}
=9f. 
\end{equation*}
\par
Setting 
\begin{equation*}
(X,Y)=
\left(
\theta^4
\left[
\begin{array}{c}
1 \\
\frac13
\end{array}
\right], 
\theta^4
\left[
\begin{array}{c}
1 \\
\frac23
\end{array}
\right]  
\right), \,\,
\end{equation*}
we obtain 
\begin{equation*}
4 \pi i
\left(
\frac{d X}{d\tau} Y
-
X \frac{d Y}{d\tau}
\right)
=
-
\frac
{
\left\{
\theta^{\prime}
\left[
\begin{array}{c}
1 \\
1
\end{array}
\right] 
\right\}^2
\theta^2
\left[
\begin{array}{c}
1 \\
\frac13
\end{array}
\right] 
}
{
3
\theta^2
\left[
\begin{array}{c}
1 \\
0
\end{array}
\right] 
\theta^2
\left[
\begin{array}{c}
1 \\
\frac23
\end{array}
\right] 
}
\left( X^2-10XY+9Y^2\right),
\end{equation*}
which implies that 
\begin{equation*}
\frac{d}{d\tau} 
\left(
9 f
\right)
=
\frac{d}{d\tau}
\left(
\frac{X}{Y}
\right)
=
\frac
{
\displaystyle \frac{dX}{d\tau}Y-X\frac{dY}{d\tau}
}
{
Y^2
}
=
\frac{2 \pi i}{2^3}
\frac
{
\eta^4(\tau) \eta^4(3\tau)
}
{
\eta^2(2\tau) \eta^2(6\tau)
}
\left(
81 f^2-90f+9
\right).
\end{equation*}
The theorem can be obtained by considering that 
$
\displaystyle
\frac{d}{d\tau}=\frac{dq}{d\tau}\frac{d}{dq}=2\pi i q \frac{d}{dq}. 
$
\end{proof}

\subsubsection*{Remark}
By Theorem \ref{thm:E_2(q)-level-6-(1,2/3)}, 
we have 
\begin{align*}
\frac
{
(q; q)_{\infty}^4 (q^3; q^3)_{\infty}^4
}
{
(q^2; q^2)_{\infty}^2 (q^6; q^6)_{\infty}^2
}
=&
\frac{   (q^2; q^2)_{\infty}   (q^3; q^3)_{\infty}^6  }{   (q; q)_{\infty}^2 (q^6; q^6)_{\infty}^3  }  
\cdot
\frac{(q ;q)_{\infty}^6 (q^6;q^6)_{\infty}  }{(q^2;q^2)_{\infty}^3 (q^3;q^3)_{\infty}^2  } 
=P_1(q)\cdot Q_1(q)  \\
=&\frac16
\left\{
E_2(q)
-4
E_2(q^2)
-
3
E_2(q^3)
+
12
E_2(q^6)
\right\},
\end{align*}
and
\begin{equation*}
f(q)=\frac
{
(q^2; q^2)_{\infty}^4 (q^3; q^3)_{\infty}^8
}
{
(q;q)_{\infty}^8
(q^6; q^6)^4
}
=\frac{P_1(q)}{Q_1(q)}. 
\end{equation*}

\subsection{Derivation of a Riccati equation (2)}

\begin{proposition}
\label{prop:2nd-order-theta(1/3,1)-(2/3,1)}
{\it
For every $\tau\in\mathbb{H}^2,$ we have 
\begin{equation*}
\frac
{
\theta^{\prime \prime}
\left[
\begin{array}{c}
\frac23 \\
1
\end{array}
\right] 
}
{
\theta
\left[
\begin{array}{c}
\frac23 \\
1
\end{array}
\right] 
}
-
\frac
{
\theta^{\prime \prime}
\left[
\begin{array}{c}
\frac13 \\
1
\end{array}
\right] 
}
{
\theta
\left[
\begin{array}{c}
\frac13 \\
1
\end{array}
\right] 
}
=
\frac{1}{12}
\frac
{
\left\{
\theta^{\prime}
\left[
\begin{array}{c}
1 \\
1
\end{array}
\right] 
\right\}^2
}
{
\theta^2
\left[
\begin{array}{c}
0 \\
1
\end{array}
\right] 
\theta^2
\left[
\begin{array}{c}
\frac13 \\
1
\end{array}
\right] 
\theta^6
\left[
\begin{array}{c}
\frac23 \\
1
\end{array}
\right] 
}
\left\{
9
\zeta_3
\theta^{8}
\left[
\begin{array}{c}
\frac23 \\
1
\end{array}
\right] 
-10\zeta_3^2
\theta^4
\left[
\begin{array}{c}
\frac13 \\
1
\end{array}
\right] 
\theta^4
\left[
\begin{array}{c}
\frac23 \\
1
\end{array}
\right] 
+
\theta^{8}
\left[
\begin{array}{c}
\frac13 \\
1
\end{array}
\right] 
\right\}. 
\end{equation*}
}
\end{proposition}

\begin{proof}
The proposition can be obtained by 
considering the derivative formulas in Theorem \ref{thm:analogue-Jacobi-(1,2/3)-(2/3,1)} and Corollary \ref{coro:analogue-Jacobi-(1,1/3)-(1/3,1)}, and 
by comparing the coefficients of the term $z^2$ in equation  (\ref{eqn:two-theta-sum-2-deno-3-(1)}). 
\end{proof}

\begin{theorem}
\label{thm:q-Riccati-diff-eq-modulus6-(2)}
{\it
For $q\in\mathbb{C}$ with $|q|<1,$ set 
\begin{equation*}
g(q)=
q\frac{
(q;q)_{\infty}^4 (q^6; q^6)_{\infty}^8
}
{
(q^2; q^2)_{\infty}^8 (q^3; q^3)_{\infty}^4
},
\,\,q=\exp(2\pi i \tau).
\end{equation*}
Then, $g$ satisfies the following Riccati equation:
\begin{equation}
\label{eqn:Riccati-q-diff-eq-modulus6-(2)}
q
\frac{d}{dq} g
=
q\frac
{
(q^2; q^2)_{\infty}^4 (q^6; q^6)_{\infty}^4
}
{
(q;q)_{\infty}^2 (q^3; q^3)_{\infty}^2
}
\left(
9g^2-10g+1
\right). 
\end{equation}
}
\end{theorem}

\begin{proof}
The heat equation (\ref{eqn:heat}) and 
Proposition \ref{prop:2nd-order-theta(1/3,1)-(2/3,1)}
imply that 
\begin{align*}
&
4\pi i 
\frac
{
d
}
{
d\tau
}
\log 
\theta
\left[
\begin{array}{c}
\frac23 \\
1
\end{array}
\right] 
-
4\pi i 
\frac
{
d
}
{
d\tau
}
\log 
\theta
\left[
\begin{array}{c}
\frac13 \\
1
\end{array}
\right]   \\
=&
\frac{1}{12}
\frac
{
\left\{
\theta^{\prime}
\left[
\begin{array}{c}
1 \\
1
\end{array}
\right] 
\right\}^2
}
{
\theta^2
\left[
\begin{array}{c}
0 \\
1
\end{array}
\right] 
\theta^2
\left[
\begin{array}{c}
\frac13 \\
1
\end{array}
\right] 
\theta^6
\left[
\begin{array}{c}
\frac23 \\
1
\end{array}
\right] 
}
\left\{
9
\zeta_3
\theta^{8}
\left[
\begin{array}{c}
\frac23 \\
1
\end{array}
\right] 
-10\zeta_3^2
\theta^4
\left[
\begin{array}{c}
\frac13 \\
1
\end{array}
\right] 
\theta^4
\left[
\begin{array}{c}
\frac23 \\
1
\end{array}
\right] 
+
\theta^{8}
\left[
\begin{array}{c}
\frac13 \\
1
\end{array}
\right] 
\right\}. 
\end{align*}
which shows that 
\begin{align*}
&
4\pi i 
\frac
{
d
}
{
d \tau
}
\log 
\zeta_3^2
\theta^4
\left[
\begin{array}{c}
\frac23 \\
1
\end{array}
\right] 
-
4\pi i 
\frac
{
d
}
{
d\tau
}
\log 
\theta^4
\left[
\begin{array}{c}
\frac13 \\
1
\end{array}
\right]   \\
=&
\frac{1}{3}
\frac
{
\left\{
\theta^{\prime}
\left[
\begin{array}{c}
1 \\
1
\end{array}
\right] 
\right\}^2
}
{
\theta^2
\left[
\begin{array}{c}
0 \\
1
\end{array}
\right] 
\theta^2
\left[
\begin{array}{c}
\frac13 \\
1
\end{array}
\right] 
\theta^6
\left[
\begin{array}{c}
\frac23 \\
1
\end{array}
\right] 
}
\left\{
9
\zeta_3
\theta^{8}
\left[
\begin{array}{c}
\frac23 \\
1
\end{array}
\right] 
-10\zeta_3^2
\theta^4
\left[
\begin{array}{c}
\frac13 \\
1
\end{array}
\right] 
\theta^4
\left[
\begin{array}{c}
\frac23 \\
1
\end{array}
\right] 
+
\theta^{8}
\left[
\begin{array}{c}
\frac13 \\
1
\end{array}
\right] 
\right\}. 
\end{align*}
\par
Jacobi's triple product identity (\ref{eqn:Jacobi-triple}) yields 
\begin{equation*}
\zeta_3^2
\frac
{
\left\{
\theta^{\prime}
\left[
\begin{array}{c}
1 \\
1
\end{array}
\right] 
\right\}^2
\theta^2
\left[
\begin{array}{c}
\frac13 \\
1
\end{array}
\right] 
}
{
\theta^2
\left[
\begin{array}{c}
0 \\
1
\end{array}
\right] 
\theta^2
\left[
\begin{array}{c}
\frac23 \\
1
\end{array}
\right] 
}
=-4\pi^2
y
\prod_{n=1}^{\infty} \frac{(1-y^{2n})^4 (1-y^{6n})^4 }{(1-y^n)^2 (1-y^{3n})^2     } 
=-4\pi^2\frac{\eta^4(\tau/3) \eta^4(\tau)  }{\eta^2(\tau/6) \eta^2(\tau/2)}, \,\,y=\exp(2 \pi i \tau/6). 
\end{equation*}
\par
Setting 
\begin{equation*}
(X,Y)=
\left(
\zeta_3^2
\theta^4
\left[
\begin{array}{c}
\frac23 \\
1
\end{array}
\right](0,\tau),   
\theta^4
\left[
\begin{array}{c}
\frac13 \\
1
\end{array}
\right](0,\tau) 
\right), 
\end{equation*}
we obtain 
\begin{equation*}
y
\frac{d}{d y}
g(y)
=
y
\frac{d}{d y}
\left(
\frac{X}{Y}
\right)
=
y \frac{ (y^2; y^2)_{\infty}^4 (y^6; y^6)_{\infty}^4 }{ (y; y)_{\infty}^2  (y^3; y^3)_{\infty}^2   }
\left(
9g^2(y)-10 g(y)+1
\right).
\end{equation*}
The theorem can be obtained by changing $\tau\longrightarrow 6\tau.$ 
\end{proof}

\subsubsection*{Remark}
By Theorem \ref{thm:E_2(q)-level-6-(2/3,1)},
we have 
\begin{align*}
q\frac
{
(q^2; q^2)_{\infty}^4 (q^6; q^6)_{\infty}^4
}
{
(q;q)_{\infty}^2 (q^3; q^3)_{\infty}^2
}
=&
\frac{ (q^3; q^3)_{\infty} (q^2; q^2)_{\infty}^6  }{ (q; q)_{\infty}^3 (q^6; q^6)_{\infty}^2 } \cdot
q \frac{(q; q)_{\infty}  (q^6; q^6)_{\infty}^6   }{  (q^2; q^2)_{\infty}^2 (q^3; q^3)_{\infty}^3   }
=P_2(q)\cdot Q_2(q)  \\
=&\frac{1}{24}
\left\{
-
E_2(q)
+
E_2(q^2)
+
3
E_2(q^3)
-3
E_2(q^6)
\right\},
\end{align*}
and
\begin{equation*}
g(q)=q\frac{
(q;q)_{\infty}^4 (q^6; q^6)_{\infty}^8
}
{
(q^2; q^2)_{\infty}^8 (q^3; q^3)_{\infty}^4
}
=\frac{Q_2(q)}{P_2(q)}.
\end{equation*}

\subsection{Derivation of a Riccati equation (3)}

\begin{proposition}
\label{prop:2nd-order-theta(1,1/3)-(0,1/3)}
{\it
For every $\tau\in\mathbb{H}^2,$ we have 
\begin{equation*}
\frac
{
\theta^{\prime \prime}
\left[
\begin{array}{c}
0 \\
\frac13
\end{array}
\right] 
}
{
\theta
\left[
\begin{array}{c}
0 \\
\frac13
\end{array}
\right] 
}
-
\frac
{
\theta^{\prime \prime}
\left[
\begin{array}{c}
1 \\
\frac13
\end{array}
\right] 
}
{
\theta
\left[
\begin{array}{c}
1 \\
\frac13
\end{array}
\right] 
}
=
\frac{1}{12}
\frac
{
\left\{
\theta^{\prime}
\left[
\begin{array}{c}
1 \\
1
\end{array}
\right] 
\right\}^2
}
{
\theta^2
\left[
\begin{array}{c}
0 \\
1
\end{array}
\right] 
\theta^2
\left[
\begin{array}{c}
1 \\
\frac13
\end{array}
\right] 
\theta^6
\left[
\begin{array}{c}
0 \\
\frac13
\end{array}
\right] 
}
\left\{
\theta^{8}
\left[
\begin{array}{c}
1 \\
\frac13
\end{array}
\right] 
-10
\theta^4
\left[
\begin{array}{c}
1 \\
\frac13
\end{array}
\right] 
\theta^4
\left[
\begin{array}{c}
0\\
\frac13
\end{array}
\right] 
+
9
\theta^{8}
\left[
\begin{array}{c}
0 \\
\frac13
\end{array}
\right] 
\right\}. 
\end{equation*}
}
\end{proposition}

\begin{proof}
The proposition can be proved in the same way as Proposition \ref{prop:2nd-order-theta(1,1/3)-(1,2/3)}. 
\end{proof}

\begin{theorem}
\label{thm:q-Riccati-diff-eq-modulus6-(3)}
{\it
For $q\in\mathbb{C}$ with $|q|<1,$ set 
\begin{equation*}
t(q)=
q^{\frac12} 
\frac{
(q^{\frac12} ;q^{\frac12}  )_{\infty}^4   (q^3; q^3)_{\infty}^8
}
{
(q; q)_{\infty}^8 (q^{\frac32}; q^{\frac32})_{\infty}^4
},
\,\,q=\exp(2\pi i \tau).
\end{equation*}
Then, $t=t(q)$ satisfies the following Riccati equation:
\begin{equation}
\label{eqn:Riccati-q-diff-eq-modulus6-(3)}
q
\frac{d}{dq} t
=
\frac12
q^{\frac12}
\frac
{
(q; q)_{\infty}^4 (q^3; q^3)_{\infty}^4
}
{
(q^{\frac12} ;q^{\frac12}  )_{\infty}^2 (q^{\frac32}; q^{\frac32})_{\infty}^2
}
\left(
9t^2-10t+1
\right). 
\end{equation}
}
\end{theorem}

\begin{proof}
The theorem can be proved in the same way as Theorem \ref{thm:q-Riccati-diff-eq-modulus6-(2)}. 
Set
\begin{equation*}
(X,Y)=
\left(
\theta^4
\left[
\begin{array}{c}
1 \\
\frac13
\end{array}
\right], 
\theta^4
\left[
\begin{array}{c}
0 \\
\frac13
\end{array}
\right]  
\right).
\end{equation*}
\end{proof}

\subsubsection*{Remark}
By Theorem \ref{thm:E_2-(1,1/3)-(0,1/3)-(0,2/3)}, 
we have 
\begin{align*}
q^{\frac12}
\frac
{
(q; q)_{\infty}^4 (q^3; q^3)_{\infty}^4
}
{
(q^{\frac12} ;q^{\frac12}  )_{\infty}^2 (q^{\frac32}; q^{\frac32})_{\infty}^2
}
=&
q^{\frac12}
\frac{ (q^{\frac12}; q^{\frac12})_{\infty} (q^3; q^3)_{\infty}^6  }{ (q; q)_{\infty}^2 (q^{\frac32}; q^{\frac32})_{\infty}^3 } 
\cdot
\frac{(q; q)_{\infty}^6  (q^{\frac32}; q^{\frac32})_{\infty}   }{  (q^{\frac12}; q^{\frac12})_{\infty}^3 (q^3; q^3)_{\infty}^2   }
=Q_3(q)\cdot T(q)  \\
=&Q_3(q)\left\{ P_3(q)+3 Q_3(q) \right\}  \\
=&\frac{1}{24}
\left\{
-
E_2(q^{\frac12})
+
E_2(q)
+
3
E_2(q^{\frac32})
-3
E_2(q^3)
\right\},
\end{align*}
and
\begin{equation*}
t(q)=q^{\frac12} 
\frac{
(q^{\frac12} ;q^{\frac12}  )_{\infty}^4   (q^3; q^3)_{\infty}^8
}
{
(q; q)_{\infty}^8 (q^{\frac32}; q^{\frac32})_{\infty}^4
}
=\frac{Q_3(q)}{T(q)}=g(q^{\frac12}).
\end{equation*}

\subsection{Derivation of a Riccati equation (4)}

\begin{proposition}
\label{prop:2nd-order-theta(1,1/3)-(0,2/3)}
{\it
For every $\tau\in\mathbb{H}^2,$ we have 
\begin{equation*}
\frac
{
\theta^{\prime \prime}
\left[
\begin{array}{c}
0 \\
\frac23
\end{array}
\right] 
}
{
\theta
\left[
\begin{array}{c}
0 \\
\frac23
\end{array}
\right] 
}
-
\frac
{
\theta^{\prime \prime}
\left[
\begin{array}{c}
1 \\
\frac13
\end{array}
\right] 
}
{
\theta
\left[
\begin{array}{c}
1 \\
\frac13
\end{array}
\right] 
}
=
\frac{1}{12}
\frac
{
\left\{
\theta^{\prime}
\left[
\begin{array}{c}
1 \\
1
\end{array}
\right] 
\right\}^2
}
{
\theta^2
\left[
\begin{array}{c}
0 \\
0
\end{array}
\right] 
\theta^2
\left[
\begin{array}{c}
1 \\
\frac13
\end{array}
\right] 
\theta^6
\left[
\begin{array}{c}
0 \\
\frac23
\end{array}
\right] 
}
\left\{
\theta^{8}
\left[
\begin{array}{c}
1 \\
\frac13
\end{array}
\right] 
+10
\theta^4
\left[
\begin{array}{c}
1 \\
\frac13
\end{array}
\right] 
\theta^4
\left[
\begin{array}{c}
0\\
\frac23
\end{array}
\right] 
+
9
\theta^{8}
\left[
\begin{array}{c}
0 \\
\frac23
\end{array}
\right] 
\right\}. 
\end{equation*}
}
\end{proposition}

\begin{proof}
The proposition can be proved in the same way as Proposition \ref{prop:2nd-order-theta(1,1/3)-(1,2/3)}. 
\end{proof}

\begin{theorem}
\label{thm:q-Riccati-diff-eq-modulus6-(4)}
{\it
For $q\in\mathbb{C}$ with $|q|<1,$ set 
\begin{equation*}
u(q)=
q^{\frac12} 
\frac{
(q ;q )_{\infty}^4   (q^{\frac32}; q^{\frac32})_{\infty}^4   (q^6; q^6)_{\infty}^4
}
{
(q^{\frac12}; q^{\frac12})_{\infty}^4        (q^2 ;q^2 )_{\infty}^4  (q^3 ;q^3 )_{\infty}^4  
},
\,\,q=\exp(2\pi i \tau).
\end{equation*}
Then, $u=u(q)$ satisfies the following Riccati equation:
\begin{equation}
\label{eqn:Riccati-q-diff-eq-modulus6-(4)}
q
\frac{d}{dq} u
=
\frac12
q^{\frac12}
\frac
{
 (q^{\frac12} ;q^{\frac12}  )_{\infty}^2     (q^{\frac32}; q^{\frac32})_{\infty}^2   (q^2; q^2)_{\infty}^2         (q^6; q^6)_{\infty}^2
}
{
  (q; q)_{\infty}^2    (q^3; q^3)_{\infty}^2
}
\left(
9u^2+10u+1
\right). 
\end{equation}
}
\end{theorem}

\begin{proof}
The theorem can be proved in the same way as Theorem \ref{thm:q-Riccati-diff-eq-modulus6-(2)}. 
Set
\begin{equation*}
(X,Y)=
\left(
\theta^4
\left[
\begin{array}{c}
1 \\
\frac13
\end{array}
\right], 
\theta^4
\left[
\begin{array}{c}
0 \\
\frac23
\end{array}
\right]  
\right).
\end{equation*}
\end{proof}

\subsubsection*{Remark}
By Theorem \ref{thm:E_2-(1,1/3)-(0,1/3)-(0,2/3)}, 
we have 
\begin{align*}
&
q^{\frac12}
\frac
{
 (q^{\frac12} ;q^{\frac12}  )_{\infty}^2     (q^{\frac32}; q^{\frac32})_{\infty}^2   (q^2; q^2)_{\infty}^2         (q^6; q^6)_{\infty}^2
}
{
  (q; q)_{\infty}^2    (q^3; q^3)_{\infty}^2
}  \\
=&
q^{\frac12} \frac{   (q;q)_{\infty}   (q^{\frac32}; q^{\frac32})_{\infty}^3    (q^6;q^6)_{\infty}^3    }{   (q^{\frac12}; q^{\frac12})_{\infty}   (q^2;q^2)_{\infty}    (q^3;q^3)_{\infty}^3 }
\cdot
\frac{(q^{\frac12} ; q^{\frac12} )_{\infty}^3  (q^2; q^2)_{\infty}^3   (q^3; q^3)_{\infty}   }{  (q; q)_{\infty}^3 (q^{\frac32}; q^{\frac32})_{\infty}  (q^6; q^6)_{\infty}   }
=R_3(q)\cdot U(q)  \\
=&R_3(q)\left\{ P_3(q)-3 R_3(q) \right\}  \\
=&\frac{1}{24}
\left\{
-
E_2(q^{\frac12})
+
5
E_2(q)
-
4
E_2(q^2)
+
3
E_2(q^{\frac32})
-15
E_2(q^3)
+
12
E_2(q^6)
\right\},
\end{align*}
and
\begin{equation*}
u(q)=q^{\frac12} 
\frac{
(q ;q )_{\infty}^4   (q^{\frac32}; q^{\frac32})_{\infty}^4   (q^6; q^6)_{\infty}^4
}
{
(q^{\frac12}; q^{\frac12})_{\infty}^4        (q^2 ;q^2 )_{\infty}^4  (q^3 ;q^3 )_{\infty}^4  
}
=\frac{R_3(q)}{U(q)}.
\end{equation*}




\end{document}